\documentclass{lms}

\usepackage{latexsym}

\usepackage[utf8]{inputenc} 
\usepackage{hyperref}       
\usepackage{url}            
\usepackage{booktabs}       
\usepackage{amsmath,amsfonts,amssymb,bm}
\usepackage[small,nohug]{diagrams}
\usepackage{nicefrac}       
\usepackage{microtype}      

\newcommand{\Z}{\mathbb{Z}}

\newcommand{\R}{\mathbb{R}}
\newcommand{\C}{\mathbb{C}}

\newcommand{\OC}{\mathcal{O}}
\newcommand{\PR}{\mathbb{P}}
\newcommand{\iv}[1]{\iota_{#1}}
\newcommand{\poltype}{\underline{d}} 
\newcommand{\AbelM}{\mathcal{A}} 
\newcommand{\Cu}{\mathcal{C}} 
\newcommand{\Du}{\mathcal{D}} 
\newcommand{\Surf}{\mathcal{S}}  
\newcommand{\Jac}{\mathcal{J}} 

\newcommand{\TD}{\mathcal{T}} 

\newcommand{\BITAN}{\mathbb{T}_1} 
\newcommand{\TAN}{\mathbb{T}_2} 

\newcommand{\IMh}[1]{\Im m \ #1} 
\newcommand{\RC}{\kappa} 
\newcommand{\WEIL}{W} 
\newcommand{\GAU}{\mathit{G}} 
\newcommand{\GrpJ}{\mathcal{G}} 

\newcommand{\NDiv}{\mathcal{D}}  
\newcommand{\GRAD}{\nabla} 

\newcommand{\CAN}[1]{\phi_{#1}} 

\newcommand{\ETA}[1]{\eta_{#1}} 
\newcommand{\OMEGA}[1]{\omega_{#1}} 
\newcommand{\PUNKT}{\ \ \text{.}}  
\newcommand{\BEISTRICH}{\ \ \text{,}}  
\newcommand{\SET}[2]{\{#1 \ | \ #2 \}}  
\newcommand{\B}{B}
\newcommand{\E}{E}
\newcommand{\tauB}{\tau_{\B}}
\newcommand{\tauE}{\tau_{\E}}
\newcommand{\PN}{\mathcal{N}}      
\newcommand{\SA}{S}                
\newcommand{\prBSA}{{\pi_{\B}}}          
\newcommand{\PMp}{{\mathcal{M}^{+}}}  
\newcommand{\PMm}{{\mathcal{M}^{-}}}  
\newcommand{\PMpij}[1]{{\mathcal{M}_{#1}^{+}}}  
\newcommand{\PMmij}[1]{{\mathcal{M}_{#1}^{-}}}  

\newcommand{\prTA}{p}           
\newcommand{\ratA}{{\psi}}      
\newcommand{\cf}{\nu}           
\newcommand{\cg}{\xi}           
\newcommand{\cfc}{\Cu}          
\newcommand{\cgc}{\Du}          

\newcommand{\inv}{{\iota}}      

\newcommand{\EU}[1]{\mathcal{U}_{#1}} 
\newcommand{\EV}[1]{\mathcal{V}_{#1}}
\newcommand{\x}[1]{x_{#1}}
\newcommand{\y}[1]{y_{#1}}
\newcommand{\Eeq}[1]{g_{#1}}             
\newcommand{\Eeqinf}[1]{h_{#1}}             
\newcommand{\vinf}[1]{v_{#1}}
\newcommand{\winf}[1]{w_{#1}}
\newcommand{\inftyplus}{\infty_{+}}
\newcommand{\inftyminus}{\infty_{-}}
\newcommand{\pa}[1]{\delta_{#1}}
\newcommand{\LF}{\mathcal{L}} 
\newcommand{\cb}{b_{23}}
\newcommand{\cc}{b_{13}}
\newcommand{\cd}{b_{12}}

\newcommand{\BBF}[1]{\bm{#1}}

\newcommand{\DER}[2]{\frac{\partial #2}{\partial #1}}   

\newcommand{\SKMULT}[1]{\DER{#1}{}}
\newcommand{\DERP}[1]{\overset{\text{\tiny$\bullet$}}{#1}}     
\newcommand{\DTH}[1]{\DERP{\theta_{#1}}}         

\newcommand{\SiegH}{{\mathcal{H}}}
\newcommand{\HD}{{\mathcal{H}_{3\Delta, \tau_0}}} 
\newcommand{\DEFP}{\mathcal{P}} 

\newcommand{\secV}{v} 
\newcommand{\xc}{x}
\newcommand{\yc}{y}
\newcommand{\zc}{z}
\newcommand{\wc}{w}

\newcommand{\rhos}[2]{\rho_{#1}} 

\def\daschmapsto{\mathrel{\mapstochar\dashrightarrow}} 

\newcommand{\leftexp}[2]{{\vphantom{#2}}^{#1}{#2}}
\newcommand{\traspose}[1]{\leftexp{t}{#1}}  
\newcommand{\bigslant}[2]{{\raisebox{.2em}{$#1$}\left/\raisebox{-.2em}{$#2$}\right.}}  
\renewcommand\labelitemi{$\bullet$}
\newcommand{\noindentTX}{\noindent} 

\newtheorem{thm1}{theorem}[section]
\newtheorem{thm1intro}{theorem} 

\newtheorem{lem1}[thm1]{Lemma}

\newtheorem{prop1}[thm1]{Proposition}

\newnumbered{que1}[thm1]{Question}			
\newnumbered{que1intro}[thm1intro]{Question}

\newnumbered{def1}[thm1]{Definition}
\newnumbered{ex1}[thm1]{Example}
\newnumbered{not1}[thm1]{Notation}
\newnumbered{oss1}[thm1]{Remark}
\newnumbered{rmk}[thm1]{Remark}
\newnumbered{note}[thm1]{Note}
\numberwithin{equation}{thm1}
\newnumbered{sketch}[thm1]{On the proof of theorem \ref{teoremafinale} - Sketch}      

\newnumbered{assertion}{Assertion}    
\newnumbered{conjecture}{Conjecture}  
\newnumbered{definition}{Definition}
\newnumbered{hypothesis}{Hypothesis}
\newnumbered{remark}{Remark}
\newnumbered{observation}{Observation}
\newnumbered{problem}{Problem}
\newnumbered{question}{Question}
\newnumbered{algorithm}{Algorithm}
\newnumbered{example}{Example}
\newnumbered{notation}{Notation} 

\title[Very ampl. can. bundle of surf. $(1,2,2)$ on ab. $3$-folds]
 {Very ampleness of the canonical bundle of surfaces of type $(1,2,2)$ on abelian threefolds} 

\author{L. Cesarano}

\classno{14J29, 14J70, 14K25, 14K20}

\extraline{The present work was supported by the ERC Advanced grant n. 340258, 'TADMICANT'.}

\begin{document}
\maketitle

\begin{abstract}
	The present work deals with the canonical map of smooth, compact complex surfaces of general type in a 
	polarization of type $(1,2,2)$ on an abelian threefold. A natural and classical question is whether the canonical system of such surfaces is very ample in the general case. In this work, we provide a positive answer 
	to this question. First, we describe the structure of the canonical map of those smooth ample surfaces of type $(1,2,2)$ in an abelian threefold which are bidouble cover of principal polarizations. Then, we study the general behavior of the canonical map of general ample surfaces $\Surf$, yielding a $(1,2,2)$-polarization on an abelian threefold $A$ which is isogenous to a product.  
	By combining these descriptions, we show that the canonical map yields a holomorphic embedding when $A$ and $\Surf$ are both sufficiently general. 
	It follows, in particular, a proof of the existence of canonical irregular surfaces in $\PR^5$ with numerical 
	invariants $p_g = 6$, $q = 3$ and $K^2 = 24$.
\end{abstract}



\section*{Introduction}

The present work aims at studying the canonical map of smooth, compact complex surfaces of general type $\Surf$ in a polarization of type $(1, 2, 2)$ on an abelian threefold $A$. A classical question to ask is whether the canonical system of such surfaces is very ample, at least when $A$ and $\Surf$ are sufficiently general. The following theorem, which is the main result of this work, provides an affirmative answer to this question.
\begin{thm1intro} \label{teoremafinaleintro}
	Let be $(A, \mathcal{L})$ a general $(1,2,2)$-polarized abelian threefold and let be $\Surf$ a general surface in the linear system $|\mathcal{L}|$. Then the canonical system of $\Surf$ yields a holomorphic embedding in $\PR^5$. 
\end{thm1intro} 

 A first motivation underlying this research question is the well-known existence problem of canonically embedded surfaces in $\PR^5$ with prescribed numerical invariants:
\begin{que1intro} \label{questionE1}
Which are the possible values of $K^2$ for the smooth surfaces of general type $\Surf$ which are canonically embedded in $\PR^5$? 
\end{que1intro}
\noindentTX It's important to recall that the requirement for the canonical map to be birational onto its image leads to some numerical constraints. From the inequality of Bogomolov-Miyaoka-Yau $K^2 \leq 9 \chi$ and the Debarre version of Castelnuovo's inequality $K^2 \geq 3p_g(\Surf) + q - 7$ (see \cite{Debarre1982}), it follows that, if the canonical map of an algebraic surface $\Surf$ with $p_g = 6$ is birational, then 
\begin{equation*} 
11 + q \leq K^2 \leq 9(7 - q) \PUNKT
\end{equation*}

Several construction methods have been considered in the attempt to give a satisfactory answer to Question \ref{questionE1} at least for some values of $K^2$. One first method we want to mention here relies on the 
existence of a very special determinantal structure for the defining equations of a canonical surface in $\PR^5$. More precisely, 
Walter proved (see \cite{Walter}) that every codimension $3$ locally Gorenstein subscheme $X$ of $\PR^{n+3}$, which is $l$-subcanonical 
(i.e. $\omega_X \cong \mathcal{O}_X(l))$ and such that $\chi(\mathcal{O}_X(l/2))$ is even if $n$ is divisible by $4$ and $l$ is even, is pfaffian.
This means that there exists a locally free resolution 
\begin{equation} \label{pfaffiand}
0 \longrightarrow \OC_{\PR^5}(-t) \longrightarrow \mathcal{E}^{\vee}(-t)\overset{\alpha}{\longrightarrow} \mathcal{E} \longrightarrow \mathcal{I}_{\Surf}\longrightarrow  0 \BEISTRICH
\end{equation}
where $\mathcal{E}$ is a vector bundle on $\PR^5$ of odd rank $2k+1$, and $\alpha$ is an antisymmetric map such that $\Surf$ is defined by the Pfaffians of order $2k$ of $\alpha$.
This implies that every smooth canonical surface $\Surf$ in $\PR^5$ is pfaffian with a resolution as in (\ref{pfaffiand}) with $t = 7$.
Moreover, in this case, the cohomology of the vector bundle $\mathcal{E}$ (and hence $\mathcal{E}$ itself by using the Horrock's correspondence) is completely determined by the cohomology of the structure sheaf of $\Surf$ and its ideal sheaf in $\PR^5$, and hence by the structure of the canonical ring of $\Surf$ (see also \cite{Catanese1997})

\noindentTX By considering a vector bundle $\mathcal{E}$ which splits into a sum of line bundles and a suitable antisymmetric map of vector bundles $\alpha$ as in Diagram $\ref{pfaffiand}$, Catanese exhibited examples of regular canonical surfaces in $\PR^5$ with $11 \leq K^2 \leq 17$ (see \cite{Catanese1997}). Moreover, G. and M. Kapustka construct via bilinkage methods a further example of canonical regular surface $\Surf$ in $\PR^5$ with $K^2=18$ (see \cite{GMKapustka2016}). 

The same two authors conjectured in a previous work (cf. \cite{GMKapustka2015}) that this degree would be an upper bound for the existence problem \ref{questionE1}. However, this conjecture has turned out to be false:
 Catanese showed (see \cite{Catanese1999}) that examples of surfaces of general type with birational canonical map of higher degree arise by considering ramified bidouble covers with prescribed branch type, and he proved that there is a family of surfaces of general type which are bidouble covers of $\PR^1 \times \PR^1$, branched on three smooth curves $D_1$, $D_2$, $D_3$ of bidegree respectively $(2,3)$, $(2,3)$, $(1,4)$, whose general fibers are surfaces with $K^2 = 24$, $q=0$, $p_g = 6$ and birational canonical map. Furthermore, in a recent paper (see \cite{Catanese2017}) he proved that, in this case, the canonical system is base point free and yields an embedding in $\PR^5$. 

\noindentTX However, it is remarkable that the given examples in \cite{Catanese1997,Catanese1999,Catanese2017,GMKapustka2015} are examples of regular surfaces, and thus they leave unanswered the existence problem \ref{questionE1} for canonical irregular surfaces of higher degree.
 
\noindentTX Examples of irregular algebraic projective varieties can be easily produced using transcendental methods. Nice examples of irregular projective varieties arise naturally by considering abelian varieties or subvarieties of abelian varieties. However, in the case of abelian varieties, the description of the equations of their projective model and the problem of determining whether an ample line bundle is very ample have been considered challenging problems 
(see \cite{MumfordEquation} and \cite{DebarreTor}), although their underlying analytic structure, as well as the structure of their ample line bundles, are well understood.

A first important class of subvarieties of abelian varieties is that of their ample divisors. Since the canonical sheaf of an ample divisor $\mathcal{D}$ is an abelian variety $A$ is, by adjunction, just the restriction of the polarization $\mathcal{O}_A(\mathcal{D})$ to $\mathcal{D}$, some natural questions arise:


\begin{que1intro}\label{questionPolarizationE}
Let us consider a couple $(A, \mathcal{D})$, where $A$ a $g$-dimensional $\poltype := (d_1, \cdots, d_g)$ abelian variety, and $\mathcal{D}$ smooth ample divisor of type $\poltype$. For which polarization types $\poltype$ the canonical map is birational, at least for the general couple $(A, \mathcal{D})$?
\end{que1intro}

\begin{que1intro}\label{questionPolarizationF}
In the notations of Question \ref{questionPolarizationE}, for which polarization types $\poltype$ the canonical map is a holomorphic embedding for the general couple $(A, \mathcal{D})$?
\end{que1intro}

The results contained in a recent work (see \cite{CatCes2021}) answer to Question \ref{questionPolarizationE}: the canonical map is birational, provided that the polarization given by $\mathcal{D}$ is non-principal (i.e $d_g \neq 1$) and $\mathcal{D}$ is sufficiently general within its polarization class. \newline
Concerning Question \ref{questionPolarizationF} above, in the same work, it is proven that if the canonical bundle of a smooth ample divisor $\mathcal{D}$ is very ample, then the following inequality holds
\begin{equation*} \label{NECCOND}
d:=d_1 \cdot d_2 \cdot \cdots \cdot d_g \geq g +1 \PUNKT
\end{equation*}

In consideration of the results which we present in this work, Catanese has conjectured that the necessary condition above for the very ampleness of the canonical bundle of a smooth ample divisor $\mathcal{D}$ is also sufficient, provided that the couple $(A, \mathcal{D})$ is sufficiently general.

The present paper provides the first known result in the direction of this conjecture, at least in the case of abelian threefolds when $d=4$. However, it leaves the case of a polarization of type $(1,1,4)$ open. In this case, whether the canonical system yields an embedding in $\PR^5$ in the general case remains unknown. 
Nevertheless, in the case of the polarization type $(1,2,2)$, we provide very explicit examples and descriptions of the canonical map of smooth surfaces within a $(1,2,2)$-polarization class which are interesting for their own sake. \newline

This work is organized as follows.

\noindentTX For the sake of completeness, we focus throughout Section 1 on some theoretical background about the structure of the canonical map of an ample divisor on an abelian variety, and we give a proof that this map can be expressed in terms of Theta functions (see Proposition \ref{desccanonicalmap}).\newline 

\noindent Section 2 treats the case of a smooth surface, which is a bidouble cover of a principal polarization of a general Jacobian threefold. The study of the canonical map of surfaces in a polarization of type $(1,1,2)$, which has been carried out by Catanese and Schreier in their joint
work \cite{CataneseSchreyer1}, plays here a crucial role. It turns out that the canonical map of a surface $\Surf$, which is the pullback of a principal polarization of a general Jacobian threefold along an isogeny $p$ with kernel isomorphic to $\Z_2^2$, is a local holomorphic embedding. However, its canonical map is never injective. Indeed, we show that, on the invariant canonical curves under the action of an involution contained in the kernel of the isogeny $p$, the canonical map has degree $2$ and it factors through the involution itself. 
\noindentTX The (apparently negative) outcome of this first attempt to prove Theorem \ref{teoremafinaleintro} suggests that the canonical map may fail to be injective only for those surfaces $\Surf$, which are invariant under the action of an involution on $A$. In such a case, on an invariant canonical curve in $\Surf$, the restriction of the canonical map should factor through the involution itself, and in particular it is not injective. \newline

\noindentTX Section 3 deals with the properties of the morphisms defined by ample line bundles of type $(1,2)$ and $(2,4)$ on abelian surfaces, and treats some results, which we apply to prove Theorem \ref{teoremafinaleintro}. 

\noindentTX In Section 4, we deal with general ample surfaces in general $(1,2,2)$-abelian threefolds which are isogenous to the product of elliptic curve. 
Such abelian threefolds $A$ admit changing-sign involutions, and under their action all smooth surfaces $\Surf$ within the polarization on $A$ are invariant. In particular, for every such involution, we see that the canonical map factors through it precisely on some invariant canonical curves in $\Surf$. \newline

The discussion of the above cases in which $\Surf$ is a bidouble cover of a principal polarization or $A$ is isogenous to a product constitute the background for the proof of Theorem \ref{teoremafinaleintro}, which is contained in Section 5, the last of this work. The proof of it proceeds by assuming that the canonical map of a general member $\Surf$ within the $(1,2,2)$-polarization class in a general abelian threefold $A$ always fails to be injective on a couple of points which varies when $A$ moves in the moduli space of $(1,2,2)$-polarized abelian threefolds. By specializing to the different loci of abelian threefolds which are isogenous to a product, we find that the original couple of points also specialize to a couple of points which are conjugated with respect to many different involutions. This argument will show that the original assumption leads to a contradiction, and hence it will follow that the canonical map is, in general, injective and with injective differential everywhere, according to the results of the second section.

\setcounter{section}{0}
\section{Preliminaries and Notations} 
 Throughout this work, a polarized abelian variety will be a couple $(A, \mathcal{L})$, where $\mathcal{L}$ is an ample line bundle on a complex torus $A$ and we denote by $|\mathcal{L}|$ the linear system of effective ample divisors which are zero loci of global holomorphic sections of $\mathcal{L}$. \newline

\noindentTX The first Chern Class $c_1(\mathcal{L}) \in H^2(A, \Z)$ is an integral valued alternating bilinear form on the lattice $H_1(A, \Z)$. Applying the elementary divisors theorem, we obtain that there exists a basis $\lambda_1, \cdots, \lambda_g, \mu_1, \cdots, \mu_g$ of $\Lambda$ with respect to which the matrix of $c_1(\mathcal{L})$ is
\begin{equation*}
\begin{pmatrix}
0 & D \\
-D & 0
\end{pmatrix} \BEISTRICH
\end{equation*}
where $D$ is a diagonal matrix $diag(d_1, \cdots, d_g)$ of positive integers with the property that every integer in the sequence divides the next. We call the sequence of integers $\poltype = (d_1, \cdots, d_g)$ the type of the polarization on $A$ induced by $\mathcal{L}$. 

\noindentTX Moreover, we will say that an ample divisor $\mathcal{D}$ in an abelian variety $A$ yields a polarization of type $\poltype = (d_1, \cdots, d_g)$ on $A$, or just that $\mathcal{D}$ is of type $\poltype$ on $A$ if the type of the polarization of $\mathcal{L}(\mathcal{D})$ is $\poltype$. \newline

	Let us consider a polarized abelian variety $(A, \mathcal{L})$ of dimension $g$, where $A := \bigslant{\C^g}{\Lambda}$, and $\Lambda$ denote a sublattice in $\C^g$. Consider $\NDiv$ a smooth ample divisor in the linear system $|\mathcal{L}|$. Denoted by $[\{\phi_{\lambda}\}_{\lambda}] \in H^1(\Lambda, H^0(\mathcal{O}_V^*)$ the factor of automorphy corresponding to the ample line bundle $\mathcal{L}$ according to the Appell-Humbert theorem 
	(see \cite[p. 32]{BLange}), the vector space $H^0(A, \mathcal{L})$ is isomorphic to the vector space of the holomorphic functions $\theta$ on $\C^g$ which satisfy, for every $\lambda$ in the lattice $\Lambda$, the functional equation
\begin{equation*}    
\theta(z + \lambda) =  \phi_{\lambda}(z) \theta(z) \PUNKT
\end{equation*}

\noindentTX When $\NDiv = div(\theta_0)$ is a divisor which is the zero locus of a holomorphic global section $\theta_0$ of $\mathcal{L}$, we have, by the adjunction formula, that
\begin{equation} \label{canbdl}
\omega_{\NDiv} = (\mathcal{O}_A(\NDiv) \otimes \omega_A)|_{\NDiv} = \mathcal{O}_{\NDiv}(\NDiv) \PUNKT
\end{equation}
The derivatives $\frac{\partial \theta_0}{\partial z_j}$ are global holomorphic sections of $\mathcal{O}_{\NDiv}(\NDiv)$ for every $j = 1, \cdots, g$. Indeed, for every $\lambda$ in $\Lambda$ and for every $z$ on $\NDiv$ we have
\begin{equation*}
\frac{\partial \theta_0}{\partial z_j}(z+\lambda) = \phi_{\lambda}(z) \frac{\partial \theta_0}{\partial z_j}(z) + \frac{\partial \phi_{\lambda}}{\partial z_j}(z) \theta_0(z) = \phi_{\lambda}(z) \frac{\partial \theta_0}{\partial z_j}(z) \PUNKT
\end{equation*}
This leads naturally to a description of the canonical map of a smooth ample divisor in an abelian variety only in terms of theta functions.

\begin{prop1} \label{desccanonicalmap}
Let $A=\bigslant{\C^g}{\Lambda}$ be an abelian variety and $\mathcal{L}$ an ample line bundle. Let $\NDiv$ be a smooth divisor defined as the zero locus of a holomorphic section $\theta_0$ of $\mathcal{L}$. Moreover, let us suppose $\theta_0, \cdots, \theta_n$ is a basis for the vector space $H^0(A, \mathcal{L})$. Then $\theta_1, \cdots, \theta_n, \frac{\partial \theta_0}{\partial z_1}, \cdots \frac{\partial \theta_0}{\partial z_g}$, where $z_1, \cdots, z_g$ are the flat uniformizing coordinates of  $\C^g$, is a basis for $H^0(\NDiv, \omega_{\NDiv})$.
\end{prop1}
\begin{proof}
From now on let us consider $\mathcal{L}$ to be $\mathcal{O}_A(\NDiv)$. We observe first that, for instance, by the Kodaira vanishing theorem, all the cohomology groups $H^i(A, \mathcal{O}_A(\NDiv))$ vanish, so this implies that $H^0(\NDiv, \mathcal{O}_{\NDiv}(\NDiv))$ has the desired dimension $n+g$. In order to prove the assertion of Proposition \ref{desccanonicalmap}, it is then enough to prove that the connecting homomorphism $\delta^0\colon H^0(\NDiv, \OC_A(\NDiv)) \longrightarrow H^1(A, \OC_A)$ maps $\frac{\partial \theta_0}{\partial z_1}, \cdots \frac{\partial \theta_0}{\partial z_g}$ to $g$ linearly independent elements. \newline

\noindentTX Let us denote the projection of $\C^g$ onto $A$ by $\pi$ and the divisor $\pi^*(\NDiv)$ by $\widehat{\NDiv}$. We have then the short exact sequence
\begin{equation*}
0 \longrightarrow \OC_{\C^g} \longrightarrow \OC_{\C^g}(\widehat{\NDiv}) \longrightarrow \OC_{\widehat{\NDiv}}(\widehat{\NDiv}) \longrightarrow 0 \BEISTRICH
\end{equation*} 
and we can denote the respective cohomology groups by:
\begin{equation} \label{modulesMNP}
\begin{split}
M&:=H^0(\C^g, \OC_{\C^g}) =  H^0(\C^g, \pi^*\OC_{A})\\
N&:=H^0(\C^g, \OC_{\C^g}(\widehat{\NDiv})) =   H^0(\C^g, \pi^*\OC_{A}(\NDiv)) \\
P&:= H^0(\widehat{\NDiv}, \OC_{\widehat{\NDiv}}(\widehat{\NDiv})) =   H^0(\widehat{\NDiv}, \pi^*\OC_{\NDiv}(\NDiv)) \PUNKT 
\end{split}
\end{equation}
The cohomology groups we have just defined in (\ref{modulesMNP}) are $\Lambda$-modules with respect to the following actions: for every element $\lambda$ of $\Lambda$ and every elements $s$, $t$, $u$ respectively in $M$, $N$, and $P$, the action of $\Lambda$ is defined as follows:
\begin{equation} \label{definitionlambdaaction}
\begin{split}
\lambda.s(z) &:= s(z+\lambda) \\
\lambda.t(z) &:= t(z+\lambda)\phi_{\lambda}(z)^{-1} \\
\lambda.u(z) &:= u(z+\lambda)\phi_{\lambda}(z)^{-1} \PUNKT
\end{split}
\end{equation}
According to Mumford \cite[Appendix 2]{MumfordAbelian}, there exists a natural homomorphism $\psi_{\bullet}$ from the cohomology groups sequences $H^p(\Lambda, \cdot)$ and $H^p(A, \cdot)$:
 \begin{diagram} 
		\cdots & \rTo &H^i(\Lambda, M)	&\rTo			&H^i(\Lambda, N) &   \rTo			&H^i(\Lambda, P) & \rTo & H^{i+1}(\Lambda, M) & \rTo & \cdots	\\ 	
		 &&\dTo_{\psi^M_i} &							&\dTo_{\psi^N_i} 			&			& \dTo_{\psi^P_i}  		&			 & 	\dTo_{\psi^M_{i+1}}  &	\\
		\cdots &\rTo &H^i(A, \OC_A)				&\rTo			& H^i(A, \OC_A(\NDiv))&   \rTo			&H^i(\NDiv, \OC_{\NDiv}(\NDiv)) & \rTo^{\delta^i} & H^{i+1}(A, \OC_A) & \rTo & \cdots	\\ 	
 \end{diagram}
The homomorphism $\psi_{\bullet}$ is actually an isomorphism (this means that all the vertical arrows are isomorphisms), because the cohomology groups $H^i(\C^g,  \OC_{\C^g}(\widehat{\NDiv}))$ vanish for every $i>0$, being $\C^g$ a Stein manifold. 
Another possible method to prove that the cohomology sequences $H^p(\Lambda, \cdot)$ and $H^p(A, \cdot)$ are isomorphic, is to use the following result: if $X$ is a variety, $G$ is a group acting freely on $X$ and $\mathcal{F}$ is a $G$-linearized sheaf, then there is a spectral sequence with $E_1$ term equal to $H^p(G, H^q(X, \mathcal{F}))$ converging to $H^{p+q}(Y, \mathcal{F})^{G}$.

\noindentTX The natural identification of these cohomology group sequences allows us to compute $\delta^0\bigl(\frac{\partial \theta_0}{\partial z_j} \bigr)$ using the explicit expression of the connecting homomorphism $H^0(\Lambda, P) \longrightarrow H^1(\Lambda, M)$: given an element $s$ of $P^{\Lambda}$, there exists an element $t$ in $N$ such that $t|_{\widehat{\NDiv}} = s$. Then, by definition of $d\colon N \longrightarrow \mathcal{C}^{1}(\Lambda; N)$, we have
\begin{equation*}
(dt)_{\lambda} = \lambda.t - t \PUNKT 
\end{equation*}
where $\lambda.t$ is defined according to (\ref{definitionlambdaaction}).
Now, from the invariance of $s$ under the action of $\Lambda$, we get
\begin{equation*}
 (\lambda.t - t)|_{\widehat{\NDiv}} = \lambda.s - s = 0 \PUNKT 
\end{equation*}
Hence, for every $\lambda$ there exists a constant $c_{\lambda} \in \C$ such that $\lambda.t - t = c_{\lambda} \theta_0$, and it follows that, by definition,
\begin{equation} \label{HOMCONN}
\delta^0(s)_{\lambda} = c_{\lambda} = \frac{\lambda.t - t}{\theta_0} \PUNKT 
\end{equation}
If we apply now (\ref{HOMCONN}) to the elements $\frac{\partial \theta_0}{\partial z_j}$, we obtain 
\begin{equation*}
\delta^0\left(\frac{\partial \theta_0}{\partial z_j} \right) = \left[ \left( \pi H(e_j, \lambda) \right)_{\lambda} \right] \in H^1(\Lambda; M) \BEISTRICH 
\end{equation*}
where $H$ is the positive definite hermitian form on $\C^g$ which corresponds to the ample line bundle $\mathcal{O}_A(\NDiv)$ by applying the Appell-Humbert theorem. We prove now that these images are linearly independent in $H^1(\Lambda; M)$. Assume that we have coefficients $a_1, \cdots, a_g \in \C$ such that:
\begin{equation*}  
\left[\left(a_1 H(e_1, \lambda) + \cdots + a_g H(e_g, \lambda)\right)_{\lambda}\right] = 0	\PUNKT 
\end{equation*}
This means that there exists $f \in C^0(\Lambda, M)$ such that, for every $\lambda \in \Lambda$, we have:
\begin{equation*}
a_1 H(e_1, \lambda) + \cdots + a_g H(e_g, \lambda) = \lambda.f(z) - f(z) = f(z+\lambda) - f(z)	\PUNKT 
\end{equation*}
For such $f$, the differential $df$ is a holomorphic $\Lambda$-invariant $1$-form. Hence, for some complex constant $c$ and a certain $\C$-linear form $L$, we can write
\begin{equation*}
f(z) = L(z) + c \PUNKT 
\end{equation*}
Hence, for every $\lambda \in \Lambda$ the following holds
\begin{equation*}
H(a_1e_1 + \cdots + a_ge_g, \lambda) = L(\lambda) \BEISTRICH
\end{equation*}
so the same holds for every $z \in \C^g$. This allows us to conclude that $L=0$, $L$ being both complex linear and complex antilinear. But the form $H$ is non-degenerate, so we conclude that $a_1 = \cdots = a_g = 0$. The proposition is proved.
\end{proof}

We can now easily compute the invariants of a smooth ample divisor $\NDiv$ on an abelian variety $A$ of dimension $g$. We will denote throughout this work by:
\begin{itemize}
	\item $p_g := h^0(\NDiv, \omega_{\NDiv})$ the geometric genus of $\NDiv$,
	\item  $q := h^1(\NDiv, \mathcal{O}_{\NDiv})$ the irregularity $\NDiv$,
	\item  $K_{\NDiv}$ a canonical divisor on $\NDiv$. \newline
\end{itemize}

\noindentTX The following formulas in Proposition \ref{invariantdivisors} can be easily verified by applying the Kodaira vanishing theorem and Proposition \ref{desccanonicalmap}.
\begin{prop1} \label{invariantdivisors}
Let $\NDiv$ a smooth divisor in a polarization of type $(d_1, \cdots ,d_g)$ on an abelian variety $A$. Then the invariants of $\NDiv$ are the following:
\begin{align*}
p_g &=  \prod_{j=1}^g d_j + g - 1 \\
q &= g \\
K_{\NDiv}^{g-1} &= g! \prod_{j=1}^g d_j \PUNKT 
\end{align*}
\end{prop1}

\begin{def1} (The Gauss map)
\noindentTX Even though the canonical map of a smooth divisor $\NDiv$ in the linear system $|\mathcal{L}|$ of a (polarized) abelian variety $(A, \mathcal{L})$ can be explicitly expressed in terms of theta functions (as we saw in Proposition \ref{desccanonicalmap}), its image is not always easy to describe. 
If we allow $\NDiv$ to be any divisor, not necessarily reduced and irreducible, the same Proposition \ref{desccanonicalmap} provides for us a basis for the space of holomorphic sections of $\mathcal{L}|_{\NDiv}$.  It makes sense, however, to consider the map defined as follows:
\begin{align*}
\GAU : \NDiv &\dashrightarrow  \PR (V)^{\vee} \\
x &\mapsto \PR (T_x \NDiv) \PUNKT
\end{align*}
This map is called the \textbf{Gauss Map}, and it is clearly defined on the smooth part of the support of $\NDiv$. In particular, if $\NDiv$ is defined as the zero locus of a holomorphic non-zero section $\theta \in H^0(A, \mathcal{L})$, the map $\GAU: \NDiv \dashrightarrow \PR (V)^{\vee} \cong \PR^{g-1}$ is defined by the linear subsystem of $|\mathcal{L}|_{\NDiv}|$ generated by $\frac{\partial \theta}{\partial z_1}, \cdots, \frac{\partial \theta}{\partial z_g}$. 
\end{def1}
\begin{ex1} \label{Gaussjacob}
\noindentTX Let us consider the well-known case of a principal polarization $\Theta$ of the Jacobian $\Jac$ of a smooth curve $\Cu$ of genus $g$. In this case, the Gauss Map coincides with the map defined by the complete linear system $|\mathcal{\Jac}(\Theta)|_{\Theta}|$ and it can be geometrically described as follows: the Abel-Jacobi theorem induces an isomorphism $\Jac \cong Pic^{g-1}(\Cu)$, so $\Theta$ can be viewed, after a suitable translation, as a divisor of $Pic^{g-1}(\Cu)$. The Riemann Singularity Theorem states 
(see \cite[Chapter VI]{ACGH}):
\begin{equation*}
mult_L \Theta = h^0(\Cu, L) \PUNKT
\end{equation*}
By a geometric interpretation of the Riemann-Roch theorem for algebraic curves, it follows that a point $L$ on the Theta divisor represented by the divisor $D = \sum_{j=1}^g{P_j}$ is smooth precisely when the linear span $\left<\phi(P_1), \cdots, \phi(P_g)\right>$ in $\PR H^0(\Cu,\mathcal{\omega_{\Cu}})$ is a hyperplane, where $\phi: \Cu \longrightarrow \PR H^0(\Cu,\mathcal{\omega_{\Cu}})^{\vee}$ denotes the canonical map of $\Cu$. Viewing now $\Jac$ as the quotient of $H^0(\Cu, \omega_{\Cu}^{\vee})$ by the lattice $H_1(\Cu, \Z)$, the Gauss Map associates to $L$ the tangent space $\PR(T_{L} \Theta)$, which is a hyperplane of $\PR(T_{L} \Jac) =  \PR H^0(\Cu, \omega_{\Cu})^{\vee}$ defined as follows:
\begin{align*} 
\GAU \colon &\Theta \dashrightarrow \PR H^0(\Cu,\mathcal{\omega_{\Cu}})^{\vee} \\
\sum_{j=1}^g{P_j}& \mapsto \left<\phi(P_1), \cdots, \phi(P_g)\right> \PUNKT
\end{align*}
It is then easy to conclude that, in this case, the Gauss map is dominant and generically finite, with degree $\binom{2g-2}{g-1}$. 

\noindentTX Furthermore, in the particular case in which $\Cu$ is a genus $3$ non-hyperelliptic curve, which we assume to be embedded in $\PR^2$ via the canonical map, $\Theta$ is smooth and the Gauss map is then nothing but the map which associates, 
to a divisor $P+Q$, the line in $\PR^2$ spanned by $P$ and $Q$ if $P\neq Q$ and the tangent line $T_P(\Cu)$ at $P$ if $P=Q$. In particular, the Gauss map $\GAU$ is a covering of degree $6$ of $\PR^2$ branched on $\Cu^{\vee}$, the dual curve of $\Cu$, which has $28$ nodes corresponding to the bitangent lines of $\Cu$, and $24$ cusps corresponding to the tangent lines passing through a Weierstrass points of $\Cu$.
\end{ex1}
However, we will see in Remark \ref{rmkGauss} that this good behavior of the Gauss map arises in more general situations.
There is furthermore a close connection between the property for the Gauss map of a reduced and irreducible divisor $\NDiv$ of being dominant, and the property for $\NDiv$ of being ample and of general type. More precisely,
it is known that a divisor $\NDiv$ on an abelian variety is of general type if and only if there is no non-trivial abelian subvariety whose action on $A$ by translation stabilizes $\NDiv$. Indeed, the following theorem holds:  %


\begin{thm1} \cite[Theorem 4 - Ueno's theorem]{KAWVIE} \label{Uenothm} 
Let $V$ a subvariety of an abelian variety $A$. Then there exist an abelian subvariety $B$ of $A$ and an algebraic variety $W$ which is a subvariety of an abelian variety such that   
\begin{itemize}
\item $V$ is an analytic fiber bundle over $W$ whose fiber is $B$,
\item $\kappa(W) = dim W = \kappa(B)$. 
\end{itemize}
$B$ is characterized as the maximal connected subgroup of $A$ such that $B + V \subseteq V$.
\end{thm1}

\begin{rmk} \label{rmkGauss} We can conclude that, for a reduced and irreducible divisor $\NDiv$ on an abelian variety $A$, the following are equivalent:
\begin{itemize}
\item[1)] The Gauss map of $\NDiv$ is dominant and hence generically finite. 
\item[2)] $\NDiv$ is an algebraic variety of general type. 
\item[3)] $\NDiv$ is an ample divisor. 
\end{itemize}
Indeed, we recall that a divisor $\NDiv$ on an abelian variety is ample if and only if it is not translation invariant under the action of any non-trivial abelian subvariety of $A$. (see \cite[p. 60]{MumfordAbelian}).
The equivalence of $1)$ and $3)$ follows by (\cite[Proposition 4.4.2]{BLange}), and the idea is that if the Gauss map is not dominant, then $\NDiv$ is not ample because it would be invariant under the action of a non-trivial abelian subvariety. The equivalence of $2)$ and $3)$ follows now easily by applying the Ueno's Theorem \ref{Uenothm}.
\end{rmk}
\begin{rmk}
If $\NDiv$ is a smooth, ample divisor, then the Gauss map of $\NDiv$ is a finite morphism. Indeed, if $\NDiv = div(\theta)$ and the Gauss map $G$ contracted a curve $\Cu$, then without loss of generality we could suppose that $\frac{\partial \theta}{\partial z_j}$ is identically $0$ on $\Cu$ for every $j = 1 \cdots g-1$, and that $\frac{\partial \theta}{\partial z_g}$ has no zeros on $\Cu$. Thus, it would follow that $\omega_{\NDiv}|_{\Cu} \cong \OC_{\Cu}$, which would contradict the fact that $\omega_{\NDiv}$ is ample on $\NDiv$, $\omega_{\NDiv}$ being the restriction of the ample line bundle $\OC_{A}(\NDiv)$ to $\NDiv$.

\noindentTX In particular, an ample smooth surface in an abelian $3$-fold $A$ is a minimal surface of general type.
\end{rmk}

\section{On the canonical map of the bidouble cover of a principal polarization} \label{basicsectionssect}
In this section, we consider surfaces in the polarization of a general $(1,2,2)$-polarized Abelian 3-fold $(A, \mathcal{L})$, 
which are étale bidouble covers of the principal polarization $\Theta$ of a Jacobian threefold. 
More precisely, we are interested in the canonical map of a surface $\Surf$ in the polarization of a general $(1,2,2)$-polarized 
Abelian 3-fold $A$, such that there exists a pullback diagram 
\begin{align} \label{diagrampullbackTheta}
\begin{diagram} 
\Surf & \rInto & A&  \\
\dTo_{p} &  & \dTo_{p} \\
\Theta & \rInto & \Jac
\end{diagram}
\end{align} 
where $p: A \longrightarrow \Jac$ denotes an isogeny onto the Jacobian variety of a smooth quartic plane curve $\Du$, with kernel 
\begin{equation} \label{KERGEQ}
 \GrpJ := Ker(p) \cong \Z_2 \times \Z_2 \PUNKT
\end{equation}
\begin{rmk} \label{CHOICENOTGENERAL}
We remark that those surfaces are not general in their linear system. Indeed, the number of moduli for the couples $(A, \Surf)$, where $(A, \mathcal{L})$ is a $(1,2,2)$-polarized abelian threefold and $\Surf$ belongs to the linear system $|\mathcal{L}|$ is $6$. On the other hand, an ètale bidouble cover as in Diagram \ref{diagrampullbackTheta} is defined by a couple of $2$-torsion line bundles on the basis variety $\Jac$. Since there are only finitely many possible choices for the couples of $2$-torsion line bundles on $\Jac$, the number of moduli for étale bidouble covers as in Diagram \ref{diagrampullbackTheta} is $3$. 
\end{rmk}

\noindentTX In the setting of Diagram \ref{diagrampullbackTheta}, it is convenient to introduce some notations. We express $\Jac$ as a quotient in the form $\Jac \cong \bigslant{\C^3}{\Lambda}$, 
where $\Lambda$ is the lattice in $\C^3$ defined as
\begin{equation} \label{choosendecom}
 \Lambda := \Lambda_{1} \oplus \Lambda_2 \BEISTRICH
\end{equation}
with $\Lambda_1 := \tau \Z^3$, $\Lambda_2:=\Z^3$ and $\tau$ is a general point in the Siegel upper half-space $\mathcal{H}_3$ of $3 \times 3$ symmetric 
matrices with positive definite imaginary part. 
The polarization $\mathcal{O}_{\Jac}(\Theta)$ on $\Jac$, which can be assumed up to translation to be of characteristic $0$ with respect to the decomposition above in (\ref{choosendecom}),
is defined by the hermitian form $H$ whose matrix with respect to the standard basis of $\C^3$ is $(\IMh{\tau})^{-1}$. In particular, since the line bundle on $\Jac$ is assumed to be of characteristic $0$, the divisor $\Theta$ is the zero locus of the Riemann theta function:
\begin{equation} \label{RIEMANNTHETA}
\theta_0(z, \tau) := \sum_{n \in \Z^3} e^{\pi i \cdot \traspose{n}\tau n + 2 \pi i \cdot \traspose{n}z} \PUNKT
\end{equation}
By virtue of the pullback diagram in (\ref{diagrampullbackTheta}), we may assume that
\begin{equation*}
(A, \mathcal{L}) =  \left(\bigslant{\C^3}{\Gamma},  p^* \mathcal{O}_{\Jac}(\Theta) \right) \BEISTRICH
\end{equation*}
where $\Gamma$ is the lattice
\begin{equation*}
\Gamma := \SET{\lambda \in \Lambda}{\IMh{H}(\lambda, \eta_j) \in \Z, \ \ j = 1,2} \BEISTRICH
\end{equation*}
where $\eta_1, \eta_2$ is a couple of $2$-torsion points in $\Jac[2]$. We recall that $\IMh{H}$ is an alternating bilinear form which is moreover integrally valued on the lattice $\Lambda$: for every $u$, $v$ in $\Z_3$ it holds: 
\begin{equation} \label{bHpair}
\IMh{H}(\tau \cdot u, v) = \traspose{u}\cdot v
\end{equation}

\noindentTX The decomposition (\ref{choosendecom}) induces a real decomposition of $\C^3$ into the direct sum of two real subvector spaces $V_j:= \Lambda_j \otimes_{\Z} \R$ defined by extending the scalars, and
by intersection we have an analogous decomposition of each of the sublattices 
\begin{equation*}
\Gamma \subseteq \Lambda \subseteq \Gamma(\mathcal{L}) := \Lambda + \eta_1 \Z + \eta_2 \Z  \PUNKT
\end{equation*}
With $j = 1,2$, we denote the members of this decomposition as follows
\begin{align*}
\Gamma_j:= \Gamma \cap V_j \subseteq \Lambda_j \subseteq \Gamma(\mathcal{L})_j := \Gamma(\mathcal{L}) \cap V_j \PUNKT
\end{align*}
Without loss of generality, we can assume that both $\eta_1$ and $\eta_2$ belong to $V_1$. Moreover, since we assume the polarization of $A$ to be of type $(1,2,2)$, the elements $\eta_1$ and $\eta_2$ must be orthogonal to each other with respect to the Weil pairing:
\begin{equation} \label{bpair}
\WEIL \colon \Jac[2] \times \Jac[2] \longrightarrow \bigslant{\frac{1}{4} \Z}{\frac{1}{2} \Z} \cong \Z_2 \PUNKT
\end{equation}
For the reader convenience, we recall that the Weil pairing in (\ref{bpair}) is the natural symplectic pairing on $\Jac[2]$ induced by the alternating bilinear form $\IMh{H}$. The latter orthogonality condition ensures indeed that the quotient 
$\bigslant{\Gamma(\mathcal{L})_1}{\Gamma_1}$ is isomorphic to $\Z_2 \times \Z_2$ and generated by $\eta_1$ and $\eta_2$.

The kernel of the isogeny $p$ in the pullback diagram in (\ref{diagrampullbackTheta}), which we denoted by $\GrpJ$, is then isomorphic to $\bigslant{\Gamma(\mathcal{L})_2}{\Gamma_2}$, and we can find two involutions $a$ and $b$ in $\GrpJ$ such that $\WEIL(\eta_1, a) = \WEIL(\eta_2, b) = 0$ and $\WEIL(\eta_0, b) = \WEIL(\eta_1, a) = 1$.

The Weil pairing $\WEIL$, being a perfect pairing of $\Z_2$-modules, allows us to identify $\bigslant{\Gamma(\mathcal{L})_1}{\Gamma_1}$ with the character group $\GrpJ^*$: indeed,
every element $\tau \cdot \gamma$ of $\bigslant{\Gamma(\mathcal{L})_1}{\Gamma_1}$ can be identified with the character defined on $\GrpJ$ as follows:
\begin{equation} \label{chardef}
\gamma(g) := e^{\pi i \WEIL(\tau \gamma, g)} = e^{2\pi i \traspose{\gamma} \cdot g}
\end{equation}

\noindent The vector space $H^0(A, \mathcal{L})$ has dimension $4$, and it is generated by the theta functions 
$\theta_{0}$, $\theta_{\alpha}$, $\theta_{\beta}$, $\theta_{\alpha + \beta}$ which are defined, for every element $\gamma$ of $\GrpJ^{*}$ as follows
\begin{equation} \label{thetadef1}
\theta_{\gamma}(z, \tau) := e^{\pi i \traspose{\gamma} \tau \gamma +2\pi i \traspose{\gamma}z }\theta_0(z + \gamma) \PUNKT
\end{equation}

In conclusion, one can easily see that, for every $\gamma \in \GrpJ$ and every $g$ in $\GrpJ$, we have
\begin{equation} \label{THETAGAMMA}
\theta_{\gamma}(z + g) = e^{\pi i \lambda(\gamma, g)} = \gamma(g) \theta_{\gamma}(z) \PUNKT
\end{equation}


\begin{not1}\label{THETAGAMMANOTATION}
From the above definition in (\ref{thetadef1}), it follows that the zero locus of each section $\theta_{\gamma}$ in $A$ is just a translated of $\Surf$, where
\begin{equation*}
(\theta_{\gamma})_0 =  \Surf + \gamma
\end{equation*}
We will denote this locus by $\BBF{\Surf_g}$, where $g$ is the only element of $\GrpJ$ whose action leaves the section $\theta_{\gamma}$ invariant according to Equation \ref{THETAGAMMA}.
\end{not1}

\begin{rmk}
The multiplication by $-1$ on $A$ leaves each of these surfaces invariant, since each function $\theta_{\gamma}$ is an even function. 
Furthermore, it can be easily seen 
that the base locus $\mathcal{B}(\mathcal{L})$ of the line bundle $\mathcal{L}$ is a set of cardinality $16$ which containes precisely those $2$-torsion points on $A$ on which every function $\theta_{\gamma}$ vanishes with odd multiplicity.


\noindent The following theorem shows that, on a surface $\Surf$ which is an étale bidouble cover of a theta divisor in a general principally polarized abelian threefold 
the differential of the canonical map is injective at every point.
However, the same theorem states that the canonical map fails to be injective precisely on the canonical curves of $\Surf$ of the type $\Surf \cap \Surf_{\gamma}$.

\end{rmk}
\begin{thm1} \label{injectivedifferential}
Let $\Surf$ be an étale bidouble cover of a theta divisor in a general principally polarized abelian threefold, and let be
\begin{align*} 
\begin{diagram} 
\Surf & \rInto & A&  \\
\dTo_{p} &  & \dTo_{p} \\
\Theta & \rInto & \Jac
\end{diagram}
\end{align*} 
the corresponding pullback diagram.

\begin{itemize}
\item[1)] If $P$ and $Q$ are two distinct points on $\Surf$ whose image with respect to the canonical map of $\Surf$ is the same, then one of the following cases occurs: 
\begin{itemize}
\item[$\bullet$] $Q = -g.P$ for some non-trivial element $g$ of $\GrpJ$. This case arises precisely when $P$ and $Q$ belong to the canonical curve $\Surf \cap \Surf_{g}$, where $\gamma(g)=1$. 
\item[$\bullet$] $Q = g.P$  for some non-trivial element $g$ of $\GrpJ$.  This case arises precisely when $P$ and $Q$ belong to the translate $\Surf_{h}$, for every $h \in \GrpJ - \{g\}$.
\item[$\bullet$] $P$ and $Q$ are two base points of $|\Surf|$ which belong to the same $\GrpJ$-orbit.
\end{itemize}
\item[2)] The differential of the canonical map of $\Surf$ is injective at every point. 
\end{itemize}
\end{thm1}
To prove Theorem \ref{injectivedifferential}, we begin by fixing some notations and recalling some known facts about the behavior of the 
canonical map of the surfaces in a polarization of type $(1,1,2)$, which has been developed by F. Catanese and F-O Schreier in their joint work \cite{CataneseSchreyer1}.
The strategy we use to prove the first claim is essentially then to consider the canonical projections which naturally arise from the isogenies onto $(1,1,2)$-polarized abelian threefolds.\newline 

\noindentTX To this purpose, for every non-trivial element $g$ of $\GrpJ$, we denote by $A_g$ the $(1,1,2)$-polarized abelian threefold obtained as the quotient of $A$ by $g$, 
by $q_{g}$ the projection of $A$ onto $A_g$, and by $\TD_g$ the image $q_g(\Surf)$ in $A_g$.  
We have then the following diagram:
\begin{align} \label{diagramproj}
\begin{diagram} 
		A &\lInto&\Surf 					&	\rTo^{\phi_{\Surf}} 	&	\Sigma	  &\rInto& \PR^5				\\ 	
		\dTo_{q_g} && \dTo_{q_g} 	&			&\dTo &&\dDashto		\\
		A_{g}&\lInto&\TD_{g} & \rTo^{\phi_{\TD_g}}			&\Sigma_g &\rInto&\PR^3		
 \end{diagram}
 \end{align}
where $\phi_{\Surf}$ denotes the canonical map of $\Surf$ and $\Sigma$ its image in $\PR^5$. \newline

\noindentTX The canonical map $\phi_{\TD_g}$ 
is defined by the theta functions $[\theta_{\gamma}, \frac{\partial \theta_{0}}{\partial z_1},  \frac{\partial \theta_{0}}{\partial z_2},  \frac{\partial \theta_{0}}{\partial z_3}]$, where $\theta_{\gamma}$ is the only holomorphic section of $\Surf$ which is invariant under the action of $g$ and anti-invariant under the action of every other non-trivial element.



\noindentTX Under this setup, we can state the following theorem (see \cite[Theorem 6.4]{CataneseSchreyer1})
\begin{thm1} \label{CATSCHREYER}
Let $\TD$ be a smooth divisor yielding a polarization of type $(1,1,2)$ on an Abelian threefold. Then the canonical map of $\TD$ is, in general, a birational morphism onto a surface $\Sigma$ of degree $12$ in $\PR^3$.

\noindentTX In the special case  where $\TD$ is the inverse image of the theta divisor in a principally polarized Abelian threefold, 
the canonical map is a degree $2$ morphism onto a sextic surface $\Sigma$ in $\PR^3$. In this case the singularities of $\Sigma$ are in general: 
a plane cubic $\Gamma$ which is a double curve of nodal type for $\Sigma$ and, according to \cite[Definition 2.5]{Catanese1981}, 
a strictly even set of $32$ nodes for $\Sigma$. Also, in this case, the normalization of $\Sigma$ is in fact the quotient of $\TD$ 
by an involution $i$ on $A$ having only isolated fixed points (on $A$), of which exactly $32$ lie on $\TD$.
\end{thm1}
 
For the sake of exposition, it is useful to recall here the basic ideas of the proof of Theorem \ref{CATSCHREYER}. We assume that 
$\TD$ is defined as the zero locus of the Riemann Theta function $\theta_0$ (as defined in \ref{thetadef1}) in a general $(1,1,2)$-polarized abelian
threefold, which we shall denote by $A'$. In this case, $\TD$ is the pullback of a divisor $\Theta$ via an isogeny of degree $2$ whose kernel is generated by an involution $\delta$.

The vector space of the global holomorphic section of the polarization $\mathcal{L}'$ on $A'$ is generated by two even functions:
\begin{equation} \label{secondDEFgamma}
 H^0(A', \mathcal{L}') = \left<\theta_0, \theta_{\gamma}\right> \BEISTRICH
\end{equation}
where $\theta_{\gamma}$ is anti-invariant with respect to the action of $\delta$ by translations.
On the other hand, the derivatives $\frac{\partial \theta_0}{\partial z_j}$ are odd functions, and by virtue of Proposition \ref{desccanonicalmap}, we conclude that every holomorphic section of the canonical bundle of $\TD$ is invariant under the action of the involution:
\begin{equation} \label{involutioniota}
\iota \colon z \mapsto -z + \delta \PUNKT
\end{equation}
Hence we can conclude that the canonical map of $\TD$ factors through the quotient $Z := \bigslant{\TD}{\iota}$ and it cannot be, 
in particular, birational. 

\noindentTX The quotient map on the surface $Z$ is however birational, as we are going to prove. Firstly, we count the number of fixed points of $\iota$ on $\TD$. Let us consider the commutative diagram
\begin{align}\label{put}
\begin{diagram} 
A' & \lInto & \TD & \rTo & Z \\
\dTo &  & \dTo &  & \dTo \\
\Jac & \lInto & \Theta & \rTo & Y := \bigslant{\Theta}{\left<-1\right>}
\end{diagram}
\end{align} 
By our generality assumption on $A'$, we can assume that $\Jac$ is the Jacobian variety of a smooth quartic plane curve $\Du$. \newline 

\noindentTX We saw in Example \ref{Gaussjacob} that every point of the Theta divisor $\Theta$ of a quartic $\Du$ can be represented by an effective divisor of degree $2$ on $\Du$, and the Gauss map $G \colon \Theta \longrightarrow \PR^2$ is the map which associates, to a divisor $u+v$, the line in $\PR^2$ spanned by $u$ and $v$ 
if $u\neq v$ and the tangent line to $\Du$ at $u$ if $u=v$. 

\noindentTX Moreover, the multiplication by $-1$ on $\Jac$ corresponds to the Serre involution $\mathcal{L} \mapsto \omega_{\Du} \otimes \mathcal{L}^{\vee}$ 
on $Pic^2(\Du)$, which can be expressed on $\Theta$ as the involution which associates, to the divisor  $u+v$ on $\Du$, the unique divisor $u'+v'$ 
such that $u+v+u+v'$ is a canonical divisor on $\Du$.

\noindentTX It can be now easily seen that the projection of $\Theta$ onto $Y$ in Diagram (\ref{put}) is a covering branched on $28$ points, which correspond
precisely to the $28$ bitangents of $\Du$, and in particular we conclude that $Y$ has precisely $28$ nodes. 
This set corresponds to the set of odd $2$-torsion points of the Jacobian on which $\Theta$ has odd multiplicity, and we have
:
\begin{equation} \label{twopoints}
\begin{split} 
 \Jac[2]^{-}(\Theta) &= \SET{z \in \Jac[2]}{mult_{z}\Theta \text{ is odd}} \\ 
 &= \SET{z \in \Jac[2]}{E(2z_1, 2z_2) = 1 \ \ \text{$(mod \  2)$}}   \BEISTRICH
 \end{split}
 \end{equation}
 where $E := \IMh H$ is the alternating bilinear form associated to the polarization on $\Jac$, which is $\Z$-valued on the lattice $\Lambda$ 
 Among the odd $2$-torsion points in (\ref{twopoints}) above, there are some which belong to the image in $\Jac$ of the fixed locus of the involution $\iota$ in $A$. The set containing them is precisely
\begin{align*}
p(\mathcal{F}ix(\iota)) \cap \Jac[2]^{-}(\Theta) &=  \SET{[v ]\in \Jac[2]^{-}(\Theta)}{[2v] = \delta \ \ \ \text{in  $A'$}} \\ &= \SET{z \in \Jac[2]^{-}(\Theta)}{E(2z, 2\eta) = 1 \ \ \text{$(mod \  2)$}} \PUNKT
\end{align*}
This latter set has cardinality equal to $16$, and thus the involution $\iota$ has exactly $32$ fixed points on $\Surf$. 
Hence the singular locus of $Z$ is a strictly even set of $32$ nodes, and it holds (see \cite[Proposition 2.11]{Catanese1981}) that 
\begin{align} \label{equationEulercover}
\chi(\TD) = 2\chi(Z) - 8 = 2 \PUNKT
\end{align}
Thus, $\chi(Z) = 5$. On the other hand, since the canonical map of $Z$ factors through $\iota$, we have $p_g(Z) = 4$, and hence we conclude that $q(Z) = 0$. 
\noindentTX The map $Z \longrightarrow Y$ in Diagram \ref{put} is a double cover which is unramified except over the remaining $12$ nodes of $Y$. 

\noindentTX We have just proved that the degree of the canonical map of $\TD$ is at least $2$. The following lemma ensures that the degree is
exactly $2$ (see also \cite{CataneseSchreyer1} for a proof).

\begin{lem1} \label{thetafunctions112}
The canonical map of $Z$ is birational. Hence, the canonical map of $\TD$ is of degree $2$, and its image is a 
surface $\Sigma$ of degree $6$ in $\PR^3$.
\end{lem1}
\begin{proof}
We recall that the canonical map $\CAN{\TD}$ of $\TD$ factors through the involution $\iota \colon z \mapsto -z + \delta$ as in \ref{involutioniota}, and all canonical sections of $\TD$ are anti-invariant with respect to $\iota$.


\noindentTX It is easily seen that the degree of the canonical map of the quotient $Z := \bigslant{\TD}{\iota}$ is at most $3$. Indeed, the Gauss map $\GAU$ of $\TD$ factors through $Z$ via the Gauss map on the quotient $\GAU_Z: Z \longrightarrow {\PR^2}$. 
The latter map is of degree $6$ and invariant with respect to the involution $(-1)_Z$ induced by the multiplication by $(-1)$ of $\TD$, while the canonical map is not.  
Moreover, $\GAU$ is ramified on a locus of degree $24$ which represents the dual curve $\Du^{\vee}$ of $\Du$ (of degree $12$) counted with multiplicity $2$.
If we denote by $\pi: \PR^3 \dashrightarrow \PR^2$ the projection which forgets the first coordinate, 
its restriction to $\Sigma$ must have a ramification locus $R$ whose degree is divisible by $12$, and we have:
\begin{equation*}
deg \Sigma = deg K_{\Sigma} = -3 deg(\pi|_{\Sigma}) + deg R \PUNKT
\end{equation*}
However, $deg(\pi|_{\Sigma}) = deg \Sigma$, and hence $4deg\Sigma = deg R$. In particular, the canonical map of $Z$ can have degree $1$ or $2$. \newline

\noindentTX Let us suppose by contradiction that $\CAN{Z}$ is of degree $2$. Then $\CAN{Z}$ would be invariant under the action of an involution $j$ on $Z$, which would imply that the map $\GAU_Z$ is invariant under the action of the whole group generated by $(-1)_Z$ and $j$. This latter group has a natural faithful representation as a subgroup of the monodromy group of Gauss map of $Z$, which is isomorphic to $S_3$, the symmetric group of degree $3$. Hence, 
\begin{equation*}
H := \left< (-1)_Z , j \right> \cong S_3
\end{equation*}
\noindentTX The ramification locus in $Y$ of $\phi_{Y}$ has two components: 
\begin{align*} 
\BITAN &:= \SET{[2p]}{p \in \Du} \\
\TAN &:= \SET{[p+q]}{\left< p,q \right>.\Du = 2p+q+r \ \text{for some $q$,$r$ on $\Du$}} \PUNKT
\end{align*}
The component $\TAN$ has to be counted twice on $Y$, while $\BITAN$ has multiplicity $1$. 
The components of the counterimage of $\BITAN$ in $Z$ are both of multiplicity $1$ and the group $H$ acts on them. 
Clearly, the involution $(-1)_Z$ exchanges them, while both components must be pointwise fixed under the action of $j$. 
But the product $(-1)_Z \cdot j$, has order $3$ (because $H \cong \mathcal{S}_3$),
hence it must fix both components, and we reach a contradiction.
\end{proof}

\noindentTX The canonical models of surfaces with invariants $p_g=4$, $q=0$, $K^2 = 6$ and birational canonical map have been studied extensively by F. Catanese 
(see \cite{Catanese1984}). In this case, 
there is a symmetric homomorphism of sheaves:
\begin{equation*}
\alpha = \begin{bmatrix}\alpha_{00} & \alpha_{01} \\
				    \alpha_{01} & \alpha_{11} 
 \end{bmatrix}: (\mathcal{O}_{\PR^3} \oplus \mathcal{O}_{\PR^3}(-2))^{\vee}(-5) \longrightarrow \mathcal{O}_{\PR^3} \oplus \mathcal{O}_{\PR^3}(-2) \BEISTRICH
\end{equation*}
where $\alpha_{00}$ is contained in the ideal $I$ generated by $\alpha_{01}$ and $\alpha_{11}$. 
The canonical model $Y$ is defined by $det (\alpha)$ and the closed curve $\Gamma$ in Theorem \ref{CATSCHREYER} is a cubic defined by the ideal $I$ and contained in the projective plane $\alpha_{11} = 0$. 
\begin{rmk} \label{rmk112}
\noindentTX In our situation, the double curve $\Gamma$ in the canonical image of the quotient $Z$ is nothing but the image of the canonical curve $\mathcal{K}$ defined as the zero locus of $\theta_{\gamma}$ in $\TD$ (refer again to (\ref{secondDEFgamma}) for the definition of $\theta_{\gamma}$). Indeed, let us denote by $\bar{\mathcal{K}}$ the image of $\mathcal{K}$ in $Z$. 
The curve $\mathcal{K}$ does not contain fixed points of $\iota$, whence $\mathcal{K}$ and $\bar{\mathcal{K}}$ are isomorphic. 
Furthermore, the curve $\bar{\mathcal{K}}$ is stable under the involution $(-1)$ on $\TD$, and the canonical map of $Z$ is of degree $2$ on $\bar{\mathcal{K}}$. In conclusion, the image of $\bar{\mathcal{K}}$ in $\Sigma$ is a curve of nodal type, and then it must be exactly $\Gamma$. Moreover, the set $\mathcal{P}$ of pinch points on $\Gamma$ is exactly the image of the set of the $2$-torsion 
points of $\TD$ which lie on the canonical curve $\mathcal{K}$. This latter set is in bijection with the set
\begin{equation*} 
\left\{(x,y) \in \Z_2^3 \times \Z_2^3\ \ | \ x_2y_2 + x_3y_3 =  1  \ \right\} \BEISTRICH
\end{equation*}
which consists precisely of $24$ points. The set $\mathcal{P}$ is then precisely the branch locus of the map $\bar{\mathcal{K}} \longrightarrow \Gamma$, 
which has degree $2$ and factors with respect to the involution $(-1)_Z$. 
\end{rmk}

\noindentTX We are now in position to prove the first claim of Theorem \ref{injectivedifferential}.
\begin{proof}[of Theorem \ref{injectivedifferential} point 1)]
According to Diagram \ref{diagramproj}, for every involution in the kernel $\GrpJ$ of the isogeny $p: A \longrightarrow \Jac$ there is a diagram
\begin{align} \label{diagramproj2}
\begin{diagram} 
		A &\lInto&\Surf 					&	\rTo^{\phi_{\Surf}} 	&	\Sigma	  &\rInto& \PR^5				\\ 	
		\dTo_{q_g} && \dTo_{q_g} 	&			&\dTo &&\dDashto		\\
		A_{g}&\lInto&\TD_{g} & \rTo^{\phi_{\TD_g}}			&\Sigma_g &\rInto&\PR^3		\\
&&&\rdTo(1,1)^{} Z_g  \ruTo(1,1)&& 
 \end{diagram}
 \end{align}
where $Z_g$ denotes the quotient of $\TD_{g}$ by the involution $z \mapsto -z+h$ in $A_g$, and where $h \in \GrpJ - \{1,g\}$. By Theorem \ref{CATSCHREYER}, we obtain that,
for every non-trivial element $g$ of $\GrpJ$ only one of the two possible cases can occur:
\begin{itemize}
 \item $q_g(U) = q_g(V)$. In this case, $V = g.U$, and thus both $U$ and $V$ lie on $\Surf \cap \Surf_{h} \cap \Surf_{g+h}$, where $h \in \GrpJ - \{1,g\}$, and the 
 claim of the theorem follows.
 \item $q_g(U) \neq q_g(V)$ but $(\pi_{g} \circ q_g)(U) = (\pi_{g} \circ q_g)(V)$, where $\pi_g$ denotes the projection of $\TD_{g}$ onto $Z_g$ in Diagram \ref{diagramproj2}. 
 In this case, we have that $V = -U + h$ for some $h \in \GrpJ - \{1,g\}$, and both $U$ and $V$ belong to $\Surf \cap \Surf_{h}$.
 \item $(\pi_{g} \circ q_g)(U) \neq (\pi_{g} \circ q_g)(V)$. In this case, by applying Theorem \ref{CATSCHREYER} together with
 Remark \ref{rmk112}, we have that $U$ and $V$ belong to $\Surf_{g}$ and $V = -U + g$, and we conclude again that $\Surf \cap \Surf_{g}$.
\end{itemize}
\end{proof}
The proof of the second claim of Theorem \ref{injectivedifferential} requires some more efforts. The first important step toward proving it, is to describe the pullback surface $\Surf$ as a quotient of the symmetric product of a certain genus $9$ curve $\Cu$, which is an ètale bidouble cover of a general algebraic curve $\Du$ of genus $3$.\newline

\noindentTX The next lemma shows that an ètale bidouble cover $\Cu$ of a general algebraic curve $\Du$ of genus $3$ is a genus $9$ tetragonal curve with a very ample theta characteristic with four global independent sections.´

\begin{lem1} \label{fundamentallemmagenus9}
Let $(A, \mathcal{L})$ be a general $(1,2,2)$-polarized Abelian 3-fold, let $p: A \longrightarrow \Jac$ be an isogeny onto the
Jacobian of a general algebraic curve $\Du$ of genus $3$, defined by two elements $\eta_1$ and $\eta_2$ belonging to $\Jac[2]$ 
which satisfy the condition that $\WEIL(\eta_1, \eta_2) = 0$. Let us moreover consider the algebraic curve $\Cu$ obtained by pulling back to $A$ the curve $\Du$ via the isogeny $p$ (here we consider $\Du$ to be embedded in its Jacobian $\Jac$).


Then, the following hold true:
\begin{itemize}
\item The genus $9$ curve $\Cu$ admits $\mathcal{E}$ and $\mathcal{F}$ two distinct $\GrpJ$-invariant $g_{4}^1$'s, with 
$\mathcal{E}^{2} \ncong \mathcal{F}^{2}$ and
\begin{equation*}
 h^0(\Cu, \mathcal{E}) = h^0(\Cu, \mathcal{F}) = 2 \PUNKT
\end{equation*}
\item The line bundle $\mathcal{M} := \mathcal{E} \otimes \mathcal{F}$ is a very ample theta characteristic of type $g_8^3$.
\item The image of $\Cu$ in $\PR^3 = \PR(\mathcal{M})$ is a complete intersection of the following type:
\begin{equation} \label{corollaryscrolleq}
\Cu : \begin{cases}
X^2 + Y^2 + Z^2 + T^2 = 0 \\
q(X^2,Y^2,Z^2,T^2) = XYZT \BEISTRICH
\end{cases}
\end{equation}
where $q$ is a quadric, and $[X,Y,Z,T]$ are coordinates on $\PR^3$ such that two generators $a$ and $b$ of kernel of the isogeny $\GrpJ:= Ker(p)$ act as follows:
\begin{align*}
\begin{split}
a.[X,Y,Z,T] &= [X,Y,-Z,-T]\\
b.[X,Y,Z,T] &= [X,-Y,Z,-T] \PUNKT
\end{split}
\end{align*} 
\item The étale bidouble cover $p: \Cu \longrightarrow \Du$ can be identified with the restriction to
$\Cu$ of the rational map $\psi: \PR^3 \dashrightarrow \PR^3$ defined by the squares of the coordinates:
\begin{equation*}
\psi: [X,Y,Z,T] \daschmapsto [x,y,z,t] := [X^2,Y^2,Z^2,T^2]  \PUNKT
\end{equation*}
and the equations of $\Du$ in $\PR^3 = \PR[x,y,z,t]$ are, according to the defining equations of $\Cu$ in (\ref{corollaryscrolleq}), in the following form:
\begin{equation*} 
\Du : \begin{cases}
x + y + z + t = 0 \\
q(x,y,z,t)^2 = xyzt \PUNKT
\end{cases}
\end{equation*}
\end{itemize}
\end{lem1}
\begin{proof}
 We show that $\Cu$ is tetragonal. 
\noindentTX Let us fix a point $Q_0 \in \Cu$ whose image $q_0 \in \Du$ with respect to $p$ does not lie on a bitangent line of $\Du$, and let us denote by $\AbelM$ the Abel map defined respect to a point different from $q_0$. \newline

\noindentTX We claim first that there exists a point $\zeta$ on $\Theta$ such that:
\begin{itemize}
 \item The point $\zeta$ does not belong to any translated of $\Theta$ with the 2-torsion elements $\eta_1$, $\eta_2$ and $\eta_1 + \eta_2$. This condition is equivalent to require that $\zeta$ does not belong to the image with respect to $p$ of the base locus of the polarization on $A$.
 \item For every special divisor of degree $3$ and every $2$-torsion element $\eta$ we have that 
 \begin{equation*}
  \zeta \neq \AbelM(q_0) - \RC - \eta - \AbelM(D) \BEISTRICH
 \end{equation*}
where $\RC$ is the vector of Riemann constants (See \cite[Theorem 2 p. 100]{Nara}). \newline
\end{itemize}



\noindentTX Indeed, if the second condition does not hold for a certain point $\zeta$ of Theta divisor, then it exists a $2$-torsion point $\eta \in \Jac[2]$ and a special divisor $D$ of degree $3$ on $\Du$ such that:
\begin{equation*}
\zeta = \AbelM(q_0) - \RC - \eta - \AbelM(D) \BEISTRICH
\end{equation*}
and it follows that, in particular (recall that $\AbelM(K) = -2\RC$, where $K$ is a canonical divisor on $\Du$)
\begin{equation*}
0 = 2 \zeta = \AbelM(2q_0 + K - 2D) - 2\eta =  \AbelM(2q_0 + K - 2D) \PUNKT
\end{equation*}
Hence, by Abel's theorem, the divisor $2(D-q_0)$ is a canonical divisor. But $D$ is supposed to be a special divisor of degree $3$, hence linearly equivalent to $K-r$ where $r$ is some point on $\mathcal{D}$. ($D$ is a $g_3^1$ on the curve $\mathcal{D}$). That means, in particular, that:
\begin{align*}
K \equiv 2(D-q_0) \equiv 2(K - r - q_0) \BEISTRICH
\end{align*}
and we conclude that $r + q_0$ is an odd theta-characteristic. However, because the odd theta characteristics correspond to the bitangent lines on $\mathcal{D}$ we assume that $q_0$ does not lie on any bitangent line of $\mathcal{D}$, we can exclude this case. This proves the first claim. \newline 

\noindentTX Now we prove that there exists a theta characteristic $\mathcal{M}$ on $\Cu$ which yields an embedding of $\Cu$ in $\PR^3$, such that the image as a complete intersection as in the statement of the lemma. We are going to prove that $\mathcal{M}$ is the mobile part of a line bundle on $\Cu$ naturally induced by restricting to $\Cu$ a certain translated of the polarization $\mathcal{L}$ on $A$, which has $4$ sections.

To prove the lemma, we pick now $\zeta$ which satisfy the claim we have just proved. With such a point $\zeta$, we can define:
\begin{equation}  \label{definitionT}
\mathcal{T} := (t_{\AbelM(q_0)-\zeta}^*\mathcal{L})|_{\Cu}
\end{equation}
where $t_x$ denotes in general the translation by a point $x$ of $A$. By Definition \ref{definitionT}, for every $\gamma$ the following are then easily seen to be holomorphic sections of $\mathcal{T}$:
\begin{equation*}
s^{q_0}_{\gamma} := t_{\AbelM(q_0)-\zeta}^* \theta_{\gamma}|_{\Cu} \PUNKT
\end{equation*}
\noindentTX Our goal is then to show that $|\mathcal{T}|$ is a linear system on $\Cu$ with $4$ base points and of degree $12$, and that its mobile part defines a $g^3_8$ on $\Cu$.
\noindentTX To this purpose, we first notice that, if we consider a point $X \in \mathcal{C}$ and its image $x$ in $\mathcal{D}$ with respect to $p$, we have that:
\begin{equation} \label{EXPSECT}
s^{q_0}_{\gamma}(X) = \theta_{\gamma}(\AbelM(x)-\AbelM(q_0)-\zeta) \PUNKT
\end{equation}
Hence, $s^{q_0}_{\gamma}(X)$ vanishes if and only if $ \theta_{0}(\AbelM(x)-\AbelM(q_0)-(\zeta + \gamma))= 0$. On the other side, one can show that (see \cite[Lemma 2, p.112]{Nara}):
\begin{align*}
div (\theta_0(\AbelM(x)-\AbelM(q_0)-\zeta - \gamma) ) &= q_0 + D_{\gamma} 
\end{align*}
where $D_{\gamma}$ is a divisor of degree $2$ on the quartic curve $\Du$ which is independent on $q_0$, for which:
\begin{align*}
\AbelM(D_{\gamma}) &= \zeta + \eta - \RC \PUNKT 
\end{align*}
Thus:
\begin{equation*}
div (s^{q_0}_{\gamma}) = \GrpJ.Q_0 + p^*(D_{\gamma}) \BEISTRICH
\end{equation*}
where $\GrpJ.Q_0$ is the orbit of $Q_0$ with respect to the action of $\GrpJ$.

\noindentTX This means that $\mathcal{T}$, which is of degree $12$, has a fixed part of degree $4$, and its mobile part $|\mathcal{M}|$ has degree $8$. We can conclude now that 
\begin{equation*}
H^0(A, \mathcal{I}_{\Cu} \otimes t_{\AbelM(q_0)-\zeta}^*\mathcal{L}) = \sum_{\eta} H^0(\Jac, \mathcal{I}_{\Du} \otimes t_{\AbelM(q_0)-\zeta-\gamma}^*\Theta) = 0 \PUNKT
\end{equation*}
Indeed, if $h^0(\Jac, \mathcal{I}_{\Du} \otimes t_{\AbelM(q_0)-\zeta-\gamma}^*\Theta) = 1$ for some $\gamma$ then, for every $p \in \Du$, we would have:
\begin{equation*}
\theta_0(\AbelM(q_0) - \zeta - \gamma - \AbelM(p)) = 0 \BEISTRICH
\end{equation*}
and we would find a special divisor $D$ such that:  
\begin{equation*}
\AbelM(q_0) - \eta - \zeta - \RC = \AbelM(D) \PUNKT
\end{equation*}
However, this would contradict the conditions on $\zeta$. We have then the following exact sequence:
\begin{equation*}
0 \longrightarrow H^0(A, t_{A(q_0)-\zeta}^*\mathcal{L})  \longrightarrow H^0(\Cu, \mathcal{T} ) \longrightarrow H^1(A, \mathcal{I}_{\Cu} \otimes t_{A(q_0)-\zeta}^*\mathcal{L}) \longrightarrow 0 \BEISTRICH
\end{equation*}
from which we can conclude that $h^0(\Cu, \mathcal{M} ) = h^0(\Cu, \mathcal{T} ) \geq 4 $. We can apply the Riemann Roch theorem to conclude that the linear system $|\mathcal{M}|$ is special on $\Cu$. Moreover, its Clifford index is:
\begin{equation*}
Cliff(\mathcal{M}) = 8 - 2|\mathcal{M}| \leq 2 \PUNKT
\end{equation*}
On the other side, the Clifford index of $\mathcal{M}$ cannot be $0$ by the Clifford theorem, so $\mathcal{M}$ is of Clifford index $2$. This implies that $\Cu$ is an
algebraic curve of genus $9$ with Clifford Index $2$, and we conclude that $\Cu$ is a tetragonal curve, 
with $\mathcal{M}$ a linear system $g_3^{8}$. Moreover, $\mathcal{M}$ is very ample. Indeed, if it were not the case, then $\Cu$ would be hyperelliptic, which would imply that $\Du$ also is. This, however, would lead to a contradiction, since we assume $\Du$ to be general. \newline

\noindentTX The sections of $\mathcal{T}$, defined in (\ref{EXPSECT}), define coordinates $[X,Y,Z,T]$ on $\PR^3 = \PR H^0(\Cu, \mathcal{M})$ on which $\GrpJ$ acts as claimed in the statement (according to Equation \ref{THETAGAMMA}, which describes how the group $\GrpJ$ acts on the sections of the polarization of $A$). 

The image of $\Cu$ in $\PR^3$ is the complete intersection of a quadric $\mathcal{Q}_2$ and 
a quartic surface $\mathcal{Q}_4$, both $\GrpJ$-invariant, with $\mathcal{Q}_2 = (q_2)_0$ and $\mathcal{Q}_4 = (q_4)_0$ for certain homogeneous polynomials $q_2$ and $q_4$. 

Without loss of generality, we can assume our quadric to be in the form´ $q_2 = X^2 + Y^2 + Z^2 + T^2$. We claim that also the quartic 
$q_4$ is $\GrpJ$-invariant. 
Indeed, if it were not the case, 
If it were not the case, then we could suppose (using the equation of $Q_2$) that the quartic $q_4$ is of the form
\begin{equation*}   
XYr^2 = ZTs^2 \BEISTRICH
\end{equation*}
where $r$ and $s$ are polynomials in the vector space generated by the three squares $X^2, Y^2, Z^2$. This means that, considered $\PR^3$ with 
coordinates $[x,y,z,t] := [X^2, Y^2, Z^2, T^2]$,
we could write the equation of $\Du$ in the following form:
\begin{equation*}
\Du : \begin{cases}
x + y + z + t = 0 \\
xyr^2 = zts^2 \BEISTRICH
\end{cases}
\end{equation*}
where $r$ and $s$ are two lines. In this case, $\Du$ would be singular in the point in which the lines $r$ and $s$ intersect. 
Hence, $\Du$ would be hyperelliptic, which would contradict our generality assumptions on $\Du$. Hence $q_4$ is $\GrpJ$-invariant, and in particular of the form described in the claim of this lemma.

To conclude the proof, is now enough to observe that the two rulings of the non-singular quadric $\mathcal{Q}_2$ induce two distinct $\GrpJ$-invariant $g_1^4$ which we shall denote by $\mathcal{E}$ and $\mathcal{F}$, whose tensor product is clearly isomorphic to $\mathcal{M}$. This finishes the proof of the lemma.
\end{proof}
We can now easily see that, for a general quartic plane curve $\Du$, an ètale bidouble cover $\Surf$ of its Theta divisor $\Theta$, defined by two $2$-torsion line bundles on $\Du$ which are orthogonal with respect to the Weil pairing, can be geometrically described as a quotient of the form 
\begin{equation*} \label{RepresentS}
\Surf = \bigslant{\Cu \times \Cu}{\Delta_{\GrpJ} \times \Z_2} \BEISTRICH
\end{equation*}
where $\Cu$ is a smooth curve of genus $9$ in $\PR^3$ as in Lemma \ref{fundamentallemmagenus9}, $\GrpJ \cong \Z_2 \times \Z_2$ is the Galois group of the cover $p: \Surf \longrightarrow \Theta$, $\Delta_{\GrpJ}$ is the diagonal subgroup of $\GrpJ \times \GrpJ$ acting naturally on $\Cu \times \Cu$ and $\Z_2$ denotes the group acting by switching the two factors of $\Cu \times \Cu$.

\begin{not1} \label{NOTcanbdlS}
\noindentTX The global holomorphic sections of the canonical bundle of $\Surf$ are clearly those of $\Cu \times \Cu$ which are invariant under the action of the group $\Delta_{\GrpJ} \times \Z_2$: 
\begin{align*}
H^0(\Surf, \omega_{\Surf}) &= H^0(\Cu \times \Cu, \omega_{\Cu} \boxtimes \omega_{\Cu})^{\Delta_{\GrpJ} \times \Z_2}
\end{align*}
These can be seen as quadrics in the projective coordinates $[X_1, Y_1, Z_1, T_1]$ and $[X_2, Y_2, Z_2, T_2]$ of $\PR^3 \times \PR^3$ which are invariant under the same action. Thus, the following is easily seen to be a basis for $H^0(\Surf, \omega_{\Surf})$:
 \begin{equation} \label{polynomialexpressioncanonicalmap}
 \begin{split}
\ETA{12} &:= 		 \begin{vmatrix} 					X_1^2	&	X_2^2\\
                                                        Y_1^2	&	Y_2^2
 			 \end{vmatrix} \ \ \ \ \ \ \ \ 
\ETA{13} := 			  \begin{vmatrix} 	X_1^2	&	X_2^2\\
  											Z_1^2	&	Z_2^2
 				 \end{vmatrix} \ \ \ \ \ \ \ \ \		
\ETA{23} := 			  \begin{vmatrix} 	Y_1^2	&	Y_2^2\\
  											Z_1^2	&	Z_2^2
 				 \end{vmatrix} \\		 			 
\OMEGA{45} &:=				\begin{vmatrix} 					X_1Y_1	&	X_2Y_2\\
  																Z_1T_1	&	Z_2T_2\\
 					 \end{vmatrix} \ \ \		
\OMEGA{67} := 		
					 \begin{vmatrix} 							X_1Z_1	&	X_2Z_2\\
  																Y_1T_1	&	Y_2T_2\\
 					 \end{vmatrix} \ \			 
\OMEGA{89} := 		
					\begin{vmatrix} 							X_1T_1	&	X_2T_2\\
  																Y_1Z_1	&	Y_2Z_2\\
 					 \end{vmatrix}  \PUNKT
\end{split}
\end{equation}


\end{not1}
\begin{rmk}\label{rmk_base_points}[\textit{The 16 base points of the $(1,2,2)$-polarization}]
With the notation of Lemma \ref{fundamentallemmagenus9}, we consider the quartic curve $\Du$ in $\PR^3$ with coordinates $[x,y,z,t]=[X^2,Y^2,Z^2,T^2]$ defined by
\begin{align*}
\Du: \begin{cases}
x + y + z + t &= 0 \\
q(x,y,z,t)^2&= xyzt \PUNKT
\end{cases}
\end{align*}
\noindentTX The lines $x,y,z,t$ in the plane $H: x+y+z+t=0$ in $\PR^3$ are clearly bitangent lines on $\Du$, and for each bitangent line $l$ in this set
we can denote by $l_1 + l_2$ the corresponding degree $2$ on $\mathcal{D}$ such that:
\begin{equation*}
l.\mathcal{D} = 2(l_1 + l_2) \BEISTRICH
\end{equation*}
\noindentTX If we consider now $L_1$ and $L_2$ points of $\mathcal{C}$ in the preimage of $l_1$ and $l_2$ respectively, then it is now easily seen that $\GrpJ.[(L_1, L_2)]$ is a $\GrpJ$-orbit of base points of the linear system $\mathcal{L}$ in $A$. Indeed, by using the polynomial expression of the holomorphic canonical section of $\Surf$ in (\ref{polynomialexpressioncanonicalmap}) one can see directly that $\GrpJ.[(L_1, L_2)]$ is a $\GrpJ$-orbit of base points of the sublinear system generated by $\OMEGA{45}$, $\OMEGA{67}$, $\OMEGA{89}$ (which we defined in Notation \ref{NOTcanbdlS}). Since the set of the base points of the polarization on the $(1,2,2)$-polarized abelian threefold $A$ containing $\Surf$ has cardinality $16$, all $\GrpJ$-orbit of base points must be in the form $\GrpJ.[(L_1, L_2)]$ for some bitangent line.
\end{rmk}
\noindentTX We are now in the position to prove the second claim of Theorem \ref{injectivedifferential}.
\begin{proof}[of the Theorem \ref{injectivedifferential}, point 2)]


\noindentTX We consider the projections $q_g$ of $A$ onto the $(1,1,2)$-polarized abelian threefolds $A_g:=\bigslant{A}{\left<g\right>}$ associated to the non-trivial elements $g$ of $\GrpJ$, and Diagram \ref{diagramproj2}:
\begin{align*}
\begin{diagram} 
		A &\lInto&\Surf 					&	\rTo^{\phi_{\Surf}} 	&	\Sigma	  &\rInto& \PR^5				\\ 	
		\dTo_{q_g} && \dTo_{q_g} 	&			&\dTo &&\dDashto		\\
		A_{g}&\lInto&\TD_{g} & \rTo^{\phi_{\TD_g}}			&\Sigma_g &\rInto&\PR^3		\\
&&&\rdTo(1,1)^{} Z_g  \ruTo(1,1)&& 
 \end{diagram}
 \end{align*}
Recall that here $\TD_g$ denotes the image $q_g(\Surf)$ in $A_g$.  
If the differential $d_{z}\phi_{\Surf}$ at a point $z$ on $\Surf$ is not injective, then we have that the differential at $q_g(z)$ of the canonical map of $Z_g$ is not injective for every non-trivial element $g$ of $\GrpJ$. By applying Theorem \ref{CATSCHREYER}, the image in $\PR^3$ of $q_{g}(z)$ with respect to $\phi_{\mathcal{T}_g}$ 
must be one of the pinch points belonging to the double nodal curve $\Gamma_{g}$ in the canonical image of $Z_g$ in $\PR^3$. By Remark \ref{rmk112}, the latters are contained in the plane $\theta_{\gamma}=0$.
Hence, $z$ must be a base point of the linear system $|\mathcal{O}_A(\Surf)|$ in $A$. 

\noindentTX Thus, it is enough to prove the claim for the base points of the linear system $|\mathcal{O}_A(\Surf)|$ in $\Surf$.

\noindentTX Let us consider in particular a base point $z_0$. We have to prove that, for every tangent vector $\nu$ to $\Surf$ in $z_0$, there exists a divisor $D$ in the canonical class $|K_\Surf|$ such that $D$ contains $z_0$, but $\nu$ is not tangent to $D$ in $z_0$. Then, to conclude the proof of the proposition it suffices to prove the following lemma.
\end{proof}

\begin{lem1}
Let $b$ be a base point of the polarization $\mathcal{L}$ on a general $(1,2,2)$ abelian 3-fold $A$ with flat uniformizing coordinates $(z_1, z_2, z_3)$. Moreover, let us consider $p: A \longrightarrow \Jac$ an isogeny onto the Jacobian of a quartic plane curve $\Du$ with kernel $\GrpJ$ and
$(\theta_{\gamma})_{\gamma \in \GrpJ^*}$ a basis for $H^0(A, \mathcal{L})$ (recall that $\GrpJ^*$ is the group generated by the $2$-torsion bundles $\eta_1$ and $\eta_2$ defining the cover $p$). 

\noindentTX Then the four points in $\PR^2(\C)$ represented by the columns of following matrix are in general position:
\begin{equation*}
\begin{pmatrix} \frac{\partial}{\partial z_1}\theta_{0}(b) &  \frac{\partial}{\partial z_1}\theta_{\eta_1}(b) &  \frac{\partial}{\partial z_1}\theta_{\eta_2}(b) &  \frac{\partial}{\partial z_1}\theta_{\eta_1 + \eta_2}(b) \\
			   \frac{\partial}{\partial z_2}\theta_{0}(b) &  \frac{\partial}{\partial z_2}\theta_{\eta_1}(b) &\frac{\partial}{\partial z_2}\theta_{\eta_2}(b)  &  \frac{\partial}{\partial z_2}\theta_{\eta_1 + \eta_2}(b) \\
			   \frac{\partial}{\partial z_3}\theta_{0}(b) &  \frac{\partial}{\partial z_3}\theta_{\eta_1}(b) &\frac{\partial}{\partial z_3}\theta_{\eta_2}(b) &  \frac{\partial}{\partial z_3}\theta_{\eta_1 + \eta_2}(b)
\end{pmatrix} \PUNKT
\end{equation*} 
\end{lem1}
\begin{proof}

With the notation of Remark \ref{rmk_base_points}, to each bitangent $l \in \{x,y,z,t\}$ of $\Du$ we fix $b_{l} := [L_1,L_2]$ a representative of a $\GrpJ$-orbit of base points.
The unramified bidouble covering $p: \Cu \longrightarrow \Du$ is then defined by the $2$-torsion points in the Jacobian $\Jac(\Du)$:
\begin{align*}
\eta_1 &= \OC_{\Du}(y_1 + y_2 - x_1 - x_2) \\ 
\eta_2 &= \OC_{\Du}(z_1 + z_2 - x_1 - x_2) \\
\eta_1 \otimes \eta_2 &= \OC_{\Du}(t_1 + t_2 - x_1 - x_2)  \PUNKT
\end{align*}
With this notation, it follows that, as projective points,
\begin{align*}
[\GRAD \theta_0(b_y)] &= [\GRAD \theta_0(b_x + \eta_1)] = [\phi_{\eta_1}(b_x) \cdot \GRAD \theta_{\eta_1}(b_x)]
= [ \GRAD \theta_{\eta_1}(b_x)]\\
[\GRAD \theta_0(b_z)] &= [\GRAD \theta_0(b_x + \eta_2)] = [\phi_{\eta_2}(b_x) \cdot \GRAD \theta_{\eta_2}(b_x)]
= [\GRAD \theta_{\eta_2}(b_x)] \\
[\GRAD \theta_0(b_t)] &=[\GRAD \theta_0(b_x + \eta_1 + \eta_2)] = [\phi_{\eta_1+\eta_2}(b_x) \cdot \GRAD \theta_{\eta_1 + \eta_2}(b_x)] 
= [\GRAD \theta_{\eta_1 + \eta_2}(b_x)]
\BEISTRICH
\end{align*}
where $\phi_{\gamma}(b_x)$ is the non-zero automorphy factor in Equation \ref{thetadef1}, which defines the section $\theta_{\gamma}$.
This proves the claim, since the projective points $[\GRAD \theta_0(b_l)]$ represent image in $\PR^2$ of the divisors $l_1 + l_2$ on $\Du$ with respect to the Gauss map of $\Theta$. The Gauss map associates to these divisors the coordinates in $\PR^{2\vee}$ of the bitangents $l$ passing through them. Hence, since the four bitangents are in general position, also their projective coordinates representation must be in general position.
\end{proof}

\newcommand{\IM}[1]{\Im \text{m} \ #1}
\newcommand{\ct}[1]{\frac{\partial \y{#1}}{\partial \x{#1}}}

\section{Polarizations of type $(1,2)$ and $(2,4)$ on abelian surfaces}

When we consider the case of a $(1,2,2)$-polarized abelian threefold which is isogenous to a product of an abelian surface and an elliptic curve, we can easily see that each holomorphic section of the polarization also splits into a product. For this reason, in such a case, the behaviour of the canonical map is strictly related to the properties of polarizations of type $(1,2)$ and $(2,4)$ on abelian surfaces.
These polarization types have been explicitly described by W. Barth (see \cite{Barth}) who was interested in the study of the quadratic equations of abelian surfaces embedded in $\PR^7$ via a polarization of type $(2,4)$.
In the first part of this section, we give a short summary about some remarkable properties of polarizations of type $(1,2)$ which will be of importance in the next section.
In the second part, we will focus on polarizations of type $(2,4)$ with a different approach from the one followed by Barth in his work. 

Given a $(2,2)$-polarized surface $(\B, \PN)$ with an isogeny $\pi \colon \B \longrightarrow \SA$ of degree $2$ on an abelian surface $\SA$, we consider the vector bundle ${\pi_{\B}}_*(\PN)$, which splits into the direct sum of its invariant and anti-invariant parts, denoted by $\PMp$ and $\PMm$. The vector bundle $\PMp \otimes \PMm$ is of type $(2,4)$, and we will show that the multiplication map $H^0(\SA, \PMp) \otimes H^0(\SA, \PMm) \longrightarrow H^0(\SA, \PMp \otimes \PMm)$ and then a skew-symmetric map
\begin{equation*}
H^0(\SA, \PMp) \land H^0(\SA, \PMm) \longrightarrow H^0(\SA, \PMp \otimes \PMm)
\end{equation*}
defined by a derivation are both injective, when $\SA$ is not isogenous to a product.

\begin{def1}[\textit{(Even and odd $2$-torsion points)}]\label{twotorsS}
Let us assume that $(\SA, \mathcal{L})$ is a general $(1,2)$-polarized abelian surface, with $\mathcal{L}$ of characteristic $0$. Considered $\mathcal{D}$ a symmetric divisor in the linear system $|\mathcal{L}|$, the set $S[2]$ of the $2$-torsion points on $S$ is the union of the sets
\begin{align*}
\SA[2]^{+}(\mathcal{D}) &:= \{x \in \SA[2] \  | \ mult_{x}\mathcal{D}\  \text{is even} \} \\
\SA[2]^{-}(\mathcal{D}) &:= \{x \in \SA[2] \  | \ mult_{x}\mathcal{D}\  \text{is odd} \}  \PUNKT
\end{align*}
The order of these sets do not depend on $\mathcal{D}$, and we have that (refer also to \cite[p. 97]{BLange}):
\begin{align*}
\#\SA[2]^{+} &:= \SA[2]^{+}(\mathcal{D}) = 12 \\
\#\SA[2]^{-} &:= \SA[2]^{-}(\mathcal{D}) = 4 \PUNKT
\end{align*}

\noindentTX For simplicity, we say that a $2$-torsion point is even (resp. odd) if it belongs to $\SA[2]^{+}$ (resp. $\SA[2]^{-}$). Since a basis for $H^0(\SA, \mathcal{L})$ is given by four even theta functions which vanish on the set of odd $2$-torsion points, we have that the set $\mathcal{B}$ of base points of $\mathcal{L}$ is $\SA[2]^{-}$.
\end{def1}

\begin{oss1} [\textit{(Singular curves in a $(1,2)$ polarization)}]\label{osservationsingularfibers}
The rational map $\phi_{\mathcal{L}}: \SA \dashrightarrow \PR^1$ can be extended to a morphism $\psi: Blow_{\mathcal{B}}(\SA) \longrightarrow \PR^1$, 
where $Blow_{\mathcal{B}}(\SA)$ denotes the blow-up of $\SA$ at the points of $\mathcal{B}$. 
The morphism $\psi$ is a fibration, and its general fiber is a smooth non-hyperelliptic curve of genus $3$. 
By applying the Zeuthen Segre's Formula, we can determine the number $\mu$ of singular fibers of $\psi$, which is:
\begin{equation*} 
\mu = 8 + \chi_{top}(Blow_{\mathcal{B}}(\SA)) = 8 + 4 = 12 \PUNKT
\end{equation*}
 
By Definition \ref{twotorsS} it follows that:
\begin{itemize}
\item[a)] No divisor $\mathcal{D}$ of $|\mathcal{L}|$ is singular at a base point.
\item[b)] If a divisor $\mathcal{D}$ of $|\mathcal{L}|$ contains a $2$-torsion point $x$ which is not contained in $\mathcal{B}$, 
then $\mathcal{D}$ is singular in $x$. Moreover, if a divisor $\mathcal{D}$ of $|\mathcal{L}|$ contains two distinct points of 
$\SA[2]^{+}(\mathcal{D})$, then $\mathcal{D}$ is reducible, or $S$ is isomorphic to a polarized product of elliptic curves. Indeed, 
let us suppose $x$ and $y$ are distinct $2$-torsion points for which $mult_{x}(\mathcal{D}) = mult_{y}(\mathcal{D}) = 2$. 
If $\mathcal{D}$ were irreducible, then its normalization $\mathcal{C}$ would be a smooth elliptic curve, and we would have a homomorphism 
from $\mathcal{C}$ to $S$, which means that $\SA$ is isomorphic to a polarized product of elliptic curves. \newline
\end{itemize}

\noindentTX We can now describe all possible configurations of reducible fibers of $\psi$. Let us consider a reducible fiber $\mathcal{D}$ of $\psi$, 
 which we write as the union of its reducible components
 \begin{equation} \label{reduciblefiberpsi}
 \mathcal{D} = E_1 + \cdots + E_s \PUNKT
 \end{equation}
We have clearly, by adjunction, that
\begin{equation*} 
2p_a(E_j) - 2 = E_j^2	\PUNKT
\end{equation*}
On the other side, none of the divisors $E_s$ in the decomposition in (\ref{reduciblefiberpsi}) can be ample: the polarization on $S$, which is of type $(1,2)$, would be in this case the tensor product of two polarizations. Thus, we can conclude that $E_j^2 = 0$ and $E_j$ is a curve of genus $1$ for every $j$. On the other side, we have that
\begin{align*}
3 = p_a(D) &= \sum_{j}p_a(E_j) + \sum_{i \neq j} E_i.E_j - s + 1 \\
&=  \sum_{i \neq j} E_i.E_j + 1 \PUNKT
\end{align*}
Furthermore, because $\mathcal{D}$ is connected, the only possible configurations are the following:
\begin{itemize}
\item[a)] $\mathcal{D} = E_1 + E_2 + E_3$ with $E_1.E_3 = E_2.E_3 = 1$ and $E_1.E_2=0$.
\item[b)] $\mathcal{D} = E_1+E_2$ and $E_1.E_2 = 2$. \newline
\end{itemize}

\noindentTX Note that, in both cases, each irreducible component $E_i$ is a smooth elliptic curve. 

In the case a), we have that $E_1^2 = 0$ and $E_1.\mathcal{D} = 1$, and it follows  that (see also \cite[cf. Lemma 10.4.6]{BLange}):
\begin{equation*}
(\SA, \mathcal{L}) \cong (E_1,\mathcal{O}_{E_1}(O_{E_1}))  \boxtimes (E_3,\mathcal{O}_{E_3}(2O_{E_3})) \cong (E_2,\mathcal{O}_{E_2}(O_{E_2}))  \boxtimes (E_3,\mathcal{O}_{E_3}(2O_{E_3})) \PUNKT
 \end{equation*}
 This means that $\SA$ itself is a product of (polarized) elliptic curves.
 
 Finally, in the case b), we can consider the difference morphism $\phi \colon E_1 \times E_2 \longrightarrow A$ defined by $\phi(p,q) := p-q$. The kernel of $\phi$ consists of the two points in which $E_1$ and $E_2$ intersect. Thus, $\phi$ is an isogeny, and we conclude that
\begin{equation*}
\phi^*(\mathcal{L}) = (\mathcal{D}.E_1, \mathcal{D}.E_2) = (2,2) \PUNKT
\end{equation*}
In conclusion, we have in this case that $S$ is isogenous to the product of two $(2)$-polarized elliptic curves, which carries a natural polarization of type $(2,2)$.
\end{oss1}

\begin{not1} \label{Surftwotwo}
\noindentTX Let us consider a $(2,2)$-abelian surface $\B := \bigslant{\C^2}{\tauB \Z^2 \oplus 2\Z^2}$, with $\tauB$ a matrix in the two-dimensional Siegel half space $\SiegH_2$. 
We denote by $\PN$ the ample line bundle of characteristic $0$ of type $(2,2)$ on $\B$. As we did in the definition in (\ref{THETAGAMMA}), we shall denote the holomorphic sections of the polarizations on $\B$ by
\begin{align} \label{BASISSECTIONPN}
H^0(\B, \PN) = \left<\theta_{00}, \theta_{01}, \theta_{10}, \theta_{11}\right> \BEISTRICH
\end{align}
where $\theta_{uv}$ is invariant (resp. anti-invariant) under the translation by the $2$-torsion point $ae_1 + be_2$ precisely when $ua + vb = 0$ (resp. $ua + vb = 1$). In particular, if we consider the quotient by $2$-torsion element $e_1+e_2$ we have a quotient map
\begin{equation*} 
\prBSA \colon \B \longrightarrow \SA :=\bigslant{\B}{\left<e_1+e_2\right>}
\end{equation*}
onto a $(1,2)$-polarized surface $\SA$. We say in the introduction that $\SA$ the rank $2$ vector bundle ${\prBSA}_*(\PN)$ on $\SA$ splits into a sum of two ample line bundles $\PMp$ and $\PMm$, which are respectively the invariant part and the anti-invariant part of ${\prBSA}_*(\PN)$ with respect to the translation by $e_1+e_2$. The line bundles $\PMp$ and $\PMm$ are algebraic equivalent (i.e. they are equivalent up to translation) and they define a $(1,2)$-polarization on $\SA$. According to (\ref{BASISSECTIONPN}) we have:
\begin{equation*} 
\begin{split}
H^0(\SA, \PMp) &= \left<\theta_{00}, \theta_{11} \right> \\
H^0(\SA, \PMm) &= \left<\theta_{10}, \theta_{01} \right>  \PUNKT
\end{split}
\end{equation*}
\end{not1}
The line bundle $\PMp \otimes \PMm$ is of type $(2,4)$, hence its space of global holomorphic sections has dimension $8$ and it is generated by $4$ odd and $4$ even holomorphic sections, as we are going to show.


\begin{def1}\label{ANTISIMMsections} 
Given an abelian variety $A = \bigslant{\C^g}{\Lambda}$ with $\mathcal{L}$ and $\mathcal{M}$ two ample line bundles, there is a natural multiplication map on the space of holomorphic sections:
\begin{equation*}
\mu \colon H^0(A, \mathcal{L}) \otimes H^0(A, \mathcal{M}) \longrightarrow H^0(A, \mathcal{L} \otimes \mathcal{M}) \BEISTRICH
\end{equation*}
and every derivation $\DER{v}{}$ given by a vector $v$ in the vector space $\C^2 \cong H^0(A, T_A)$, defines a skew-symmetric multiplication map
\begin{equation*}
\SKMULT{v} \colon H^0(A, \mathcal{L}) \land H^0(A, \mathcal{M}) \longrightarrow H^0(A, \mathcal{L} \otimes \mathcal{M})
\end{equation*}
as follows:
\begin{align*}
\SKMULT{v}(s \wedge t) &= \begin{vmatrix}
                              s & t \\
                              \DER{v}{s} & \DER{v}{t}
                             \end{vmatrix} = s \cdot \DER{v}{t} - \DER{v}{s} \cdot t \PUNKT
\end{align*}
It is easily seen that this map is well defined and $\SKMULT{v}(s \wedge t)$ is a holomorphic section of $H^0(A, \mathcal{L} \otimes \mathcal{M})$. To show this, let us assume that $[\{\phi_{\lambda}\}_{\lambda}]$ and $[\{\psi_{\lambda}\}_{\lambda}]$ are the factor of automorphy of the ample line bundles $\mathcal{L}$ and $\mathcal{M}$ respectively. Let us consider two global holomorphic sections $s$ and $t$ of $\mathcal{L}$ and $\mathcal{M}$. Then we have, for every vector $z$ in $\C^g$ and every element $\lambda$ of the lattice $\Lambda$, that:
\begin{align*}
        \begin{vmatrix}
        s(z + \lambda)  & t(z + \lambda) \\
        \DER{v}{s}(z + \lambda)& \DER{v}{t}(z + \lambda)
        \end{vmatrix} &= 
        \begin{vmatrix}
        \phi_{\lambda}(z)s(z) & \psi_{\lambda}(z)t(z) \\
        \DER{v}{\phi_{\lambda}}(z)s(z) + \phi_{\lambda}(z)\DER{v}{s}(z)& \DER{v}{\psi_{\lambda}}(z)t(z) + \psi_{\lambda}(z)\DER{v}{t}(z)
        \end{vmatrix} \\
        &=\begin{vmatrix}
        \phi_{\lambda}(z)s(z) & \psi_{\lambda}(z)t(z) \\
         \phi_{\lambda}(z)\DER{v}{s}(z)& \psi_{\lambda}(z)\DER{v}{t}(z)
        \end{vmatrix} \\
        &=\phi_{\lambda}(z)\psi_{\lambda}(z)\begin{vmatrix}
        s(z) & t(z) \\
        \DER{v}{s}(z)&\DER{v}{t}(z)
        \end{vmatrix}        
\end{align*}
This proves the claim, since $[\{\phi_{\lambda}\cdot \psi_{\lambda}\}_{\lambda}]$ is the factor of automorphy of the line bundle $\mathcal{L} \otimes \mathcal{M}$.
\end{def1}

\begin{ex1}\label{basicexampleEll}
Every $(2)$-polarized elliptic curve $\E = \bigslant{\C}{\left<\tau \Z \oplus 2\Z \right>}$ can be viewed as the Riemann surface defined in an affine neighborhood $\EU{}$ with coordinates ($\x{}$,$\y{}$) by an equation $\Eeq{} = 0$, where:
\begin{equation}\label{LegendreNormalForm}
\Eeq{} := \y{}^2 - (\x{}^2 - 1)(\x{}^2 - \pa{}^2) \BEISTRICH
\end{equation}
and parameter $\pa{}$ in Equation \ref{LegendreNormalForm} only depends on the element on $\tau$. 
\noindentTX The affine model in Equation (\ref{LegendreNormalForm}) is often referred to as the Legendre normal form of $\E_{}$ (see \cite{Catanese2016}). \newline

\noindentTX The Riemann surface defined by the equation $\Eeq{} = 0$ has two points at infinity, which we shall denote by $\inftyplus$ and $\inftyminus$, at which the function $\x{}$ has a simple pole. Thus, we can consider $\vinf{} := \frac{1}{\x{}}$ to be a local parameter at these points in an affine neighborhood $\EV{}$ with coordinates $(\vinf{}, \winf{})$, where $\E$ is defined by the equation
\begin{equation}\label{LegendreNormalForm2}
\Eeqinf{} := \winf{}^2 - (1 - \vinf{}^2)(1- \pa{}^2\vinf{}^{2}) \PUNKT
\end{equation}
The change of coordinates between $\EU{}$ and $\EV{}$ is the map $(\x{}, \y{}) \mapsto (\vinf{}, \winf{}) = (\x{}^{-1},  \y{}\x{}^{-2})$. \newline

\noindentTX The function $\x{}$ defines a double cover of $\PR^1$ of degree $2$, which is ramified over the points of the set $\{1, -1, \pa{}, -\pa{}\}$. Such a double cover $\E \longrightarrow \PR^1$ has been called by Inoue a Legendre function for $\E$ (see \cite[p. 60]{BCatFrapp} and \cite{INOUE})). In general, a Legendre function $\LF$ on an elliptic curve $\E = \bigslant{\C}{\tau \Z \oplus 2 \Z}$ is a meromorphic function on $\E$ which satisfies, for every point $z$ on $\E$, the following properties:
\begin{itemize}
\item $\LF(-z) = \LF(z)$.
\item $\LF \left(z+1\right) = -\LF(z)$.
\item $\LF \left(z+\frac{\tau}{2}\right) = \frac{\pa{}}{\LF{(z)}}$, $\LF(\frac{1}{2}) = -1$, $\LF(0) = 1$, $\LF(\frac{\tau}{2}) = \pa{}$, $\LF \left(\frac{1+\tau}{2}\right) = -\pa{}$.
\item  $\LF'(z) = 0$ if and only if $z \in \{0, 1, \frac{\tau}{2}, 1 + \frac{\tau}{2}  \}$. 
\end{itemize}
Denoted by $\theta_0$ and $\theta_1$ the canonical generators of the space $H^0(\E, \mathcal{O}_{\E}(2 \cdot O_{\E}))$, since a Legendre function on an elliptic curve is unique, (see \cite[Lemma 3.2]{INOUE}), we have that:
\begin{equation} \label{LegendreRelationsEQ}
\frac{\theta_0(0)\theta_1(z)}{\theta_1(0)\theta_0(z)} = \x{} \PUNKT
\end{equation}
Indeed, $\x{}$ behaves as a Legendre function, while the left member satisfies all properties which characterize Legendre functions.

\noindentTX Note that the algebraic function $\y{}$ corresponds in the analitic setting to the derivative of the Legendre function (up to a constant), which vanishes precisely on the set of $2$-torsion points of $\E_j$, and we have:

\begin{equation} \label{etadefinition_Y}
\y{} \sim \eta := \begin{vmatrix}\theta_0 & \theta_1 \\ 
\frac{\partial}{\partial z} \theta_0 & \frac{\partial}{\partial z} \theta_1 \end{vmatrix} = \SKMULT{z}(\theta_0 \wedge \theta_1) \in H^0(\E, \OC_{\E}(4\cdot O_{\E}))  \PUNKT
\end{equation}

This shows that $\y{}$ corresponds, in the analytic setting, precisely to $\SKMULT{z}(\theta_0 \wedge \theta_1)$.
\end{ex1}

According to Definition \ref{ANTISIMMsections}, every tangent vector $\SKMULT{v}$ on $\SA$ induces a well defined map:
\begin{equation*}
\SKMULT{v} \colon H^0(\SA, \PMp) \land H^0(\SA, \PMm) \longrightarrow H^0(\SA, \PMp \otimes \PMm) \PUNKT
\end{equation*}
 It is important to notice that the image of the multiplication map $\mu \colon H^0(\SA, \PMp) \otimes H^0(\SA, \PMm) \longrightarrow H^0(\SA, \PMp \otimes \PMm)$ is generated by even sections, while the image of $\SKMULT{v}$ only contains odd holomorphic sections. Actually, the vector space $H^0(S, \PMp \otimes \PMm)$ splits into the direct sum of these images, as we are going to prove.

\begin{prop1}\label{TWOFOURTH}
Let $(\B, \PN)$ a $(2,2)$-polarized abelian surface, with $\prBSA \colon \B \longrightarrow \SA$ an isogeny of degree $2$, where $\SA$ is a simple $(1,2)$-polarized abelian surface. Considered the natural splitting of the rank $2$ vector bundle 
\begin{equation*}
{\prBSA}_*(\PN) = \PMp \oplus \PMm
\end{equation*}
into its invariant and anti-invariant parts, we have that:
\begin{itemize}
 \item[1)] For every tangent vector $\DER{v}{}$, the multiplication map $\mu \colon H^0(\SA, \PMp) \otimes H^0(\SA, \PMm) \longrightarrow H^0(\SA, \PMp \otimes \PMm)$ and the map
 \begin{equation}\label{DIFFMAPV}
\SKMULT{v} \colon H^0(\SA, \PMp) \land H^0(\SA, \PMm) \longrightarrow H^0(\SA, \PMp \otimes \PMm)
\end{equation}
 are injective. Denoted by $U^+$ and $U^-$ the respective $4$-dimensional images in $H^0(\SA, \PMp \otimes \PMm)$, we have the following splitting:
 \begin{equation*}
H^0(\SA, \PMp \otimes \PMm) = U^+ \oplus U^-
 \end{equation*}
 \item[2)] The rational map 
 \begin{equation*}
  \phi \colon \SA \dashrightarrow \PR^3 \times \PR^3
 \end{equation*}
 induced by the sublinear systems $|U^+|$ and $|U^-|$ is invariant under the action of the involution $(-1)$ on $\SA$, and it is generically finite of degree $2$ onto its image.
\end{itemize}
\end{prop1}
\begin{proof}
Since $\SA$ is not isogenous to a product of elliptic curves, $\B$ is not a product, and the linear system $(2,2)$ embeds the Kummer variety $K(\B) = \bigslant{\B}{\left<\pm 1 \right>}$ as a quartic in $\PR^3$. Hence, the multiplication map $\mu$ on $H^0(\SA, \PMp) \otimes H^0(\SA, \PMm)$ is injective, since otherwise the Kummer variety $K(\B)$ would be contained in a quadric.

To prove that every derivation $\DER{v}{}$ induces an injective map on $H^0(\SA, \PMp) \otimes H^0(\SA, \PMm)$ as in Definition \ref{ANTISIMMsections}, it is enough to prove it for a basis of derivations $\bigl\{\DER{z_1}{}, \DER{z_2}{} \bigr\}$. \newline 

\noindentTX Once fixed a derivation $\DER{v}{}$, to simplify the notation we will denote by $\DTH{ij}$ the derivative $\DER{v}{\theta_{ij}}$. In general we denote, moreover:
Let us denote by:
\begin{align*}
Z_1 &:= \SKMULT{v}(\theta_{00} \wedge \theta_{01}) = \theta_{00}\DTH{01} - \DTH{00}\theta_{01} \\
Z_2 &:= \SKMULT{v}(\theta_{00} \wedge \theta_{10}) = \theta_{00}\DTH{10} - \DTH{00}\theta_{10} \\
Z_3 &:= \SKMULT{v}(\theta_{11} \wedge \theta_{01}) = \theta_{11}\DTH{01} - \DTH{11}\theta_{01} \\
Z_4 &:= \SKMULT{v}(\theta_{11} \wedge \theta_{10}) = \theta_{11}\DTH{10} - \DTH{11}\theta_{10}
\end{align*}
Since $Z_1$ and $Z_4$ are anti-invariant for the action of the involution $z \mapsto z + (0,1)$ on $\SA$, while $Z_2$ and $Z_3$ are always invariant sections, the claim 1) follows once we prove that $div(Z_1) \neq div(Z_4)$ and $div(Z_2) \neq div(Z_3)$. \newline

\noindent Since the latter conditions are open in the moduli space, we can see what happens when $\B$ is the product of elliptic curves $\E_1 \times \E_2$. On each factor $\E_i$ we use local coordinates $(\x{i}, \y{i})$ as in Example (\ref{basicexampleEll}) with equations $\y{i}^2 - (\x{i}^2 - 1)(\x{i}^2 - \pa{i}^2) = 0$ as in (\ref{LegendreNormalForm}). If we consider now the tangent vector $\DER{z_1}{}$ then $Z_1 = 0 = Z_4$, but $div(Z_2)$ coincides with $div(y_1)$ and $div(Z_3)$ with $div(y_1) + 2div(x_2)$ on the affine open set of $\B$ where $\theta_{00}$ does not vanish. Hence, we have also that $div(Z_2) \neq div(Z_3)$ in the general case. We can argue similarly by using the derivation $\DER{z_2}{}$, and we get that $div(Z_1) \neq div(Z_4)$. The claim 1) follows. \noindent 

\noindentTX We prove now the second claim. The rational map $\phi \colon \SA \dashrightarrow \PR^3 \times \PR^3$ can be extended to a morphism
\begin{equation*}
\widehat{\phi} \colon X \longrightarrow \PR^3 \times \PR^3 \BEISTRICH
\end{equation*}
where $X$ denotes the blow-up of $\SA$ at the points of $\mathcal{B}(\PMp) \cup \mathcal{B}(\PMm)$. Indeed, the second component of $\phi$ is defined by the composition of the map $\SA \dashrightarrow |\PMp| \times |\PMm| \cong \PR^1 \times \PR^1$ with the Segre embedding. 



Moreover, its extension is finite of degree $4$ onto its image, and $\mathcal{B}(\PMp) \cup \mathcal{B}(\PMm)$ is stable under the involution $x \mapsto x + (0,1)$ 
Hence, it sufficies to show that, if $x$ and $y = x + (0,1)$ are base points of $|\PMp|$ on $\SA$, then $\widehat{\phi}_1(x) \neq \widehat{\phi}_1(y)$ (here $\widehat{\phi}_1$ denotes the first component of $\widehat{\phi}$). \newline

\noindentTX We have that:

\begin{align*}
\widehat{\phi}_1(y) = \begin{bmatrix}
                        -\DTH{00}(y)\theta_{01}(y) \\
                        -\DTH{00}(y)\theta_{10}(y) \\
                        -\DTH{11}(y)\theta_{01}(y) \\
                        -\DTH{11}(y)\theta_{10}(y)
                      \end{bmatrix} = 
                      \begin{bmatrix}
                        \DTH{00}(x)\theta_{01}(x) \\
                        -\DTH{00}(x)\theta_{10}(x) \\
                        \DTH{11}(x)\theta_{01}(x) \\
                        -\DTH{11}(x)\theta_{10}(x)
                      \end{bmatrix} \PUNKT
\end{align*}
Since $\DTH{00}(x)$ and $\DTH{11}(x)$ can not simultaneously vanish (because otherwise the divisors $div(\theta_{00})$ and $div(\theta_{11})$ would not transversally intersect at $x$), the equality of the images $\widehat{\phi}_1(x) = \widehat{\phi}_1(y)$ would imply that both $x$ and $y$ also belong to $div(\theta_{01})$ or $div(\theta_{10})$. But on the other hand, since the $(2,2)$-polarization on $\B$ embeds the Kummer variety of $\B$, we have that also on $\SA$ we have:
\begin{align*}
 div(\theta_{00}) \cap div(\theta_{11}) \cap div(\theta_{10}) &= \emptyset \\
 div(\theta_{00}) \cap div(\theta_{11}) \cap div(\theta_{01}) &= \emptyset
\end{align*}
This forces $\widehat{\phi}_1(x) \neq \widehat{\phi}_1(y)$, and the second claim follows. \newline





\end{proof}



\section{Surfaces in a $(1,2,2)$-polarization which is isogenous to a product} \label{degenerationspolarizationsection}
\noindentTX 
In Section \ref{basicsectionssect}, we investigated the behavior of the canonical map of an ètale bidouble cover $\Surf$ of type $(1,2,2)$ of a general principal polarization, and in Proposition \ref{injectivedifferential}) we proved that the canonical map has everywhere injective differential and it is injective outside of some canonical curves, on which its restriction has degree $2$. However, Remark \ref{CHOICENOTGENERAL} shows that such surfaces are quite special in their linear system, and the question whether the canonical map in the general case is an embedding still remains open. \\

\noindent In this section, we look at this question by considering sufficiently general surfaces in a polarization of type $(1, 2, 2)$ on an abelian threefold $A$, which is isogenous to a polarized product of a $(2, 2)$-polarized surface $B$ and a $(2)$-polarized elliptic curve $E$. To achieve this goal we first need to consider the special case in which $B$ itself is a product of elliptic curves. 
We begin this section by introducing some notation.


\begin{not1} \label{productSurfEll}
\noindentTX Let us consider a point $\tau = (\tauB, \tauE)$ in the product $\mathcal{H}_2 \times \mathcal{H}_1$. 
The $(2,2,2)$-polarized abelian threefold $T$ defined by $\tau$ is then the product of a $(2,2)$-abelian surface $\B := \bigslant{\C^2}{\tauB \Z^2 \oplus 2\Z^2}$ and a $(2)$-polarized elliptic curve $\E := \bigslant{\C}{\tauE \Z \oplus 2\Z}$.
We denote by $\PN$ the ample line bundle on $\B$ of characteristic $0$ of type $(2,2)$. 

As in Notation \ref{Surftwotwo} in the previous section, it is defined a projection of $\B$
\begin{equation*} 
\prBSA \colon \B \longrightarrow \SA :=\bigslant{\B}{\left<e_1+e_2\right>}
\end{equation*}
onto a $(1,2)$-polarized surface $\SA$, and we have the splitting ${\prBSA}_*(\PN) = \PMp \oplus \PMm$ into invariant and anti-invariant part. \newline 


\noindentTX We can now consider the $(1,2,2)$-polarized abelian threefold $(A, \mathcal{L})$ obtained as the quotient of $T$ by the translation $e_1 + e_2 + e_3$:
\begin{equation*} 
\prTA \colon T \longrightarrow A:= \bigslant{T}{\left<e_1+e_2+e_3\right>} \PUNKT
\end{equation*}
The vector space of global holomorphic sections of $\mathcal{L}$ is generated by those global sections of the polarization on $T$ which are invariant under the action by $e_1 + e_2 + e_3$, and we have the following decomposition: 
\begin{equation} \label{sectionsDECOMP}
 \begin{split}
H^0(A, \mathcal{L}) &= H^0(T, \PN \boxtimes \mathcal{O}_{\E}(2\cdot O_{\E}))^{e_1 + e_2 + e_3} \\
&= H^0(\SA \times \E, \prBSA_* \PN \boxtimes \mathcal{O}_{\E} (2\cdot O_{\E}))^{e_1 + e_2 + e_3} \\
&\cong H^0(\SA, \PMp)\cdot \theta^{\E}_0 \oplus H^0(\SA, \PMm)\cdot \theta^{\E}_1
 \end{split}
\end{equation}
By (\ref{sectionsDECOMP}) it also follows a similar decomposition for the base locus of $\mathcal{L}$:
\begin{equation} \label{BASEL1}
 \mathcal{B}(\mathcal{L}) = \mathcal{B}(\PMp) \times div(\theta^{\E}_1) \cup \mathcal{B}(\PMm) \times div(\theta^{\E}_0)
\end{equation}
\end{not1}



The following lemma shows that the images of two distinct base points of $\mathcal{L}$ with respect to the canonical map of a general surface in $|\mathcal{L}|$ are different. This happens also in the case in which our $(1,2,2)$ abelian variety $A$, as in Notation \ref{productSurfEll}, is the quotient of three elliptic curves.
\begin{lem1} \label{generalbehavior0} 
Let us consider a $(1,2,2)$-polarized abelian variety $(A, \mathcal{L})$ defined as the quotient of the product of three $(2)$-polarized elliptic curves $\E_j := \bigslant{\C}{\tau_{jj}\Z \oplus 2\cdot \Z}$ by the involution $e_1 + e_2 + e_3$. Then the canonical map of a general surface $\Surf$ in the polarization $|\mathcal{L}|$ on $A$ in injective on the set $\mathcal{B}(\mathcal{L}) \subseteq \Surf$ of base points of $|\mathcal{L}|$. 
\end{lem1}
\begin{proof}
According to the splitting of $H^0(A, \mathcal{L})$ in (\ref{sectionsDECOMP}), the vector space $H^0(A, \mathcal{L})$ is generated by four elements which we shall denote by
\begin{equation} \label{THETAELLP}
 \theta_{ijk}:= \theta^{\E_1}_{i}\theta^{\E_2}_{j}\theta^{\E_3}_{k} \BEISTRICH
\end{equation}
with $(ijk)$ a triplet in the set $\{(000),(011),(101),(110)\}$.

\noindentTX The subgroup $\mathcal{G}$ of $A$ generated by $e_1$ and $e_2$ acts on the base locus $\mathcal{B}(\mathcal{L})$ by translation, and splits $\mathcal{B}(\mathcal{L})$ as the union of four $\mathcal{G}$-orbits, each consisting of four points. We shall denote them as $\mathcal{B}_{ijk}$ for $(ijk)$ a triplet as above:
\begin{equation} \label{listbasepoints122}
 \mathcal{B}_{ijk} = \mathcal{G} \text{\ \textbf{.}}\left(\frac{1 + (i-1)\tau_{11}}{2}, \frac{1 + (j-1)\tau_{22}}{2}, \frac{1 + (k-1)\tau_{33}}{2} \right)
\end{equation}
The projective points in $\PR^5$ which represent the canonical images of the $\mathcal{G}$-orbits of base points only depend on their image with respect to the Gauss map of $\Surf$.
on which $\mathcal{G}$ acts by changing the signs of the coordinates. Hence, it suffices to show that, with $P_{ijk}$ representatives of the $\mathcal{G}$-orbits of base points, the columns of the matrix $[\GRAD \theta (P_{ijk})]_{ijk}$ represent projective points in $\PR^3$ in general position. An easy computation shows that, if $\theta = \theta_{000} + \cb \theta_{011} + \cc \theta_{101} + \cd \theta_{110}$, for certain complex coefficients $\cb$, $\cc$ and $\cd$, then we have: 
\begin{equation*}
(\GRAD(\theta)(P_{ijk}))_{ijk} = \begin{bmatrix}
\pm \cb \partial_{z_0} \theta_{011}(P_{000}) &     \partial_{z_0} \theta_{000}(P_{011}) & \pm \cd \partial_{z_0} \theta_{110}(P_{101}) & \pm \cc \partial_{z_0} \theta_{101}(P_{110}) \\
\pm \cc \partial_{z_1} \theta_{101}(P_{000}) & \pm \cd \partial_{z_1} \theta_{110}(P_{011}) & \partial_{z_1} \theta_{000}(P_{101}) & \pm \cb \partial_{z_1} \theta_{011}(P_{110}) \\
\pm \cd \partial_{z_2} \theta_{110}(P_{000}) & \pm \cc \partial_{z_2} \theta_{101}(P_{011}) & \pm \cb \partial_{z_2} \theta_{011}(P_{101}) & \partial_{z_2} \theta_{000}(P_{110})
                                \end{bmatrix}
\end{equation*}
Here the signs $\pm$ refer to the action of $\mathcal{G}$ on the coordinates. One can see that after the specialization in some particular cases (for instance to the case in which $(\cb, \cc, \cd) = (1,0,0)$), one can see that the columns are in general position as claimed. This shows that this holds true also for the general choice of the coefficients $b_{ij}$.
\end{proof}

When the following generality condition is satisfied,

\begin{prop1}\label{generalbehavior1} 
Let $(A, \mathcal{L})$ be a $(1,2,2)$-polarized abelian variety which is isogenous to the product of three elliptic curves with flat uniformizing coordinates $(z_1, z_2, z_3)$. Let us moreover assume that $\mathcal{L}$ is of characteristic $0$ on $A$, and $\Surf := div(\theta)$ be a general surface in the polarization class $|\mathcal{L}|$.

\noindentTX Then we have the following: 
\begin{itemize}
 \item For every $i=1,2,3$, the involution $\iv{i}\colon z_i \mapsto -z_i$ on $A$ which changes the sign only to the $i$-th coordinate $z_i$ leaves $\Surf$ invariant and induces an involution on every canonical divisor of the form $W_j := div \left( \frac{\partial \theta}{\partial z_j} \right)$ on $\Surf$.
 
 \item The canonical map of $\Surf$ is injective outside of the divisors $W_j$, on which the canonical map has degree $2$ and factors respectively through the involution $\iv{j}$.
 
 \item The differential of the canonical map is injective outside of the divisors $W_j$. For each $j$, the canonical map is not injective on the fixed points set of the involution $\iv{j}$ in $W_j$. Moreover, this latter set contains precisely $8$ points. 
\end{itemize}

\end{prop1}
\begin{proof}
We write all global sections of $\mathcal{L}$, 
as in Equation \ref{THETAELLP} in Lemma \ref{generalbehavior0}, as products of the form:
\begin{equation*}
 \theta_{ijk}= \theta^{\E_1}_{i}\theta^{\E_2}_{j}\theta^{\E_3}_{k} \PUNKT
\end{equation*}
Since the factors in the latter product are even theta functions, it is straightforward to verify that $\Surf$ is invariant under every involution which only changes the sign to some of the coordinates $z_j$. In particular, the derivative $\left( \frac{\partial \theta}{\partial z_j} \right)$ is anti-invariant with respect to $\iv{i}$ precisely when $i=j$, invariant otherwise. This also shows that the restriction to $W_j$ of the canonical map of $\Surf$ factors through the involution $\iv{j}$, being the derivative $\frac{\partial \theta}{\partial z_j}$ the only anti-invariant global holomorphic section of the canonical bundle on $\Surf$. This proves the first claim. \newline

\noindentTX Let us prove then the second claim. We start by considering two points $P$ and $Q$ on $\Surf$ such that $\CAN{\Surf}(P) = \CAN{\Surf}(Q)$, where $\CAN{\Surf}$ denotes the canonical map of $\Surf$. By Lemma \ref{generalbehavior0} and by our generality assumption we can assume that $P$ and $Q$ are not in the base locus $\mathcal{B}(\mathcal{L}) \subseteq \Surf$. Recall that the base locus $\mathcal{B}(\mathcal{L})$ is, according with the decomposition in (\ref{BASEL1}) in the Notations (\ref{productSurfEll}):
\begin{equation} \label{BASEL}
 \mathcal{B}(\mathcal{L}) = \mathcal{B}(\PMp) \times div(\theta^{\E}_1) \cup \mathcal{B}(\PMm) \times div(\theta^{\E}_0) \PUNKT
\end{equation}

\noindentTX By a slight abuse of notation we shall refer to them as points $P=(x,s)$ and $Q=(y,t)$ in the product $T = \B \times \E$ as introduced in Notation \ref{productSurfEll}. We see that the equation of $\Surf$ can be written in the form
\begin{equation}\label{product122Surfequation}
\theta = \cf \theta^{\E}_0 + \cg \theta^{\E}_1 = 0 \BEISTRICH
\end{equation}
for certain general global sections $\cf \in H^0(\SA, \PMp)$ and $\cg \in H^0(\SA, \PMm)$. Then, by Proposition \ref{desccanonicalmap}, the canonical map of $\Surf$ has the following structure:

\begin{equation*}
\CAN{\Surf} = \left[H^0(\SA, \PMp)\theta^{\E}_0, H^0(\SA, \PMm)\theta^{\E}_1, \frac{\partial \cf}{\partial z_1} \theta^{\E}_0 + \frac{\partial \cg}{\partial z_1} \theta^{\E}_1, \frac{\partial \cf}{\partial z_2} \theta^{\E}_0 + \frac{\partial \cg}{\partial z_2} \theta^{\E}_1, \cf \theta'^{\E}_0 + \cg \theta'^{\E}_1 \right] \PUNKT
\end{equation*}

\noindentTX Hence, denoted by $\cfc := div(\cf)$ and $\cgc := div(\cg)$ the curves in $\SA := \bigslant{B}{e_1 + e_2}$, respectively in the linear systems $|\PMp|$ and $|\PMm|$, we have to the following cases:
\begin{itemize}
 \item[a)] \underline{$\theta^{\E}_0(s) = \theta^{\E}_0(t) = 0$}. In this case, the surface equation in (\ref{product122Surfequation}) implies that both $x$ and $y$ represent points lying on the curve $\cgc$ in $\SA$. Thus, by considering only the non-vanishing components of the expression of $\CAN{\Surf}$ when evaluated on $P$ and $Q$, we have:
 \begin{equation}\label{CANEXP1}
 \CAN{\Surf} = \left[H^0(\SA, \PMm)\theta^{\E}_1,  \frac{\partial \cg}{\partial z_1} \theta^{\E}_1, \frac{\partial \cg}{\partial z_2} \theta^{\E}_1,\cf \theta'^{\E}_0 + \cg \theta'^{\E}_1 \right] \PUNKT
\end{equation}
Again by Proposition \ref{desccanonicalmap}, since a set of generators for $H^0(\cgc, \omega_{\cgc})$ is obtained by considering the derivatives $\frac{\partial \cg}{\partial z_1}$ and $\frac{\partial \cg}{\partial z_2}$ together with the restriction of the sections of $H^0(\SA, \PMm)$ to $\cgc$, we can conclude that $\phi_{\omega_{\cgc}}(x) = \phi_{\omega_{\cgc}}(y)$, where $\phi_{\omega_{\cgc}}$ is the canonical map of $\cgc$. By applying Remark \ref{osservationsingularfibers}, our generality assumption implies that $\cgc$ is not hyperelliptic, and we infer that $x = y$ as points in $\SA$. Thus, $x = y$ or $x = y + e_1 + e_2$ in $\B$, and we can assume that $(x,s) = (y, -t)$ or $(x,s) = (y, -t + e_3)$ as points on $A$. This latter second subcase can be excluded, since it occurs only if $x$ or $y$ belongs to $\mathcal{B}(\PMm)$. Indeed, by the expression of the base locus of $\mathcal{L}$ in (\ref{BASEL}) it would lead to the conclusion, contrary to our hypothesis, that both $P$ and $Q$ are base points. Thus, $(x,s) = (y,-t)$ (i.e. $P = \inv_3(Q)$), and from the expression of the canonical map of $\Surf$ in (\ref{CANEXP1}) we conclude that both $P$ and $Q$ belong to $W_3 = div\left(\frac{\partial \theta}{\partial z_3}\right)$.

\item[b)] \underline{$\theta^{\E}_1(s) = \theta^{\E}_1(t) = 0$.} This case can be considered to be equivalent to the previous case.
\item[c)] \underline{$\theta^{\E}_0(s) = 0$, $\theta^{\E}_0(t) \neq 0$ but $y \in \mathcal{B}(\PMp)$ as point on $\SA$.} Then, according to Equation \ref{product122Surfequation}, $x \in \cgc$.
Since we are assuming that $P$ and $Q$ are ouside of the base locus of $\mathcal{L}$ (which contains $\left(\mathcal{B}(\PMp) \times div(\theta^{\E}_1)\right)$), we have that $\theta^{\E}_1(t) \neq 0$. Then, since $\CAN{\Surf}(P) = \CAN{\Surf}(Q)$, we have that $y$ also belongs to $\cgc$. Hence $\cgc$ is a curve in the linear system $|\PMp|$ passing through $y$ which is a base point of $|\PMm|$. Since by Remark \ref{osservationsingularfibers} this implies that $\cgc$ is hyperelliptic, this case can be excluded by our generality assumption on $\cf$ and $\cg$.
\item[d)] \underline{$\theta^{\E}_0$ and $\theta^{\E}_1$ do not vanish on $s$ or $t$, and neither $x$ nor $y$ belongs to $\mathcal{B}(\PMp) \cup \mathcal{B}(\PMm)$}. In this case, we can consider the 
rational map on $\SA$ defined by the product of the linear systems $|\PMp|$ and $|\PMm|$
\begin{equation*}
\ratA \colon \SA \dashrightarrow \PR H^0(\SA, \PMp) ´\times \PR H^0(\SA, \PMm) \cong \PR^1 \times \PR^1
\end{equation*}
The map $\ratA$ is clearly defined outside the union of the base loci $\mathcal{I} = \mathcal{B}(\PMp) \cup \mathcal{B}(\PMm)$, it is dominant and generically finite of degree $4$. Since $\SA$ is the quotient of a product of two elliptic curves, this map factors through the involutions $\inv_1$ and $\inv_2$. Hence, in this case we have that $y \in \left<\inv_1, \inv_2\right>.x$ as points of $\SA$, and $t = \pm s$. It is now straightforward to prove that there is no other possibility, and that, if $Q = \inv_j(P)$, then $P$ and $Q$ lie on the divisor $W_j$. This finishes the proof of the second claim. \newline
\end{itemize}

\noindentTX We prove now the third claim of the proposition.


One can show, in general, that the complementary set of an ample divisor in an abelian variety is an affine set. Unfortunately it is not always possible to achieve an useful description of these affine sets in terms of coordinates and equations. However, our special situation suggests the idea to use affine coordinates on the factors $\E_j$ and to use them to define suitable local affine coordinates on $\Surf$ and to give an expression of the canonical map of $\Surf$ in local coordinates. 

As we have already observed throughout Example \ref{basicexampleEll}, the elliptic curves $\E_j$ are defined in an affine open set $\EU{j} = \E_j - div(\theta^{\E_j}_0)$ by an equation $\Eeq{j} = 0$, where:
\begin{equation}\label{LegendreNormalFormE2}
\Eeq{j} := \y{j}^2 - (\x{j}^2 - 1)(\x{j}^2 - \pa{j}^2) \PUNKT
\end{equation}
In the affine open set $\EV{j} = \E_j - div(\theta^{\E_j}_1)$ around the poles $\inftyplus$ and $\inftyminus$ of the function $\x{j}$, we have coordinates $(\vinf{j}, \winf{j})$, with
\begin{equation}\label{LegendreNormalForm2}
\Eeqinf{j} := \winf{j}^2 - (1 - \vinf{j}^2)(1- \pa{j}^2\vinf{j}^{2}) \PUNKT
\end{equation}
On the open subset $\mathcal{U}_{000} := T - div(\theta_{000})$ of $T = \E_1 \times \E_2 \times \E_3$ ,we have coordinates $\left( \binom{\x{1}}{\y{1}},\binom{\x{2}}{\y{2}},\binom{\x{3}}{\y{3}}\right)$, and since $\Surf = div(\theta)$ for some general holomorphic section $\theta$, we have that the pullback $\widehat{\Surf}$ of $\Surf$ in via the projection $\prTA \colon T \longrightarrow A$ can be expressed, according to Equation \ref{LegendreRelationsEQ} as the closed subset in the open set $\mathcal{U}_{000}$ defined by
\begin{equation}\label{SURFEQU}
\widehat{\Surf} \cap \mathcal{U}_{000} = \left\{\left(\binom{\x{1}}{\y{1}},\binom{\x{2}}{\y{2}},\binom{\x{3}}{\y{3}}\right) \in \mathcal{U}_{000} \ \ | \ \ f = f_{000} := 1 + \cb \x{2}\x{3}+\cc \x{1}\x{3} + \cd \x{1}\x{2} = 0  \right\} \BEISTRICH
\end{equation}
where $\cb$, $\cc$ and $\cd$ are general (in particular non-zero) complex coefficients. This generality hypothesis implies that no couple of coordinates $\x{i}$ and $\x{j}$ can simultaneously vanish. Moreover, according with Equation \ref{LegendreRelationsEQ}, the action on the elliptic curve $\E_j$ with a translation by $e_j$ changes the signs of both coordinates $(\x{j}, \y{j})$ (resp. $(\vinf{j}, \winf{j})$) while the involution $\iv{j}$ only inverts the sign of $\y{j}$ (resp $\winf{j}$). 

\noindent By applying the coordinate changes $(\x{j}, \y{j}) \mapsto (\vinf{j}, \winf{j}) = (\x{j}^{-1},  \y{j}\x{j}^{-2})$ to the elliptic curves $\E_j$, we obtain the equation of $\widehat{\Surf}$ in the other open sets of the form $\mathcal{U}_{ijk} := T - div(\theta_{ijk})$. For the reader convenience, we can express the defining equation in the open set $\mathcal{U}_{110}$. To this purpose we only need to write $\vinf{1}^{-1}$ and $\vinf{2}^{-1}$ in place of $\x{1}$ and $\x{2}$ and to multiply the resulting equality by $\vinf{1}\vinf{2}$. We obtain:
\begin{equation*}
\widehat{\Surf} \cap \mathcal{U}_{110} = \left\{\left(\binom{\vinf{1}}{\winf{1}},\binom{\vinf{2}}{\winf{2}},\binom{\x{3}}{\y{3}}\right) \in \mathcal{U}_{110} \ \ | \ \ f_{110} := \vinf{1}\vinf{2} + (\cb + \cc)\x{3} + \cd = 0  \right\} \PUNKT
\end{equation*}

\noindentTX It is easy to see that the preimage $B := \prTA^{-1}(\mathcal{B}({\mathcal{L}}))$ in $T$ of the base locus in $A$ is the union of the orbits of the points $(\inftyplus, \inftyplus, \inftyplus)$, $(\inftyplus, 0, 0)$, $(0, \inftyplus, 0)$ and $(0,0,\inftyplus)$ with respect to the group isomorphic to $\Z_2^3$ which operates on each factor as the sign change $(\x{j}, \y{j}) \mapsto (-\x{j}, -\y{j})$. \newline

\noindentTX With this local description of $\widehat{\Surf}$ in Equation \ref{SURFEQU}, we express the canonical map in local coordinates and we determine the locus where the rank of the differential drops. \newline


\noindentTX For every $(ij) \in \{(12),(13),(23)\}$, the holomorphic $2$-forms $\omega_{ij}:= dz_i \wedge dz_j$ on $A$ can be looked at as non-zero elements of $H^0(\Surf, \omega_{\Surf})$ when restricted to $\Surf$. However, around a point of $\E_j$ on which $\y{j}$ does not vanish (i.e. the $j$-th component is not $2$-torsion on $\E_j$)
the holomorphic 1-form $dz_j$ is equivalent up to a non-zero constant to $\frac{d\x{j}}{\y{j}}$, while around a point at infinity it is equivalent to $\frac{d\vinf{j}}{\winf{j}}$. 
The remaining holomorphic forms of the basis for $H^0(\Surf, \omega_{\Surf})$ can be expressed by applying the residue map $H^0(A, \OC_{A}(\Surf)) = H^0(A, \omega_{A}(\Surf)) \longrightarrow H^0(\Surf, \omega_{\Surf})$, again by using a local equation of $\Surf$. For example, we can see how the canonical map looks like at a point $P = (\x{j}, \y{j})_j$ in $\widehat{\Surf} \cap \mathcal{U}_{000}$, for which the coordinates $\y{j}$ do not vanish. Since $\Surf$ is smooth, we can assume that, in $P$,
\begin{equation} \label{temporaryassumption1proof}
\frac{\partial f}{\partial \x{3}} = \cb\x{2}+\cc\x{1} \neq 0 \PUNKT
\end{equation}
Thus, we can use $\x{1}$ and $\x{2}$ as local parameters of $\Surf$ in $P$, and we can express the global sections of the canonical bundle of $\Surf$ locally in $P$ as holomorphic forms of the type $g(x_1, x_2) d\x{1} \wedge d\x{2}$, where $g$ denotes a function around $P$ without poles. With this procedure, we can write:
\begin{equation} \label{eqdiff1}
\begin{split}
\omega_{12} &:= dz_1 \wedge dz_2 = \frac{d\x{1}}{\y{1}} \wedge \frac{d\x{2}}{\y{2}} = \frac{1}{\y{1} \y{2}} d\x{1} \wedge d\x{2} \\
\omega_{13} &:= dz_1 \wedge dz_3 = \frac{d\x{1}}{\y{1}} \wedge \frac{d\x{3}}{\y{3}} = -\frac{\cd \x{1} + \cb \x{3}}{(\cb \x{2} + \cc \x{1})\y{1}\y{3}} d\x{1} \wedge d\x{2} \\
\omega_{23} &:= dz_2 \wedge dz_3 =  \frac{d\x{2}}{\y{2}} \wedge \frac{d\x{3}}{\y{3}} = \frac{(\cc \x{3} + \cd \x{2})}{(\cb \x{2} + \cc \x{1})\y{1}\y{2}} d\x{1} \wedge d\x{2}  \PUNKT
\end{split}
\end{equation}
\noindentTX We write down also the global holomorphic differentials which arise by the residue map $H^0(A, \OC_{A}(\Surf)) = H^0(A, \omega_{A}(\Surf)) \longrightarrow H^0(\Surf, \omega_{\Surf})$. We denote, with $(ijk) \in \{(000),(011),(101),(110)\}$, 
\begin{align} \label{eqdiff1a}
\psi_{ijk} := (\theta_{ijk} \cdot d z_{1} \wedge d z_{2} \wedge d z_{3})  \neg \left( \frac{\partial \x{3}}{\theta_{000}\frac{\partial f}{\partial \x{3}}}  \right) \BEISTRICH
\end{align}
where $\neg$ is the contraction operator. 
In conclusion, we have, up to a non-zero constant:
\begin{equation}\label{eqdiff2}
\begin{split}
\psi_{000} &= \frac{1}{(\cb\x{2}+\cc\x{1})\y{1}\y{2}\y{3}} d \x{1} \wedge d \x{2} \\
\psi_{011} &= \frac{\x{2} \x{3}}{(\cb\x{2}+\cc\x{1})\y{1}\y{2}\y{3}} d \x{1} \wedge d \x{2} \\
\psi_{101} &= \frac{\x{1} \x{3}}{(\cb\x{2}+\cc\x{1})\y{1}\y{2}\y{3}} d \x{1} \wedge d \x{2} \\
\psi_{110} &= \frac{\x{1} \x{2}}{(\cb\x{2}+\cc\x{1})\y{1}\y{2}\y{3}} d \x{1} \wedge d \x{2}  \PUNKT
\end{split}
\end{equation}
Once we have multiplied all expressions in (\ref{eqdiff1}) and (\ref{eqdiff2}) by ${(\cb\x{2}+\cc\x{1})\y{1}\y{2}\y{3}}$, we obtain the following expression of the canonical map of $\Surf$, which is defined independently on the assumption in (\ref{temporaryassumption1proof}) and every point of the affine space $\mathbb{A}^6$ with coordinates $(\x{1}, \x{2}, \x{3},\y{1}, \y{2}, \y{3})$:
\begin{equation*}
\phi_{\Surf} = \begin{bmatrix}
(\cb\x{2}+\cc\x{1})\y{3} & (\cb\x{3}+\cd\x{1})\y{2} & (\cd\x{2}+\cc\x{3})\y{1} & 1 & \x{1}\x{2} & \x{1}\x{3} & \x{2} \x{3}  \PUNKT
\end{bmatrix} \PUNKT
\end{equation*}
To compute its rank, we have to consider the equations $\Eeq{j} = 0$ of $\E_j$, considered in the form as in Equation \ref{LegendreNormalFormE2}, and the affine map $\Phi: \mathbb{A}^6 \longrightarrow \mathbb{A}^9$ defined by
\begin{equation}\label{mapPsiaffinecase}
\Phi = \begin{pmatrix}
(\cb\x{2}+\cc\x{1})\y{3} & (\cb\x{3}+\cd\x{1})\y{2} & (\cd\x{2}+\cc\x{3})\y{1} & \x{1}\x{2} & \x{1}\x{3} & \x{2} \x{3} & g_1 & g_2 & g_3   \\
\end{pmatrix} \PUNKT
\end{equation}
The differential of the canonical map of $\Surf$ is injective at the point $P$ if the matrix of the differential of $\Phi$ at $P$ has maximal rank. The matrix of the differential in $P$ is obtained, as usual, by considering the derivatives with respect to the coordinates $\x{j}$ and $\y{j}$, which is:
\begin{equation} \label{matrixNequation}
\mbox{\scriptsize{$
   N := \begin{bmatrix}
      \cc\y{3}  &       \cd\y{2}  &       0  &    0 &   \x{3} &  \x{2}  &      2\x{1}(\pa{1}^2-2\x{1}^2+1)    &        0  &     0                   \\ 
      \cb\y{3}   &      0   &        \cd\y{1}  &    \x{3} &  0 &   \x{1}  &     0               &             2\x{2}(\pa{2}^2-2\x{2}^2+1)   & 0               \\
      0        &   \cb\y{2}  &       \cc\y{1}  &   \x{2} &  \x{1} &  0 &    0            &                0  &    2\x{3}(\pa{3}^2-2\x{3}^2+1)   \\
      0       &    0       &    \cd\x{2}+\cc\x{3} &   0  &  0  &  0 &        2\y{1}   &                       0 &      0                  \\
      0        &   \cd\x{1}+\cb\x{3} &    0      &   0  &  0  &  0 &        0      &                 2\y{2} &                  0           \\
      \cc\x{1}+\cb\x{2} &    0       &    0       &   0  &  0  &  0  &      0     &                       0      &     2\y{3}    \\                
\end{bmatrix}
$}}
\end{equation}

\noindentTX It is now straightforward to see that the determinant of the submatrix formed by the first six colums from left is:
\begin{equation*} 
 det(N_{1,2,3,4,5,6}) := 2\x{1}\x{2}\x{3}(\cc\x{1}+\cb\x{2})(\cd\x{1}+\cb\x{3})(\cd\x{2}+\cc\x{3})
\end{equation*}
On the other side, the determinant of the submatrix formed by deleting the central three colums from is:
\begin{equation} \label{FIRSTDETRELATION}
 det(N_{1,2,3,7,8,9}) := 16\x{1}\x{2}\x{3}\y{1}\y{2}\y{3}
\end{equation}
Hence, the rank of $N$ may decrease only on the divisors $div(\x{j})$ for some $j$. Indeed, by the defining equation of $\Surf$ in (\ref{SURFEQU}) and by our generality hypothesis on the coefficients in the equation, at most one of the factors involved in the previous determinant expressions in (\ref{FIRSTDETRELATION}) can vanish.

If $\x{3}=0$ (which implies that $\x{1} \neq 0$ and $\x{2} \neq 0$), we can assume, without loss of generality, that $\y{3} = \delta_{3}$.
and the differential matrix $N$ in (\ref{matrixNequation}) has the following form
\begin{equation*}
\mbox{\scriptsize{$
N := \begin{bmatrix}
      \cc \delta_{3}  &       \cd\y{2}  &       0  &    0 &   0 &  \x{2}  &      2\x{1}(\pa{1}^2-2\x{1}^2+1)    &        0  &     0                   \\ 
      \cb \delta_{3}  &      0   &        \cd\y{1}  &    0 &  0 &   \x{1}  &     0               &             2\x{2}(\pa{2}^2-2\x{2}^2+1)   & 0               \\
      0        &   \cb\y{2}  &       \cc\y{1}  &   \x{2} &  \x{1} &  0 &    0            &                0  &   0   \\
      0       &    0       &    \cd\x{2} &   0  &  0  &  0 &        2\y{1}   &                       0 &      0                  \\
      0        &   \cd\x{1} &    0      &   0  &  0  &  0 &        0      &                 2\y{2} &                  0           \\
      \cc\x{1}+\cb\x{2} &    0       &    0       &   0  &  0  &  0  &      0     &                       0      &     2\delta_{3}    \\                
\end{bmatrix}
$}} \PUNKT
\end{equation*}

Finally, we notice that
\begin{align*}
 det(N_{1,4,6,7,8,9}) &=  -8\cc \y{1}^2\y{2}(\cc \x{1}-\cb \x{2}) \\
 det(N_{2,4,6,7,8,9}) &=   8\cd \x{1}\y{1}\x{2}(\x{2}^2-\pa{2}^4)\\
 det(N_{3,4,6,7,8,9}) &=  -8\cd \y{2}\x{2}^2(\x{1}^2-\pa{1}^4) \BEISTRICH
\end{align*}
which, according to our generality condition can not simultaneously vanish. Indeed, if this were the case and all $\y{i}$ were non-zero, then we would have that $\cc \x{1}-\cb \x{2} = 0$ and $\x{1}^2-\pa{1}^4 = \x{2}^2-\pa{2}^4 = 0$. But also this situation can be avoided if we suppose the coefficients $\cb,\cc,\cd$ to be sufficiently general. If otherwise $\y{1} = 0$, then we have clearly that $det(N_{3,4,6,7,8,9}) \neq 0$, and the conclusion follows.
\bigskip 

On the other side, the matrix in (\ref{matrixNequation}) shows that injectivity of the differential of $\Phi$ fails at the points where $y_3 = 0 = \cc\x{1}+\cb\x{2}$, $y_2 = 0 =  \cd\x{1}+\cb\x{3}$ or $\y{1} = 0 = \cd\x{2}+\cc\x{3}$. Since in local coordinates we have that $W_1 = div(\cd\x{2}+\cc\x{3})$, $W_2 = div(\cd\x{1}+\cb\x{3})$ and $W_3 = div(\cc\x{1}+\cb\x{2})$, these points are the $24$ points which belong to some divisor $W_j$ and are fixed under the corresponding involution $\inv_j$, as claimed. \newline

\noindentTX It remains to show that the differential of the canonical map has maximal rank at a base point of the polarization, when the choice of the coefficients in (\ref{SURFEQU}) is sufficiently general. \newline

\noindentTX Let us assume that $P = (\inftyplus, \inftyplus, \inftyplus)$ (the other cases can be treated similarly). Around $P$ we use the coordinates $(\vinf{j}, \winf{j})_j$ of the open set $\mathcal{U}_{111}$ the equation of $\widehat{\Surf}$ is:
\begin{equation*}
f_{\infty} = \vinf{1}\vinf{2}\vinf{3} + \cb \vinf{1} + \cc \vinf{2} + \cd \vinf{3} = 0
\end{equation*}
According with the affine equations in (\ref{LegendreNormalFormE2}) and (\ref{LegendreNormalForm2}) of the elliptic curves $\E_j$, we have that, in $P$, all coordinates $\vinf{j}$ vanish, while $\winf{j} = 1$ for each $j$. We have that
\begin{equation*}
\frac{\partial}{\partial \vinf{3}}f_{\infty} = \vinf{1}\vinf{2} + \cd
\end{equation*}
which is non-zero if we assume that $\cd \neq 0$.
If we repeat the same procedure, 
we obtain a map $\mathbb{A}^6 = \mathbb{A}(\vinf{j}, \winf{j})_j \longrightarrow \mathbb{A}^9$ which represents the canonical map around $P$:

\begin{equation} \label{canonicalexpressionatinfty}
\Phi_{\infty} = \begin{bmatrix}
((\vinf{1}\vinf{2} + \cd)\winf{3} & (\vinf{1}\vinf{3} + \cc)\winf{2} & (\vinf{2}\vinf{3} + \cb)\winf{1}  &  \vinf{1} &  \vinf{2} & \vinf{3} &  \Eeqinf{1} &  \Eeqinf{2} & \Eeqinf{3} \\
\end{bmatrix} \PUNKT
\end{equation}
It is now easy to verify that the matrix of the differential in $P$ is of the form
\begin{equation*}
\begin{bmatrix}
0 & 0 & 0  &  1 &  0 & 0 &  * &  * & * \\
0 & 0 & 0  &  0 &  1 & 0 &  * &  * & * \\
0 & 0 & 0  &  0 &  0 & 1 &  * &  * & * \\
\cd & 0 & 0  &  0 &  0 & 0 &  * &  * & * \\
0 & \cc & 0  &  0 &  0 & 0 &  * &  * & * \\
0 & 0 & \cb  &  0 &  0 & 0 &  * &  * & *
\end{bmatrix} \BEISTRICH
\end{equation*}
which has clearly maximal rank when the coefficients $\cb$, $\cc$ and $\cd$ are non-zero. This completely proves the third claim of the proposition.
\end{proof}

\section{On the canonical map in the general case} \label{GENERALFINALSECTION}

Following the notations of the previous section (see Notation \ref{productSurfEll}), to every symmetric matrix $\tau = (\tau_{ij})$ in the Siegel half space $\SiegH_3$, there is a corresponding $(2,2,2)$-polarized abelian threefold $T = \bigslant{\C^3}{\tau \Z^3 \oplus 2\Z^3}$ with a projection onto a $(1,2,2)$-polarized abelian threefold $A$. This latter abelian threefold is isogenous to a product of three elliptic curves if $\tau$ belongs to the diagonal subspace $\SiegH_1^3 = \SiegH_1 \times \SiegH_1 \times \SiegH_1$ of $\SiegH_3$. In such a case, Proposition \ref{generalbehavior1} tells us which pairs of points are critical for the canonical map of a general surface $\Surf = div(\theta)$ in the $(1,2,2)$-polarization of $A$, according to the following general definition:
\begin{def1}\label{CRITICALPAIRS}
Let $\NDiv$ be a smooth ample divisor on an abelian variety $A$. Denoted by $\CAN{\NDiv}$ the canonical map of $\NDiv$, a \textit{\textbf{critical pair}} for $\NDiv$ is a couple of points $(P, Q)$ such that:
\begin{itemize}
 \item $\CAN{\NDiv}(P) = \CAN{\NDiv}(Q)$ if $P \neq Q$.
 \item The differential of $\CAN{\NDiv}$ at $P$ has not maximal rank, if $P = Q$.
\end{itemize}
\end{def1}
Recall that, denoted by $(z_1, z_2, z_3)$ a coordinate system for an abelian threefold $A$ which is isogenous to the product of three elliptic curves, we proved in Proposition \ref{generalbehavior1} that the canonical map of $\Surf = div(\theta) \subseteq A$ fails to be an embedding precisely on those pairs of points which lie on the canonical divisors 
\begin{equation*}
W_k := div\left(\frac{\partial \theta}{\partial z_k}\right) = div\left(dz_i \wedge dz_j \right) \BEISTRICH
\end{equation*}
where $(i,j,k)$ is a permutation of the indices $(1,2,3)$. This happens because the divisors $W_k$ are invariant for the involution $\inv_k \colon z_k \mapsto -z_k$.


\begin{rmk} (On the generality condition in the statement of Proposition \ref{generalbehavior1}) \label{BEQUADROG}
Let us consider a diagonal matrix $\tau_0 \in \SiegH_1^3$ in the Siegel upper half space $\SiegH_3$. As in Notation \ref{productSurfEll}, for ever couple of indices $(i,j)$ there is a $(1,2)$-polarized abelian surface $\SA_{ij}$ with a projection:
\begin{equation*}
\prBSA_{ij} \colon \B_{ij}:= \E_i \times \E_j \longrightarrow \SA_{ij}
\end{equation*}
For the $(1,2,2)$-abelian threefold $(A, \mathcal{L})$, which is isogenous to via a degree $2$ projection $\prTA \colon T \longrightarrow A$, we have that the vector space of sections of the polarization decomposes, for every couple of indices $(i,j)$, as:
\begin{equation} \label{sectionsDECOMP2}
H^0(A, \mathcal{L}) \cong H^0(\SA_{ij}, \PMpij{ij})\cdot \theta^{\E_k}_0 \oplus H^0(\SA_{ij}, \PMmij{ij})\cdot \theta^{\E_k}_1
\end{equation}
In the proof of Proposition \ref{generalbehavior1}, a smooth surface $\Surf := div(\theta)$ in the polarization class $|\mathcal{L}|$ of such an abelian threefold $A$ were required to fulfill some conditions in order to be considered as \textbf{sufficiently general}. There conditions can be resumed as follows:
\begin{itemize}
 \item For every permutation $(i,j,k)$ of the indices $(1,2,3)$, the divisors $\cfc_{ij} := div(\cf_{ij})$ and $\cgc_{ij} := div(\cg_{ij})$ of $\SA_{ij}$, which are defined according to the decomposition in (\ref{sectionsDECOMP2}) such that
 \begin{equation*}
  \theta = \cf_{ij}\cdot \theta^{\E_k}_0 + \cg_{ij}\cdot \theta^{\E_k}_1
 \end{equation*}
are smooth irreducible genus $3$ curves. 
\item For every point $P$ on $\Surf$ it holds that:
\begin{itemize}
\item If the sections $\cf_{ij}$ and $\cg_{ij}$ both vanish on $P$, then neither the $i$-th coordinate nor the $j$-th coordinate of $P$  are represented by $2$-torsion on $\E_i$, resp. $\E_j$.
\item If the sections $\cf_{ij}$ and $\cg_{ij}$ both vanish on $P$, then $\cf_{ik}$ and $\cg_{ik}$ do not simultaneously vanish at $P$.
\end{itemize}
\end{itemize}
If the latter generality conditions are fulfilled, then the proof of Proposition \ref{generalbehavior1} ensures that the only critical pairs on $\Surf$ are as claimed in the statement of same proposition. 

\end{rmk}

In this section, we aim to prove that, when $\tau$ is sufficiently general, there are no critical pairs on a general surface $\Surf$ within the $(1,2,2)$ polarization of $A$. That means that the canonical map yields a holomorphic embedding in $\PR^5$ in the general case, in agreement with Theorem \ref{teoremafinale}:
\begin{thm1*} 
Let be $(A, \mathcal{L})$ a general $(1,2,2)$-polarized abelian threefold and let be $\Surf$ a general surface in the linear system $|\mathcal{L}|$. Then the canonical map of $\Surf$ is a holomorphic embedding. 
\end{thm1*}



%
\begin{sketch}
For the sake of clarity and exposition, we introduce at this place the basic ideas behind the proof of Theorem \ref{teoremafinale}, and we give a complete proof of it at the end of this section. \newline

\noindentTX We fix a certain general diagonal matrix $\tau_0 = diag(\tau_{ii})_i$ in the subspace $\SiegH_1^3$ of $\SiegH_3$ and we denote by 
\begin{equation}\label{HDCENTRDEF}
 \HD := \SET{\tau \in \SiegH_3}{\text{The diagonal entries of $\tau$ coincide with those of $\tau_0$}} \PUNKT
\end{equation}
To every open neighborhood $\mathcal{U}$ of $\tau_0$ in $\HD$ corresponds a family of $(1,2,2)$-polarized abelian varieties $(A_{\mathcal{U}}, \mathcal{L}_{\mathcal{U}})$. The abelian threefold $(A_{\tau_0}, \mathcal{L}_{\tau_0})$ of the family is isogenous to the product of three elliptic curves as in Remark \ref{BEQUADROG}. Every non-zero holomorphic section $\theta$ of $H^0(A_{\tau_0}, \mathcal{L}_{\tau_0})$ gives rise to a family of surfaces $\Surf_{\mathcal{U}}$, each contained in the polarization class of the respective fiber of the family of abelian threefolds $A_{\mathcal{U}}$. We have, for every element $\tau$ of the basis $\mathcal{U}$, a natural identification:
\begin{equation} \label{IDENTIFICATION}
H^0(A_{\tau_0}, \mathcal{L}_{\tau_0}) \cong H^0(A_{\tau}, \mathcal{L}_{\tau}) \cong \PR^3 \BEISTRICH
\end{equation}
defined by the natural restriction of the holomorphic sections of $H^0(A_{\mathcal{U}}, \mathcal{L}_{\mathcal{U}})$ to the fibers of the family:
\begin{equation} 
H^0(A_{\mathcal{U}}, \mathcal{L}_{\mathcal{U}})|_{\tau} \cong H^0(A_{\tau}, \mathcal{L}_{\tau}) \PUNKT
\end{equation}


\noindentTX If we assume our claim to be false, then we expect to find, for every $\tau$ in a sufficiently small neighborhood $\mathcal{U}$ of $\tau_0$, a critical pair $(P_{\tau}, Q_{\tau})$ for the canonical map of $\Surf_{\tau}$. 

This means that we have a whole family $\DEFP$ on $\mathcal{U}$ of critical pairs $(P_{\mathcal{U}}, Q_{\mathcal{U}})$, which gives rise, by restriction to the closed loci 
\begin{equation} \label{LOCIEQ}
\HD^{(ij)} := \{\tau_{ik} = \tau_{jk} = 0\} \cong \SiegH_2 \times \SiegH_1 \BEISTRICH
\end{equation}
to three different families $\DEFP_{ij}$ of critical pairs, each defined respectively on $\mathcal{U}_{ij} := \mathcal{U} \cap \HD^{(ij)}$:
\begin{equation} \label{restrictionPdiagram}
\begin{diagram}
&  &\DEFP_{ij} &\rInto{}& \DEFP\\
&  &\dTo_{}    &      &\dTo\\
&  &\mathcal{U}_{ij} := \mathcal{U} \cap \HD^{(ij)}         &\rInto &\mathcal{U}\\
\end{diagram}
\end{equation}

It is then natural to ask whether we can classify all possible critical pairs on (general) surfaces in abelian threefolds which lie in the loci $\HD^{(ij)}$, similarly to what we did in Proposition \ref{generalbehavior1}. Note that these loci parametrize $(1,2,2)$-abelian variety which are isogenous to products of a $(2,2)$-polarized abelian surface $\B$ and a $(2)$-polarized elliptic curve $\E$, according to Definition \ref{LOCIEQ} and Notation \ref{productSurfEll}.


The following proposition, which we state for simplicity for the couples of indices $(i,j) = (1,2)$, answers to this question.
\end{sketch}

\begin{prop1}\label{generalbehavior2} 
Let $\tau_0 = diag(\tau_{11}, \tau_{22}, \tau_{33})$ be a diagonal matrix in the Siegel upper half-space $\mathcal{H}_{3}$, let $(A_{\tau_0}, \mathcal{L}_{\tau_0})$ be the associated $(1,2,2)$-polarized abelian variety with a general divisor $S_{\tau_0} = div(\theta)$ in the polarization $|\mathcal{L}_{\tau_0}|$. Then, for a sufficiently small neighborhood $\mathcal{U}$ of $\tau_0$ in the closed locus $\HD^{(12)}$, for every surface $\Surf_{\tau} = div(\theta)$ of the restricted family $\mathcal{U}_{12} := \mathcal{U} \cap \HD^{(12)}$ (see Diagram \ref{restrictionPdiagram}), it holds that

\begin{itemize}
 \item The involution $\iv{3} \colon (z_1, z_2, z_3) \mapsto (z_1, z_2, -z_3)$ on $A_{\tau}$ leaves $\Surf_{\tau}$ invariant and induces an involution on the canonical divisor $W_3 := div \left( \frac{\partial \theta}{\partial z_3} \right)$ on $\Surf_{\tau}$. (Recall that, in general, $W_k := div \left( \frac{\partial \theta}{\partial z_k} \right)$.)
 
 \item The canonical map of $\Surf$ is one-to-one, except:
 \begin{itemize}
 \item on the divisor $W_3$, on which the canonical map has degree $2$ and factors respectively through the involution $\iv{3}$.
 
 \item on the finite set $W_1 \cap W_2$, on which the canonical map of $\Surf$ factors through the involution  $\iv{1}\cdot \iv{2}: (z_1, z_2, z_3) \mapsto (-z_1, -z_2, z_3)$.
\end{itemize} 
\item The differential of the canonical map of $\Surf$ has everywhere maximal rank, except on those points of $W_3$ which are fixed points under the action of the involution $\iv{3}$ on $W_3$.
\end{itemize}
\end{prop1}

\noindentTX Proposition \ref{generalbehavior2} shows, in particular, how the behavior of the canonical map changes when, in the moduli space, we move away from the locus of abelian threefolds, which are isogenous to a product of elliptic curves, and we move 
along closed loci parametrizing abelian threefolds which are isogenous to a product of a simple $(2,2)$ abelian surface and an elliptic curve.

From Proposition \ref{generalbehavior2} also follows also that every restriction $\DEFP_{ij}$ as in Diagram \ref{restrictionPdiagram} of a family $\DEFP$ of critical pairs only contains pairs of points which are conjugate under the same involution. However, on the central fiber $\tau_0$, (that is, on the intersection of these loci) the three different families of critical pairs $\DEFP_{ij}$ specialize to a unique pair of points $(P,Q)$ which must be conjugated under many different involutions, which would lead to a contradiction.

After this rough explaination of the basic approach to our problem, we aim to formalize the steps toward the proof of Theorem \ref{teoremafinale}.

\begin{def1}\label{generaldeformationdefinition}
In the notation of Proposition \ref{generalbehavior2}, let us suppose that $(P,Q)$ is a critical pair on $\Surf_{\tau_{0}} = div(\theta)$. A \textbf{family of critical pairs around $(P,Q)$} is the datum of a couple 
$(\Surf_{\mathcal{U}} \longrightarrow \mathcal{U}, \DEFP)$, where:
\begin{itemize}
\item[i)] $\mathcal{U}$ is an open set of $\HD$ (see definition of $\HD$ in Equation \ref{HDCENTRDEF}).
\item[ii)] $\Surf_{\mathcal{U}} \longrightarrow \mathcal{U}$ is a family of surfaces defined as $\Surf_{\tau} = div(\theta)$ for every $\tau$, by applying the natural identification 
\begin{equation} 
H^0(A_{\tau_0}, \mathcal{L}_{\tau_0}) \cong H^0(A_{\mathcal{U}}, \mathcal{L}_{\mathcal{U}}) \cong H^0(A_{\tau}, \mathcal{L}_{\tau}) \PUNKT
\end{equation}
 \item[iii)] With $\Delta_{\mathcal{U}}$ the diagonal subscheme in $\PR_{\mathcal{U}}^5 \times \PR_{\mathcal{U}}^5$ and $\phi_{\Surf_{\mathcal{U}}}$ the map, which on each surface $\Surf_{\tau}$ of the family coincides with the canonical map $\phi_{\Surf_{\tau}}$, $\DEFP$ is a closed irreducible subscheme of 
\begin{equation*} 
\mathcal{K}_{\mathcal{U}} := (\phi_{\Surf_{\mathcal{U}}} \times_{\mathcal{U}} \phi_{\Surf_{\mathcal{U}}})^{-1}(\Delta_{\mathcal{U}}) \subseteq \Surf_{\mathcal{U}} \times_{\mathcal{U}} \Surf_{\mathcal{U}}
 \end{equation*}
which is dominant over $\mathcal{U}$ and such that its restriction on the central fiber $\tau_0$ is $(P,Q)$. Here $\times_{\mathcal{U}}$ denotes, as usual, the cartesian product on $\mathcal{U}$ according to the pullback diagram:
\begin{equation*}
\begin{diagram} 
\Surf_{\mathcal{U}} \times_{\mathcal{U}} \Surf_{\mathcal{U}}	&\rTo			&\Surf_{\mathcal{U}}    \\ 	
		 \dTo &							& \dTo \\
\Surf_{\mathcal{U}}			&\rInto			& \mathcal{U} 
\end{diagram}
\end{equation*}
\end{itemize}
\end{def1}

\begin{rmk} \label{BEQUADROij} 
In the notation of the previous Definition \ref{generaldeformationdefinition}, if we restrict a family of critical points $(\Surf_{\mathcal{U}} \longrightarrow \mathcal{U}, \DEFP)$ to one of the loci $\HD^{(ij)}$ (recall the definition in \ref{LOCIEQ}), we obtain a subfamily of critical points $(\Surf_{\mathcal{U}_{ij}} \longrightarrow \mathcal{U}_{ij}, \DEFP_{ij})$ defined on the open set $\mathcal{U}_{ij} := \mathcal{U} \cap \HD^{(ij)}$. Furthermore, we have in this case a family of projections
\begin{equation} \label{PProjU}
p_{\mathcal{U}} \colon T_{\mathcal{U}} = B_{\mathcal{U}, ij} \times E_k \longrightarrow A_{\mathcal{U}}
\end{equation}
and a family of $(2,2)$-polarized abelian surfaces $\B_{\mathcal{U}_{ij}} \longrightarrow \mathcal{U}_{ij}$ with a projection onto a family of $(1,2)$-polarized abelian surfaces $\SA_{\mathcal{U}_{ij}} \longrightarrow \mathcal{U}_{ij}$, and there is a decomposition as in (\ref{sectionsDECOMP2}):
\begin{equation} \label{sectionsDECOMPgen}
 H^0(A_{\mathcal{U}}, \mathcal{L}_{\mathcal{U}}) \cong H^0(\SA_{\mathcal{U}_{ij}}, \PMpij{ij})\cdot \theta^{\E_k}_0 \oplus H^0(\SA_{\mathcal{U}_{ij}}, \PMmij{ij})\cdot \theta^{\E_k}_1 
\end{equation}
\end{rmk}
\begin{not1} \label{NOTNUXI}
 If we write a non-zero section $\theta$ of $H^0(A_{\mathcal{U}}, \mathcal{L}_{\mathcal{U}})$ according to (\ref{sectionsDECOMPgen}), we have
\begin{equation}\label{thetaDECOMPOSITIONij}
  \theta = \cf_{\mathcal{U}, ij}\cdot \theta^{\E_k}_0 + \cg_{\mathcal{U}, ij}\cdot \theta^{\E_k}_1 
\end{equation}
for some $\cf_{\mathcal{U}, ij}$ and $\cg_{\mathcal{U}, ij}$ in $H^0(\SA_{\mathcal{U}_{ij}}, \PMpij{ij})$ and $H^0(\SA_{\mathcal{U}_{ij}}, \PMmij{ij})$ respectively.
We denote henceforth, similarly to what we did in Remark \ref{BEQUADROG}, the following families of divisors in $\SA_{\mathcal{U}, ij}$ as follows
\begin{align*}
 \cfc_{\mathcal{U}, ij} &:= div(\cf_{\mathcal{U}, ij}) \\
 \cgc_{\mathcal{U}, ij} &:= div(\cg_{\mathcal{U}, ij})
\end{align*}

\end{not1}



\begin{oss1} \label{firstOssBXE}
\noindentTX Let us consider a family of critical pairs $(\Surf_{\mathcal{U}} \longrightarrow \mathcal{U}, \DEFP)$ and its restriction to:
\begin{equation} \label{LOCIEQ2}
\HD^{(ij)} := \{\tau_{ik} = \tau_{jk} = 0\} \cong \SiegH_2 \times \SiegH_1 \BEISTRICH
\end{equation}
We have, following our previous discussion in Remark \ref{BEQUADROij}, a family $(\Surf_{\mathcal{U}_{ij}} \longrightarrow \mathcal{U}_{ij}, \DEFP_{ij})$ on the open set $\mathcal{U}_{ij} := \mathcal{U} \cap \HD^{(ij)}$, with $\Surf_{\mathcal{U}} = div(\theta)$ which satisfy the generality conditions stated in Remark \ref{BEQUADROij}. 

Suppose that, for some $\tau$ in $\mathcal{U}_{ij}$, the point $P_{\tau}$ belongs to the canonical divisor
\begin{equation*}
W_k = div\left(\frac{\partial \theta}{\partial z_k}\right)_{\tau} \subseteq \Surf_{\tau} \PUNKT
\end{equation*}
Then the set $\prBSA_{\tau}\left(\prTA_{\tau}^{-1}(W_k)\right)$, defined on each $\tau$ according to the following diagram 
\begin{equation} \label{restrictionPdiagram}
\begin{diagram}
&  &\B_{\mathcal{U}, ij} \times \E_k &\rTo{\prTA_{\mathcal{U}}}& A_{\mathcal{U}}\\
&  &\dTo_{\prBSA_{\mathcal{U}_{ij}} \times id_{\E_k}}    &      & \\
&  &\SA_{\mathcal{U}, ij} \times \E_k         & &\\
\end{diagram}
\end{equation}
is reducible, and we have, for every $\tau$ in $\mathcal{U}_{ij}$ (see again Remark \ref{BEQUADROij} for the notations):
\begin{equation*}
 (\prBSA_{\tau} \times id_{\E_k})\left(\prTA_{\tau}^{-1}(W_k)\right) = (\cfc_{\tau, ij} \cap \cgc_{\tau, ij}) \times \E_k \cup \SA_{\tau} \times \E_k[2] \PUNKT 
\end{equation*}

Indeed, if we write the section $\theta$ as in Equation \ref{thetaDECOMPOSITIONij}, then the coordinates of the points in the divisor $W_k$ must fulfill the conditions:


\begin{equation}\label{casesDEFeta}
\begin{cases}
\cf_{\tau, ij}\theta^{E_k}_0 &+ \ \ \cg_{\tau, ij}\theta^{E_k}_1 = 0\\
\cf_{\tau, ij}\frac{\partial}{\partial z_k} \theta^{E_k}_0 &+ \ \ \cg_{\tau, ij}\frac{\partial}{\partial z_k}\theta^{E_k}_1 = 0 \PUNKT
\end{cases}
\end{equation}
The conditions in (\ref{casesDEFeta}) can be looked at as a linear system in $\cf_{\tau, ij}$ and $\cg_{\tau, ij}$. This implies that both $\cf_{\tau}$ and $\cg_{\tau}$ vanish, or the determinant 
\begin{equation} \label{etadefinition}
\eta_k := \begin{vmatrix}\theta^{E_k}_0 & \theta^{E_k}_1 \\ 
\frac{\partial}{\partial z_k} \theta^{E_k}_0 & \frac{\partial}{\partial z_k} \theta^{E_k}_1 \end{vmatrix} \in H^0(E_k, \OC_{E_k}(4O_{E_k})) \PUNKT
\end{equation}
The conclusion follows now by Equation \ref{etadefinition_Y} in Example \ref{basicexampleEll} (since $div(\eta_k) = \E_k[2]$ ).
\end{oss1}

In the next lemma, we show that the restriction $\DEFP_{ij}$ (as in Remark \ref{BEQUADROij}) of a family of critical pairs $\DEFP$ around a couple of points $(P,Q)$ with $P \neq Q$ must be contained in a closed locus of $\Surf_{\mathcal{U}_{ij}} \times_{\mathcal{U}_{ij}} \Surf_{\mathcal{U}_{ij}}$ which, in 
\begin{equation}\label{eqproduct}
\left(\SA_{\mathcal{U}, ij} \times \E_k\right)^2 := \left(\SA_{\mathcal{U}, ij} \times_{\mathcal{U}_{ij}} \SA_{\mathcal{U}, ij} \right) \times \E_k^2 \BEISTRICH
\end{equation}
is defined by a holomorphic section of the $(4,4)$-polarization $\OC_{\E_k}(4O_{\E_k}) \boxtimes  \OC_{\E_k}(4O_{\E_k})$ on $\E_k \times \E_k$. 

\begin{lem1}\label{secondOssBXE} 
\noindentTX Let us consider a family of critical pairs $(\Surf_{\mathcal{U}} \longrightarrow \mathcal{U}, \DEFP)$ and its restriction $(\Surf_{\mathcal{U}_{ij}} \longrightarrow \mathcal{U}_{ij}, \DEFP_{ij})$ to the open set $\mathcal{U}_{ij} := \mathcal{U} \cap \HD^{(ij)}$ with $\Surf_{\mathcal{U}} = div(\theta)$ which satisfy the generality conditions of Remark \ref{BEQUADROij}. Let us assume that $P_{\tau} \neq Q_{\tau}$ for every $\tau$ in $\mathcal{U}_{ij}$ (i.e., $\DEFP$ is a family of critical pairs of distinct points). \newline

Then, in the notations of the following diagram (cf. also Diagram \ref{restrictionPdiagram})

\begin{equation} \label{restrictionPdiagram2}
\begin{diagram}
&  &\left(\B_{\mathcal{U}, ij} \times \E_k\right)^2 &\rTo{\ \ \ \  \prTA_{\mathcal{U}} \times_{\mathcal{U}} \prTA_{\mathcal{U}} \ \ \ \ }& A_{\mathcal{U}}^2 = A_{\mathcal{U}} \times_{\mathcal{U}} A_{\mathcal{U}}\\
&  &\dTo_{(\prBSA_{\mathcal{U}_{ij}} \times id_{\E_k})^2}    &      & \\
&  &\left(\SA_{\mathcal{U}, ij} \times \E_k \right)^2        & &\\
\end{diagram}
\end{equation}

there exists a non-zero holomorphic section $\secV_k$ of $H^0(\E_k \times \E_k, \OC_{\E_k}(4O_{\E_k}) \boxtimes \OC_{\E_k}(4O_{\E_k}))$ such that:
\begin{equation*}
 (\prBSA_{\mathcal{U}_{ij}} \times id_{\E_k})^2\left(\left(\prTA_{\mathcal{U}} \times_{\mathcal{U}} \prTA_{\mathcal{U}}\right)^{-1}(\DEFP_{ij})\right) \subseteq \left[(\cfc_{\mathcal{U}, ij} \cap \cgc_{\mathcal{U}, ij}) \times \E_k \right]^2 \cup \left[\SA_{\mathcal{U}, ij}^2 \times div(\secV_k) \right]
\end{equation*}
\end{lem1}
\begin{proof}
It sufficies to prove the claim for a generic fiber of the family. For this purpose, we can assume $(i,j,k) = (1,2,3)$ and we pick an element $\tau \in \mathcal{U}_{12}$. The critical pair $(P,Q)$ for the canonical map of $\Surf_{\tau} = div(\theta)$ can be represented by a couple of points $(\widetilde{P}, \widetilde{Q})$ in $\left(\SA_{\tau, 12} \times \E_3 \right)^2$ according to Diagram \ref{restrictionPdiagram2}, with:
\begin{equation*}
\begin{split}
\widetilde{P}&=(\xc, \zc) \\
\widetilde{Q}&=(\yc, \wc) \PUNKT
\end{split}
\end{equation*}
The section $\theta$ can be written, as in Equation \ref{thetaDECOMPOSITIONij}, as:
\begin{equation} \label{eqStau}
  \theta = \cf_{\tau, 12}\cdot \theta^{\E_3}_0 + \cg_{\tau, 12}\cdot \theta^{\E_3}_1 
\end{equation}
To simplify the notations, the section $\cf_{\tau, 12}$ (resp. $\cg_{\tau, 12}$) in Equation \ref{eqStau} will be denoted by $\cf$ (resp. $\cg$), and the section $\theta^{\E_3}_0$ (resp. $\theta^{\E_3}_1$) by $\theta_0$ (resp. $\theta_1$). 
The condition that the canonical image of $P$ and $Q$ coincide implies that
there exists $\lambda \in \C^*$ such that:
\begin{equation*}
\begin{cases}
\cf(\xc)\theta_0(\zc) &= \lambda \cdot \cf(\yc)\theta_0(\wc) \\
\cg(\xc)\theta_1(\zc) &= \lambda \cdot \cg(\yc)\theta_1(\wc) \\
\cf(\xc)\DTH{0}(\zc) + \cg(\xc)\DTH{1}(\wc) &= \lambda \cdot \Bigl( \cf(\yc)\DTH{0}(\wc) + \cg(\yc)\DTH{1}(\wc) \Bigr) \PUNKT
\end{cases}
\end{equation*}
(Here $\DTH{0}$ and $\DTH{1}$ are the derivatives of $\theta_0$ and $\theta_1$). This leads to
\begin{align*}
& \cf(\xc)\left(\frac{\DTH{0}(\zc)\theta_0(\wc) - \theta_0(\zc)\DTH{0}(\wc)}{\theta_0(\wc)}  \right) \ \ + \\ & \cg(\xc)\left(\frac{\DTH{1}(\zc)\theta_1(\wc) - \theta_1(\zc)\DTH{1}(\wc)}{\theta_1(\wc)}  \right) \ \ =  0 \BEISTRICH
\end{align*}
and we can conclude that:
\begin{align} \label{definitionofpsi0}
\cf(\xc)\theta_1(\wc) \rhos{0}{3}(\zc,\wc) + \cg(\xc)\theta_0(\wc)\rhos{1}{3}(\zc,\wc)   &=  0 \BEISTRICH
\end{align}
where, with $j = 0,1$:
\begin{align} \label{RHODEF}
\rhos{j}{3}(\zc,\wc) := \begin{vmatrix}\theta_j(\zc) & \theta_j(\wc) \\ \DTH{j}(\zc) & \DTH{j}(\wc) \end{vmatrix} \in H^0(E_3 \times E_3, \OC_{E_3}(2O_{E_3}) \boxtimes  \OC_{E_3}(2O_{E_3}))
\end{align}
Together with the defining equation of $\Surf_{\tau}$ in (\ref{eqStau}), we obtain the following linear system in $\cf(\xc)$ and $\cg(\xc)$:
\begin{align*} \label{systemfgrho}
&\cf(\xc)\theta_0(\zc) & &+& &\cg(\xc)\theta_1(\zc) &  &=  0 \\
&\cf(\xc)\theta_1(\wc) \rhos{0}{3}(\zc,\wc) & &+& &\cg(\xc)\theta_0(\wc)\rhos{1}{3}(\zc,\wc) &   &=  0 \PUNKT
\end{align*}
The determinant of the matrix of the linear system above is precisely
\begin{equation}\label{DEFVAB}
\secV_3(\zc,\wc) := \begin{vmatrix}\theta_0(\zc)\theta_{0}(\wc) & \theta_1(\zc)\theta_1(\wc) \\
						\rhos{0}{3}(\zc,\wc) & \rhos{1}{3}(\zc,\wc)\end{vmatrix} \in H^0(E_3 \times E_3, \OC_{E_3}(4O_{E_3}) \boxtimes  \OC_{E_3}(4O_{E_3})) \BEISTRICH
\end{equation}
and we conclude that $\xc$ belongs to $\left(div(\cf) \cap div(\cg)\right) = \cfc_{\tau, 12} \cap \cgc_{\tau, 12}$, or
\begin{equation*}
\secV_3(\zc, \wc) = 0
\end{equation*}
This proves the claim of the lemma.
\end{proof}

\begin{rmk} \label{OSScomponentZ}
In general, if we have an elliptic curve $\E$ and $\theta_0, \theta_1$ the canonical basis for $H^0(\E, \mathcal{O}_{\E}(2 \cdot O_{\E})$ it is easy to determine a polynomial equation, in a suitable affine open set, of the divisor $V := div(\secV)$ in $\E \times \E$ defined as in Equation \ref{DEFVAB} above as 
\begin{equation}\label{DEFVAB}
\secV(\zc,\wc) := \begin{vmatrix}\theta_0(\zc)\theta_{0}(\wc) & \theta_1(\zc)\theta_1(\wc) \\
						\rhos{0}{}(\zc,\wc) & \rhos{1}{}(\zc,\wc)\end{vmatrix} \PUNKT
\end{equation}
(Here $\rhos{0}{}$ and $\rhos{1}{}$ are defined as in Equation \ref{RHODEF}.)

In the notations of Example \ref{basicexampleEll}, for every point $(\zc,\wc)$ of $\E \times \E$ it easily seen that:
\begin{equation*}
\secV(\zc,\wc) = \theta_0(z)\theta_1(z)\cdot \left(\theta_0(\wc)\DTH{1}(\wc) - \DTH{0}(\wc)\theta_1(\wc) \right) - \theta_0(w)\theta_1(w)\cdot \left(\theta_0(\zc)\DTH{1}(\zc) - \DTH{0}(\zc)\theta_1(\zc) \right)
\end{equation*}
In the affine open set $U$ where $\theta_0(z)$ and $\theta_0(w)$ do not vanish, we can divide by $\theta_0(z)^2\theta_0(w)^2$ and in local affine coordinates $(\x{1}, \y{1}, \x{2}, \y{2})$ with two equations $g_1$, $g_2$ in Legendre normal form as in Equation \ref{LegendreNormalForm} in Example \ref{basicexampleEll}, we get that:
\begin{equation*}
\secV(\zc,\wc) \sim \x{1}\y{2} - \x{2}\y{1}
\end{equation*}
Moreover, by using the equations $\y{j}^2 = (\x{j}^2 - 1)(\x{j}^2 - \pa{}^2)$ of $\E$ on the components, we can easily recover the equations of the components of $V$ in the open set $U$: 
\begin{equation*}
div(V)|_U = (\x{1}^2 - \x{2}^2)(\x{1}^2\x{2}^2 - \pa{}^2)
\end{equation*}
In conclusion, the divisor $V$ splits into the sum of components:
\begin{equation*}
V = \Delta + (-1).\Delta + X + (-1).X \BEISTRICH
\end{equation*}
where $\Delta$ is the diagonal locus in $\E \times \E$. All components are invariant under the action of the involution $(\zc,\wc) \mapsto (\zc + 1, \wc+ 1)$, and the irreducible component $X$ meets $\Delta$ in $8$ points whose coordinates are not $2$-torsion on $\E$.
\end{rmk}

\begin{proof}[of Proposition \ref{generalbehavior2}]
Let us consider $(A_{\tau_0}, \mathcal{L}_{\tau_0})$ the $(1,2,2)$-polarized abelian variety associated to a diagonal matrix $\tau$ in $\mathcal{H}_1^3$, with a divisor $\Surf_{\tau_0} = div(\theta)$ which is sufficiently general in the polarization $|\mathcal{L}_{\tau_0}|$, in accordance with the conditions of Remark \ref{BEQUADROG}. 

Since these conditions are open, they still hold for a sufficiently small neighborhood $\mathcal{U}$ of $\tau_0$ in the closed locus $\HD^{(12)}$. Hence, for a suitable $\mathcal{U}$, there is a family $(A_{\mathcal{U}}, \mathcal{L}_{\mathcal{U}})$ of $(1,2,2)$-polarized abelian threefolds with $\Surf_{\mathcal{U}} = div(\theta)$, such that, considered the decomposition
\begin{equation*} 
  \theta = \cf_{\mathcal{U}}\cdot \theta_0 + \cg_{\mathcal{U}}\cdot \theta_1 \BEISTRICH
\end{equation*}
the fibers of the induced families $\cfc_{\mathcal{U}} =div \left(\cf_{\mathcal{U}}\right)$ and $\cgc_{\mathcal{U}} = div\left(\cg_{\mathcal{U}}\right)$ defined as in Notation \ref{NOTNUXI} (for the sake of readability we omit the indices $1,2$ as in the proof of Lemma \ref{secondOssBXE}), are smooth genus $3$ non-hyperelliptic curves. \newline 

By Proposition \ref{generalbehavior1}, the possibile critical pairs on the central fiber $\Surf_{\tau_0}$ are only of one of the following types: \newline
 $\bullet$ \textbf{type (a)}: The pair belongs to the divisor $W_j = div \left(\frac{\partial \theta}{\partial z_j}\right)$ for some $j$. In this case, the points are conjugated under the the sign reversing involution $\iota_j$, which changes the sign on the $j$-th coordinate $z_j \mapsto -z_j$. \newline
 $\bullet$ \textbf{type (b)}: The pair belongs to the intersection of two divisors $W_i$ and and $W_j$ for some couple of indices $i,j$. In this case, the points are conjugated under the involution $\iota_i \iota_j$. \newline

\vspace{0.25cm}
Hence, to prove the proposition, it is enough to prove the following claim. \newline 

\noindent \textbf{Claim}: Let us consider a family of critical pairs $(\Surf_{\mathcal{U}} \longrightarrow \mathcal{U}, \DEFP_{12})$ on $\mathcal{U}$ (up to shrinking $\mathcal{U}$). Then the family $\DEFP_{12}$ is one of the following types:
\begin{itemize}
 \item[1)] On every fiber (up to shrinking $\mathcal{U}$), the critical pair is of type (a) with $j=3$. In particular, the pair always belongs to the divisor $W_3$.
 \item[2)] On every fiber (up to shrinking $\mathcal{U}$), the critical pair is of type (b) with $(i,j)=(1,2$. In particular, the pair always belongs to the divisor $W_1 \cap W_2$.
\end{itemize}

\vspace{0.25cm} 
To this purpose, let us pick a fiber on a point $\tau$ of $\mathcal{U}$ and a critical pair $(P_{\tau},Q_{\tau})$ on $\Surf_{\tau}$. This critical pair can be represented (see Diagram \ref{restrictionPdiagram2}) as a couple of points $(\widetilde{P}_{\tau}, \widetilde{Q}_{\tau})$ in $\left(\SA_{\tau, 12} \times \E_3 \right)^2$ with
\begin{equation*}
\begin{split}
\widetilde{P}_{\tau}&=(\xc, \zc) \\
\widetilde{Q}_{\tau}&=(\yc, \wc) \BEISTRICH
\end{split}
\end{equation*}
where the coordinates $\widetilde{P}_{\tau}$ and $\widetilde{Q}_{\tau}$ are functions of $\tau$.
By Lemma \ref{secondOssBXE}, we have that $(\xc,\yc)$ belongs to $\cfc_{\tau} \cap \cgc_{\tau}$, or $(\zc, \wc)$ belongs to $ div(\secV_3)$.

In the first case, already the condition that $(\xc,\yc)$ belongs to $\cfc_{\tau}$ implies that $x = y$, since the curve $\cfc_{\tau}$ is non-hyperelliptic for every $\tau$, up to shrinking $\mathcal{U}$: indeed, recall that (again by Proposition \ref{desccanonicalmap}) the canonical map of $\cfc_{\tau}$ is obtained by considering the derivatives $\frac{\partial \cg}{\partial z_1}$ and $\frac{\partial \cg}{\partial z_2}$ together with the restriction of the sections of $H^0(\SA, \PMm)$ to $\cgc$, and we can conclude that $\phi_{\cgc_{\tau}}(x) = \phi_{\cgc_{\tau}}(y)$, where $\phi_{\cgc_{\tau}}$ denotes the canonical map of $\cgc_{\tau}$. \newline

\noindentTX Let us assume that $(\zc, \wc)$ belongs to $V = div(\secV_3)$.

\noindentTX The affine equations in Remark \ref{OSScomponentZ} show that $V$ splits into a sum of irreducible components:
\begin{equation*}
V = \Delta + (-1).\Delta + X + (-1).X \BEISTRICH
\end{equation*}

\noindentTX For the central fiber $\tau_0$, we described these components in Remark \ref{OSScomponentZ} and we determined their equation on an affine open set. Since the components $(\zc, \wc)$ of the points $\widetilde{P}_{\tau}$ and $\widetilde{Q}_{\tau}$ are bound to move on $V$, we have to deal with two cases:
\begin{itemize}
 \item $(\zc, \wc)$ is contained in $\Delta \cup (-1).\Delta$ but not in $X \cup (-1).X'$.
 \item $(\zc, \wc)$ is contained in $X \cup (-1).X'$. \newline
\end{itemize}

\noindentTX Let us assume that $(\zc, \wc)$ is contained in $\Delta \cup (-1).\Delta$ but not in $X \cup (-1).X'$. Then, because this latter condition is open, we can conclude that the same condition holds true on the neighborhood $\mathcal{U}$ up to shrinking it. Hence, on every fiber on $\tau$, we have that for the coordinates $(\zc, \wc) = (\zc, \pm\zc)$ or $(\zc, \wc) = (\zc, \pm\zc + 1)$. Let us assume that we are in the case $(\zc, \wc) = (\zc, \zc)$ (the other cases can be treated similarly). The ratio of the values $[s,t]:=[\theta_0(z), \theta_1(z)]$ of the $(2)$-polarization on $\E_3$ define, according to the equation of $\Surf_{\tau}$, a curve
\begin{equation} \label{EQFCICLE}
\mathcal{F} \colon s \cdot \cf  + t \cdot \cg = 0
\end{equation}
in the $(2,2)$-polarization class of $\B_{\tau}$,
on which the points $x$ and $y$ belong.

On the other side, the expression of the canonical map of $\Surf_{\tau}$ is
\begin{equation}\label{CANEXP2}
 \CAN{\Surf_{\tau}} = \left[H^0(\SA_{\tau}, \PMp)\theta_0, H^0(\SA_{\tau}, \PMm)\theta_1, \frac{\partial \cf}{\partial z_1} \theta_0 + \frac{\partial \cg}{\partial z_1} \theta_1, \frac{\partial \cf}{\partial z_2} \theta_0 + \frac{\partial \cg}{\partial z_2} \theta_1,\cf \DTH{0} + \cg \DTH{1} \right] \PUNKT
\end{equation}
This expression can be easily compared with the one of the canonical map of $\mathcal{F}$, which is
\begin{equation}\label{CANEXP3}
 \CAN{\mathcal{F}} = \left[s\cdot H^0(\SA_{\tau}, \PMp), t \cdot H^0(\SA_{\tau}, \PMm), s\frac{\partial \cf}{\partial z_1} + t\frac{\partial \cg}{\partial z_1} , s\frac{\partial \cf}{\partial z_2} + t\frac{\partial \cg}{\partial z_2}\right] \BEISTRICH
\end{equation}
and follows that $\CAN{\mathcal{F}}(x) = \CAN{\mathcal{F}}(y)$ . 

\noindentTX If $\mathcal{F}$ is smooth, then it is non-hyperelliptic (because in a simple $(2,2)$-polarization there are no irreducible hyperelliptic curves). Hence, in this case, we conclude that $x=y$ and $P=Q$. \newline

\noindentTX If otherwise $\mathcal{F}$ is not smooth, then we can assume that $\mathcal{F}$ is a nodal curve. Indeed, since the $(1,2)$-polarization class on a simple abelian surface is a pencil only containing irreducible curves (cf. Remark \ref{osservationsingularfibers}), we may assume, by generality, that the singular curves in the pencil generated by $\cf$ and $\cg$ in the $(2,2)$-polarization class of $\B$ only contains nodal curves, which also are non-hyperelliptic since they are irreducible. 

\noindentTX This implies that $x=y$ if $x$ is not a node on $\mathcal{F}$, otherwise $x =-y$, and the second possibility occurs precisely when $(P,Q)$ is a critical pair which belong to $W_1 \cap W_2$ with $Q = \iota_1 \iota_2 (P)$.

\vspace{0.25cm}






\noindentTX Let us suppose now that we are in the second case, namely the case that $(\zc, \wc)$ is contained in $X \cup (-1).X$ and $\wc \neq \pm\zc$ when $\tau$ varies away from the central point $\tau_0$ in a suitable base neighborhood. Our goal is to show by contradiction that this case actually does not occur. \newline

\noindentTX By Proposition \ref{generalbehavior1}, on the central fiber $\tau_0$ the point $(\zc, \wc)$ must also belong on $\Delta \cup (-1).\Delta$. Hence, without loss of generality (since the other cases are similar), we may assume by contradiction that on the central fiber $(\zc, \wc)$ belongs to the intersection $\Delta \cap X$, with $\zc$ and $\wc$ not $2$-torsion on $\E_3$, and on them neither $\theta_0$ nor $\theta_1$ vanish. (cf. Remark \ref{OSScomponentZ}). In this case, we can also assume that, on $\tau_0$,
\begin{equation*}
\yc = (\yc_1, \yc_2) = (\xc_1, -\xc_2) = \inv_2(\xc)
\end{equation*}
and both $P$ and $Q$ belong to $W_2$. Since all the sections of the canonical bundle of $\Surf_{\tau}$ are invariant with the only exception of $\DER{z_2}{\theta}$, we have that all the family of critical pairs $(P,Q)_{\tau}$ also belongs to the family of divisors $W_2$ when $\tau$ varies in a sufficiently small neighborhood of $\tau_0$, with
\begin{equation}\label{KEYCONDITIONsign}
\xc \neq \pm\yc
\end{equation}
for every $\tau$ in the base neighborhood. \newline
However, $x$ and $y$ are supposed here not to remain costant with respect to $\tau$, and in particular they move outside of the preimage in $\B$ of the union of the base loci $\mathcal{B}(\mathcal{M}^+) \cup \mathcal{B}(\mathcal{M}^-)$ in $\SA$. By the expression of the canonical map in (\ref{CANEXP2}), it follows that:
\begin{equation} \label{DERF}
 \phi_{|\mathcal{M}^+| \times |\mathcal{M}^-|}(x) = \phi_{|\mathcal{M}^+| \times |\mathcal{M}^-|}(y)
\end{equation}
where $\phi_{|\mathcal{M}^+| \times |\mathcal{M}^-|}$ denotes the rational map defined on $\SA$ in $\PR^1 \times \PR^1$ from the pencils $|\mathcal{M}^+|$ and $|\mathcal{M}^-|$ on each factor.

On the other side, the family $(P,Q)_{\tau}$ moves inside the divisor $W_2$. Thus $x$ and $y$ must be contained, for every $\tau$, in the divisor defined by the equation:

\begin{equation*}
\cf \DER{z_2}{\cg} \cdot \DER{z_2}{\cf}\cg = \SKMULT{z_2}{}(\cf \wedge \cg) = 0
\end{equation*}
However, by Proposition \ref{TWOFOURTH}, for the general choice of $\cf$ and $\cg$, a couple $(x,y)$ of points in the zero locus of $\SKMULT{z_2}{}(\cf \wedge \cg)$ which also fulfills the condition in (\ref{DERF}) must be of the form $(x, x)$ or $(x, -x)$. But this contradicts Equation (\ref{KEYCONDITIONsign}) above, and the claim is proved.

\end{proof}


In conclusion, we can prove the main result of this paper: 
\begin{thm1} \label{teoremafinale}
Let be $(A, \mathcal{L})$ a general $(1,2,2)$-polarized abelian threefold and let be $\Surf$ a general surface in the linear system $|\mathcal{L}|$. Then the canonical map of $\Surf$ is a holomorphic embedding. 
\end{thm1}
\begin{proof}
   
\noindentTX By Proposition \ref{injectivedifferential}, the canonical map is a local embedding in the general case. Hence, it sufficies to prove that the canonical map is injective in the general case. \newline

\noindentTX By way of contradiction, let us assume that the latter claim is false. \newline

\noindent Our assumption implies that, for every $(1,2,2)$-polarized abelian threefold $A$ and every smooth surface $\Surf$ within its polarization, the canonical map of $\Surf$ is not injective. Following Notation \ref{productSurfEll}, our assumption will also hold true for all abelian threefolds in a neighborhood $\mathcal{U}$ in the moduli space around a certain fixed threefold which is isogenous to a product of elliptic curves. Let us consider a $(1,2,2)$-polarized abelian threefold $A_{\tau_0}$, where $\tau_0 = (\tau_{jj})_{j=1,2,3}$ is a diagonal matrix in $\SiegH_1^3$. We have then a isogeny:
\begin{equation*} 
\prTA \colon T \longrightarrow A_{\tau_0}:= \bigslant{T_{\tau_0}}{\left<e_1+e_2+e_3\right>} \PUNKT
\end{equation*}
where $T: = \E_1 \times \E_2 \times \E_3$ is the product of elliptic curves, each defined as a quotient $\E_j = \bigslant{\C}{\tau_{jj}\Z \oplus 2\Z}$ carrying a natural $(2)$-polarization.

According to Proposition \ref{generalbehavior1}, on a general surface $\Surf_{\tau_{0}} = div(\theta)$ in the $(1,2,2)$-polarization class of $A_{\tau_0}$, there exists a critical pair $(P,Q)$ on $\Surf_{\tau_{0}}$ which, by our assumption, admits the existence of a family of critical pair
$(\Surf_{\mathcal{U}} \longrightarrow \mathcal{U}, \DEFP)$. According to Definition \ref{generaldeformationdefinition}:
\begin{itemize}
\item[i)] $\mathcal{U}$ is an open set of $\HD$ of matrices whose diagonal is fixed and equal to the entries of $\tau_0$
\item[ii)] $\Surf_{\mathcal{U}} \longrightarrow \mathcal{U}$ is a family of surfaces defined by $\Surf_{\tau} = div(\theta)$ for every $\tau$.
\item[iii)] $\DEFP$ is a closed irreducible subscheme of 
\begin{equation*} 
\mathcal{K}_{\mathcal{U}} := (\phi_{\Surf_{\mathcal{U}}} \times_{\mathcal{U}} \phi_{\Surf_{\mathcal{U}}})^{-1}(\Delta_{\mathcal{U}}) \subseteq \Surf_{\mathcal{U}} \times_{\mathcal{U}} \Surf_{\mathcal{U}}
 \end{equation*}
which is dominant over $\mathcal{U}$ and such that its restriction on the central fiber $\tau_0$ coincides with $(P,Q)$. \newline
\end{itemize}

\noindentTX The geometric points of the intersection of $\DEFP$ with the diagonal locus in $\Surf_{\mathcal{U}} \times_{\mathcal{U}} \Surf_{\mathcal{U}}$ represent critical pairs of infinitely near couples of points on a certain surface of the family. 
If $P = Q$, without loss of generality we can assume, by Proposition \ref{generalbehavior1}, that the third component of $P$ is represented by a $2$-torsion point of $\E_3[2]$, and $P$ belongs to the divisor $W_{3} = div \left(\DER{z_3}{\theta}\right)$. However, by Proposition \ref{generalbehavior2},  the restriction of the family of critical pairs $\DEFP$ to the close locus $\mathcal{U}_{13}$ defines a subfamily of critical pairs whose general element $(P_{\tau}, Q_{\tau})$ belongs to $W_{2}$ with $Q_{\tau} = \iota_2(P_{\tau})$, or to the intersection $W_{1} \cap W_{3}$ with $Q_{\tau} = \iota_1\iota_3(P_{\tau})$ for every $\tau$ in $\mathcal{U}_{13}$ (recall that $\iota_{h}: \C^3 \longrightarrow \C^3$ denotes the sign-changing involution of to the $h$-th coordinate, and by definition in Equation \ref{LOCIEQ} and Diagram \ref{restrictionPdiagram}, $\mathcal{U}_{ij} = \mathcal{U} \cap \HD^{(ij)}$).

In both cases we can conclude that $P_{\tau_0}$ has $2$ coordinates which are $2$-torsion on the respective elliptic curve factors of $T$, and we would reach a contradiction with the generality conditions in Remark \ref{BEQUADROG}. \newline

From now on, we can assume (up to shrinking $\mathcal{U}$) that $\DEFP$ does not intersect the diagonal locus in $\Surf_{\mathcal{U}} \times_{\mathcal{U}} \Surf_{\mathcal{U}}$. For every couple of indices $(ij)$ we can fix an irreducible component $\DEFP_{ij}$ of the restriction $\mathcal{P}$ to the closed locus $\mathcal{U}_{ij}$. By our assumption on $\DEFP$, the component $\DEFP_{ij}$ also has empty intersection with the diagonal subscheme of $\Surf_{\mathcal{U}_{ij}} \times_{\mathcal{U}_{ij}} \Surf_{\mathcal{U}_{ij}}$. 


\noindentTX By applying Proposition \ref{generalbehavior2}, we have that every fiber $\DEFP_{\tau, ij}$ of the family $\DEFP_{ij}$ must be contained in 
 \begin{align*} 
\mathcal{X}_{\tau,ij}  &:= \SET{(P, \iota_{i}\iota_{j}(P)) \in \Surf_{\tau} \times \Surf_{\tau}}{P \in W_i \cap W_j}  \\
\mathcal{W}_{\tau,k}  &:= \SET{(P, \iota_{k}(P)) \in \Surf_{\tau} \times \Surf_{\tau}}{P \in W_k} \PUNKT
 \end{align*}
By definition, we have that $\mathcal{U}_{12} \cap \mathcal{U}_{13} \cap \mathcal{U}_{23} = \{\tau_0\}$, hence on the central fiber we must have that $(P,Q)$ belongs to the intersections of of the different $\mathcal{X}_{\tau_0, ij}$ or to $\mathcal{Y}_{\tau_0, k}$, in accordance with the behavior of $\DEFP_{ij}$. 

\noindentTX Because all surfaces of the family $\Surf_{\mathcal{U}}$ are supposed to be smooth, then $(P,Q)$ belongs to an intersection of the form $\mathcal{X}_{\tau_0, ij} \cap \mathcal{W}_{\tau_0, i}  \cap \mathcal{W}_{\tau_0, j}$ for a couple of indices $i,j$, with $i \neq j$.
 On the other hand, every point of $\mathcal{X}_{\tau_0, ij} \cap \mathcal{W}_{\tau_0, i}  \cap \mathcal{W}_{\tau_0, j}$ is of the form $(P,Q)$ such that $Q = \iota_{i}\iota_{j}(P)$, $Q = \iota_{i}(P)$ and $Q = \iota_{j}(P)$. But this implies that $P=Q$, and we reach a contradiction with the assumption that $\DEFP$ does not intersect the diagonal sublocus.
 
\noindentTX The proof of the theorem is complete.
 \end{proof}


\begin{acknowledgements}\label{ackref}
 The author gratefully acknowledges Prof. Fabrizio Catanese for the numerous useful discussions 
 regarding this research topic. My grateful thanks are also extended to Dr. Stephen Coughlan and Dr. Davide
 Frapporti for their valuable suggestions.
\end{acknowledgements}


\newcommand{\BIBND}[4]{\bibitem{#1} 
{\bibname #2}, #3.
#4.}
\newcommand{\BIBA}[5]{\bibitem{#1}
{\bibname #2}, #3. 
#4, (#5).}
\newcommand{\BIBB}[6]{\bibitem{#1}
{\bibname #2}, #3. 
#4, (#5), #6.}


%
%

%
%
%
%


\affiliationone{}
%
\end{document}